\documentclass[10pt]{article}

\usepackage[dvips]{epsfig} %this package is for manipulating postscript graphics
\usepackage{amssymb,amsmath,amsthm,mathrsfs} %this package contains all the special ams fonts etc.
\usepackage{graphicx}
\usepackage{lineno}

\usepackage[authoryear,sort&compress]{natbib}
\usepackage{url}
\usepackage{enumerate}
\usepackage{colordvi,color,pspicture}
\usepackage{psfrag}
\DeclareMathOperator{\argsup}{argsup}

\DeclareMathOperator{\argmin}{argmin}

\DeclareMathOperator{\BER}{BER}

\newcommand{\E}[1]{\mathbf{E}\,#1}

\newcommand{\Var}[1]{\mathbf{Var}\,#1}
\newcommand{\Cov}[1]{\mathbf{Cov}\,#1}
\newcommand{\I}{\mathbf{I}}

\newcommand{\C}{\mathcal{C}}
\newcommand{\D}{\mathcal{D}}
\newcommand{\U}{\mathcal{U}}

\newcommand{\g}{\gamma}
\newcommand{\G}{\Gamma}
\newcommand{\V}{\mathcal{V}}
\newcommand{\A}{\mathcal{A}}
\newcommand{\Y}{\mathcal{Y}}
\newcommand{\y}{\mathsf{y}}
\newcommand{\X}{\mathcal{X}}
\newcommand{\M}{\mathcal{M}}

\newcommand{\mS}{\mathcal{S}}

\newcommand{\NY}{N_{\Y}}
\newcommand{\NS}{N_S}
\newcommand{\NAS}{N_{AS}}
\newcommand{\NCS}{N_{CS}}
\newcommand{\NCSt}{N_{CS}^{\tau}}
\newcommand{\NPE}{N_{PE}}
\newcommand{\TY}{T\left( \Y_3 \right)}
\newcommand{\N}{\mathcal{N}}
\newcommand{\R}{\mathbb{R}}
\newcommand{\T}{\mathcal{T}}

\newcommand{\RS}{\mathscr{R}_S}
\newcommand{\UT}{\U(\TY)}
\newcommand{\Tr}{\mathscr T^r}
\newcommand{\ve}{\varepsilon}

\theoremstyle{plain}
\newtheorem{theorem}{Theorem}[section]
\newtheorem{lemma}[theorem]{Lemma}

\newtheorem{corollary}[theorem]{Corollary}

\theoremstyle{definition}
\newtheorem{definition}[theorem]{Definition}

\theoremstyle{remark}

\newtheorem{remark}[theorem]{Remark}

\begin{document}

%\begin{singlespace}

\title{Technical Report \# KU-EC-09-2:\\
Extension of One-Dimensional Proximity Regions to Higher Dimensions}
\author{
Elvan Ceyhan
\thanks{Address:
Department of Mathematics, Ko\c{c} University, 34450 Sar{\i}yer, Istanbul, Turkey.
e-mail: elceyhan@ku.edu.tr, tel:+90 (212) 338-1845, fax: +90 (212) 338-1559.
}
}

\date{\today}
\maketitle

%\linenumbers
%\modulolinenumbers[2]

\begin{abstract}
\noindent
Proximity maps and regions are defined based on the relative allocation of points from two or
more classes in an area of interest and are used to construct random graphs
called proximity catch digraphs (PCDs) which
have applications in various fields.
The simplest of such maps is the spherical proximity map which maps a point from
the class of interest to a disk centered at the same point with radius being the
distance to the closest point from the other class in the region.
The spherical proximity map gave rise to class cover catch digraph (CCCD) which was
applied to pattern classification.
Furthermore for uniform data on the real line,
the exact and asymptotic distribution of the domination number
of CCCDs were analytically available.
In this article, we determine some appealing properties of the spherical proximity
map in compact intervals on the real line and use these properties as a
guideline for defining new proximity maps in higher dimensions.
Delaunay triangulation is used to partition the region of interest
in higher dimensions.
Furthermore, we introduce the auxiliary tools used for the construction
of the new proximity maps, as well as some related concepts that will be
used in the investigation and comparison of them and the resulting graphs.
We characterize the geometry invariance of PCDs for uniform data.
We also provide some newly defined proximity maps in
higher dimensions as illustrative examples.
\end{abstract}

\noindent
%{\scriptsize

%\vspace{.25 in}

%\indent
%$^\star$ This research was supported by the research agency TUBITAK via the Kariyer Project \# 107T647.\\
%\indent
%$^*$corresponding author.\\
%\indent {\it e-mail:} elceyhan@ku.edu.tr (E.~Ceyhan) }

%\end{singlespace}

\newpage

%\linenumbers
%\modulolinenumbers[2]

\section{Introduction}
\label{sec:introduction}
Classification and clustering have received
considerable attention in the statistical literature.
In recent years, a new classification approach has been developed.
This approach is
based on the proximity maps that incorporate the
relative positions of the data points from various classes.
Proximity maps and the associated (di)graphs are used in disciplines
where shape and structure are crucial.
Examples include computer vision (dot patterns), image
analysis, pattern recognition (prototype selection), geography and
cartography, visual perception, biology, etc.
\emph{Proximity graphs} were first introduced by \cite{toussaint:1980},
who called them \emph{relative neighborhood graphs}.
The notion of relative neighborhood graph has been generalized in
several directions and all of these graphs are now called proximity graphs.
From a mathematical and algorithmic point of view,
proximity graphs fall under the category of \emph{computational geometry}.

A general definition of proximity graphs is as follows:
%\begin{definition}
Let $V$ be any finite or infinite set of points in $\R^d$. Each
(unordered) pair of points $(p,q) \in V \times V$ is associated with
a neighborhood $\mathfrak N(p,q) \subseteq \R^d$.
Let $\mathfrak P$ be a property defined on $\mathfrak N=\{\mathfrak N(p,q):\; (p,q) \in V \times V\}$.
A \emph{proximity} (or \emph{neighborhood}) \emph{graph}
$G_{\mathfrak N,\mathfrak P}(V,E)$ defined by the property
$\mathfrak P$ is a graph with the vertex set $V$ and the edge set
$E$ such that $(p,q) \in E$ iff $\mathfrak N(p,q)$ satisfies property $\mathfrak P$.
%$\square$ %\end{definition}
Examples of most commonly used proximity graphs are the Delaunay tessellation,
the boundary of the convex hull, the Gabriel graph,
relative neighborhood graph, Euclidean minimum spanning tree, and
sphere of influence graph of a finite data set.
See, e.g., \cite{jaromczyk:1992} for more detail.
The \emph{Delaunay tessellation} $\D(V)$ of a finite set of points $V$,
is the dual of the Voronoi diagram generated by $V$.
See \cite{okabe:2000} for further details.
The \emph{convex hull} of a set $V \subset \R^d$, denoted as $\C_H(V)$,
%encircles $V$ as if by a rubber band so that the region inside is $\C_H(V)$. More formally,  $\C_H(V)$
is the intersection of all convex sets (there exists infinitely many of them) that contain $V$.
The boundary of $\C_H(V)$ can be viewed as a proximity graph which is
also a subgraph of $\D(V)$.
The \emph{Gabriel graph} of $V$, denoted as $GG(V)$, is defined as the
graph in which $(p,q)$ is an edge of $GG(V)$ iff the circle centered
at the midpoint of the line segment $\overline{pq}$ and with
diameter $d(p,q)$ (the distance between $p$ and $q$)
does not contain any other points from $V$.
The \emph{relative neighborhood graph} of $V$
 is a prominent representative of the family of graphs which are
defined by some sort of neighborliness. For a set of points  $V
\subset \R^d$, the relative neighborhood graph of $V$, denoted
$RNG(V)$, is a graph with vertex set $V$ and edge set which are
exactly the pairs $(p,q)$ of points for which $d(p,q) \le \min_{v
\in V} \max(d(p,v),d(q,v))$.  That is, $(p,q)$ is an edge of
$RNG(V)$ iff $\text{Lune}(p,q)$ does not contain any other points of
$V$, where $\text{Lune}(p,q)$ is defined as the intersection of two
discs centered at $p,q$ each with radius $d(p,q)$ (see, e.g.,
\cite{jaromczyk:1992}).
The \emph{Euclidean minimum spanning tree} of $V$, denoted
$EMST(V)$, is defined as the spanning tree in which the sum of the
Euclidean lengths of the edges yield the minimum over all spanning
trees with vertex set $V$.
The \emph{sphere of influence graph} on $V$, denoted as $SIG(V)$, has
vertex set $V$ and $(p,q)$ as an edge iff the circles centered at
$p$ and $q$ with radii $\min_{v \in V\setminus\{p\}}d(p,v)$ and
$\min_{v \in V\setminus\{q\}}d(q,v)$, respectively, have nonempty
intersection.
Note that $EMST(V)$ is a subgraph of $RNG(V)$ which is a subgraph of
$GG(V)$ which is a subgraph of $\D(V)$ (see \cite{okabe:2000}).
Furthermore, in the examples above, $d(x,y)$, can be any distance in $\R^d$.
Furthermore, the distance between a point $x$ and a set $A$ is
defined as $d(x,A):=\inf_{y \in A} d(x,y)$;
and the distance between two sets $A$ and $B$ is defined as
$d(A,B):=\inf_{(x,y) \in A\times B} d(x,y)$.

A \emph{digraph} is a directed graph, i.e., a graph with directed
edges from one vertex to another based on a binary relation.
Then the pair $(p,q) \in V \times V$ is an ordered pair and
is an \emph{arc} (directed edge) denoted as $pq$ to reflect
the difference between an arc and an edge.
For example, the nearest neighbor (di)graph in
\cite{paterson:1992} is a proximity digraph.
The \emph{nearest neighbor digraph}, denoted as $NND(V)$, has the vertex set $V$
and $pq$ as an arc iff $d(p,q) = \min_{v \in V \setminus\{p\}} d(p,v)$.
That is, $pq$ is an arc of $NND(V)$ iff $q$ is a nearest
neighbor of $p$. Note that if $pq$ is an arc in $NND(V)$,
then $(p,q)$ is an edge in $RNG(V)$.

The \emph{proximity catch digraphs} (PCDs) are based on the property $\mathfrak P$
that is determined by the following mapping which is defined in a
more general space than $\R^d$.
%\begin{definition}
Let $(\Omega,\M)$ be a measurable space. The \emph{proximity map}
$N(\cdot)$ is given by  $N:\Omega \rightarrow \wp(\Omega)$, where
$\wp(\cdot)$ is the power set functional, and the \emph{proximity
region} of $x \in \Omega$, denoted as $N(x)$, is the image of $ x \in
\Omega$ under $N(\cdot)$. %$\square$ \end{definition}
The points in $N(x)$ are thought of as being ``closer" to $x \in \Omega$
than are the points in $\Omega \setminus N(x)$.
Proximity maps are
the building blocks of the \emph{proximity graphs} of \cite{toussaint:1980};
an extensive survey is available by \cite{jaromczyk:1992}.
%\begin{definition}
The PCD $D$ has the vertex set
$\V=\bigl\{ p_1,\ldots,p_n \bigr\}$ and the arc set $\A$ is defined
as $p_ip_j \in \A$ iff $p_j \in N(p_i)$ for $i\not=j$. %$\square$ \end{definition}
Notice that the PCD $D$ depends on the \emph{proximity} map $N(\cdot)$,
and if $p_j \in N(p_i)$, then $N(p_i)$ is said to \emph{catch} $p_j$.
Hence the name \emph{proximity catch digraph}.
If arcs of the form $p_ip_i$ (i.e., loops) were allowed, $D$ would
have been called a \emph{pseudodigraph} according to some authors
(see, e.g.,  \cite{chartrand:1996}).

The finite and asymptotic distribution of the domination number of CCCDs
for uniform data in $\R$ are mathematically tractable.
In this article, we determine some appealing properties of the
proximity map associated with CCCD for data in a compact interval
in $\R$ and use them as guidelines for defining new proximity maps
in higher dimensions.
As CCCD behaves nicely for uniform data in $\R$
(in the sense that the exact and asymptotic distributions of the
domination number and relative arc density are available),
by emulating its properties in higher dimensions
we expect the new PCDs will behave similarly.
Furthermore, we introduce some auxiliary tools used for the
construction of the new proximity maps, as well as some
related concepts that will be used in the investigation and
comparison of the proximity maps.
Additionally, we discuss the conditions for the geometry invariance
for uniform data in triangles.

We describe the data-random PCDs in Section \ref{sec:PCD},
Voronoi diagrams and Delaunay tessellations in Section \ref{sec:voronoi-delaunay-tessellation},
the appealing properties of spherical proximity maps in $\R$
in Section \ref{sec:prox-maps},
transformations preserving uniformity on triangles in $\R^2$ in Section \ref{sec:transformations},
triangle centers in Section \ref{sec:triangle-centers},
vertex and edge regions in Section \ref{sec:vertex-edge-regions}, and
proximity regions in Delaunay tessellations in Section \ref{sec:prox-regions-in-Deltri}.
We present the results on relative arc density and the domination number
of the new PCDs in Section \ref{sec:rel-dens-dom-num},
introduce two new proximity maps in Section \ref{sec:new-prox-map},
and discussion and conclusions in Section  \ref{sec:disc-conc}.

\section{Data-Random Proximity Catch Digraphs}
\label{sec:PCD}
\cite{priebe:2001} introduced the class
cover catch digraphs (CCCDs) and gave the exact and the asymptotic
distribution of the domination number of the CCCD based on two
classes $\X_n$ and $\Y_m$ both of which are random samples from
uniform distribution on a compact interval in $\R$.
\cite{devinney:2002a}, \cite{marchette:2003},
\cite{priebe:2003b}, \cite{priebe:2003a}, and \cite{devinney:2006}
applied the concept in higher dimensions and
demonstrated relatively good performance of CCCDs in classification.
The methods employed involve \emph{data reduction}
(\emph{condensing}) by using approximate minimum dominating sets as
\emph{prototype sets} since finding the exact minimum dominating set is in
general an NP-hard problem
--- in particular, for CCCDs --- (see \cite{devinney:Phd-thesis}).
Furthermore, the exact and the asymptotic distribution of the
domination number of the CCCDs are not analytically tractable in
dimensions greater than 1.

Let $\X_n=\bigl\{X_1,\ldots,X_n \bigr\}$ and $\Y_m=\bigl\{ Y_1,\ldots,Y_m \bigr\}$
be two sets of $\Omega$-valued random variables from classes $\X$ and $\Y$
whose joint pdf is $f_{X,Y}$.
Let $d(\cdot,\cdot):\Omega \times \Omega \rightarrow [0,\infty)$ be a distance function.
The class cover problem for a target class, say $\X_n$, refers to finding a
collection of neighborhoods, $N(X_i)$ around $X_i \in \X_n$ such that
(i) $\X_n \subseteq \bigl( \bigcup_i N(X_i) \bigr)$ and
(ii) $\Y_m \bigcap \bigl( \bigcup_i N(X_i)\bigr) =\emptyset$.
A collection of neighborhoods satisfying both conditions is called a {\em class cover}.
A cover satisfying condition (i) is a {\em proper cover} of
class $\X$ while a collection satisfying condition (ii) is a {\em
pure cover} relative to class $\Y$.
From a practical point of view,
for example for classification, of particular interest are the class
covers satisfying both (i) and (ii) with the smallest collection of
neighborhoods, i.e.,  minimum cardinality cover.
This class cover problem is a generalization of the set cover
problem in \cite{garfinkel:1972} that emerged in statistical pattern
recognition and machine learning, where an edited or condensed set
(i.e., prototype set) is selected from $\X_n$ (see, e.g., \cite{devroye:1996}).

In particular, the proximity regions are constructed using data sets
from two classes.
Given $\Y_m \subseteq \Omega$, the {\em proximity map}
$\NY(\cdot): \Omega \rightarrow \wp(\Omega)$ associates a {\em
proximity region} $\NY(x) \subseteq \Omega$ with each point $x \in
\Omega$. The region $\NY(x)$ is defined in terms of the distance
between $x$ and $\Y_m$.
More specifically, our proximity maps will be
based on the relative position of points from class $\X$ with
respect to the Delaunay tessellation of the points from class $\Y$.
See \cite{okabe:2000} for more on Delaunay tessellations.

If $\X_n=\left\{ X_1,\ldots,X_n \right\}$ is a set of $\Omega$-valued
random variables then $\NY\left( X_i \right)$ are random sets.
If $X_i$ are independent identically distributed then so are the
random sets $\NY\left( X_i \right)$.
The data-random PCD $D$ ---
associated with $\NY(\cdot)$ --- is defined with vertex set
$\X_n=\{X_1,\cdots,X_n\}$ and arc set $\A$ by $X_iX_j \in \A \iff X_j \in \NY\left( X_i \right)$.
Since this relationship is not symmetric, a
digraph is used rather than a graph.
The random digraph $D$
depends on the (joint) distribution of the $X_i$ and on the map $\NY(\cdot)$.
Let $\mu(\NY):=P(X_iX_j \in \A)=P(X_j \in \NY\left( X_i \right))$;
so $\mu(\NY)$ is the probability of having an arc from $X_i$ to $X_j$,
hence is called \emph{arc probability} for the PCD based on $\NY$.

The PCDs are closely related to the {\em proximity graphs} of \cite{jaromczyk:1992} and
might be considered as a special case of {\em covering sets} of
\cite{tuza:1994} and {\em intersection digraphs} of \cite{sen:1989}.
This data random proximity digraph is a {\em
vertex-random proximity digraph} which is not of standard type
(see, e.g., \cite{janson:2000}).
The randomness of the PCDs lies in the fact that the vertices are random with joint pdf
$f_{X,Y}$, but arcs $X_iX_j$ are deterministic functions of the
random variable $X_j$ and the set $\NY\left( X_i \right)$.
For example, the CCCD of \cite{priebe:2001} can be
viewed as an example of PCD with $\NY(x)=B(x,r(x))$, where
$r(x):=\min_{\y \in \Y_m}d(x,\y)$.
The CCCD is the digraph of order
$n$ with vertex set $\X_n$ and an arc from $X_i$ to $X_j$ iff $X_j\in B(X_i,r(X_i))$.
That is, there is an arc from $X_i$ to $X_j$ iff
there exists an open ball centered at $X_i$
which is ``pure" (or contains no elements) of $\Y_m$,
and simultaneously contains (or ``catches") point $X_j$.
%Such a digraph is called \emph{$\mathscr D_{n,m}$-digraph}
%($\mathscr C_{n,m}$-graph in \cite{priebe:2001}),
%which is in fact a pseudodigraph since it permits loops.
%The change in the notation is to emphasize the fact that $\mathscr D_{n,m}$ is a digraph.

Notice that the CCCDs are defined with (open) balls only, whereas
PCDs are not based on a particular geometric shape or a functional
form; that is, PCDs admit $\NY(\cdot)$ to be any type of region,
e.g., circle (ball), arc slice, triangle, a convex or nonconvex polygon, etc.
In this sense, the PCDs are defined in a more general
setting compared to CCCDs.
On the other hand, the types of PCDs introduced in this article are well-defined for points
restricted to the convex hull of $\Y_m$, $\C_H\left( \Y_m \right)$.
Moreover, the proximity maps introduced in this article
will yield closed regions.
Furthermore, the CCCDs based on balls use proximity regions that are
defined by the obvious metric,
while the PCDs do not necessarily require a metric,
but some sort of dissimilarity measure only.

\section{Voronoi Diagrams and Delaunay Tessellations}
\label{sec:voronoi-delaunay-tessellation}
The proximity map in $\R$ defined as $B(X,r(X))$ where
$r(X)=\min_{Y\in \Y_m}d(X,Y)$ in (\cite{priebe:2001}).
Our next goal it to extend this idea to higher dimensions and investigate the
properties of the associated digraph.
Now let $\Y_m=\left \{\y_1,\ldots,\y_m \right\} \subset \R$
and $Y_{j:m}$ be the $j^{th}$ order statistic.
The above definition of the proximity map is based on the intervals
$I_j=\left( Y_{(j-1):m},Y_{j:m} \right)$ for $j=0,\ldots,(m+1)$ with
$\y_{0:m}=-\infty$ and $\y_{(m+1):m}=\infty$.
This \emph{intervalization} can be viewed as a tessellation
since it partitions $\C_H(\Y_m)$, the convex hull of $\Y_m$.
For $d>1$, a natural tessellation that
partitions $\C_H(\Y_m)$ is the Delaunay tessellation, where each Delaunay cell is a
$(d+1)$-simplex (e.g., a Delaunay cell is an interval for $d=1$, a
triangle for $d=2$, a tetrahedron for $d=3$, and so on.)
Let $\T_j$ be the $j^{th}$ Delaunay cell for $j=1,\ldots,J$
in the Delaunay tessellation of $\Y_m$.
In $\R$, the cell that contains $x$ is implicitly used
to define the proximity map.

\begin{figure}[ht]
\centering
\psfrag{Voronoi Diagram}{\Huge{Voronoi Diagram}}
\psfrag{Delaunay Triangulation}{\Huge{Delaunay Triangulation}}
\epsfig{figure=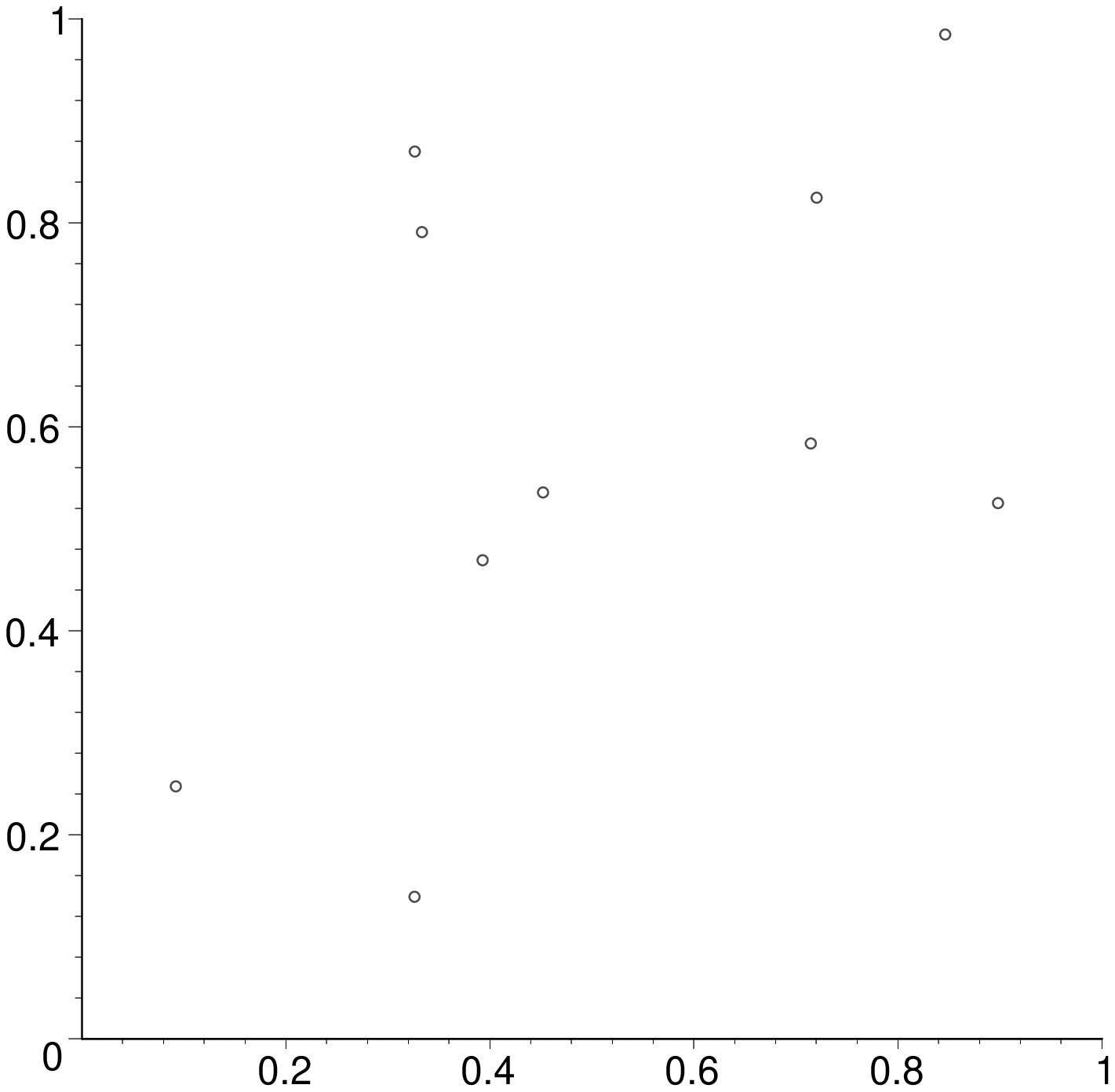, height=125pt , width=125pt}\\
\rotatebox{-90}{ \resizebox{2.2 in}{!}{ \includegraphics{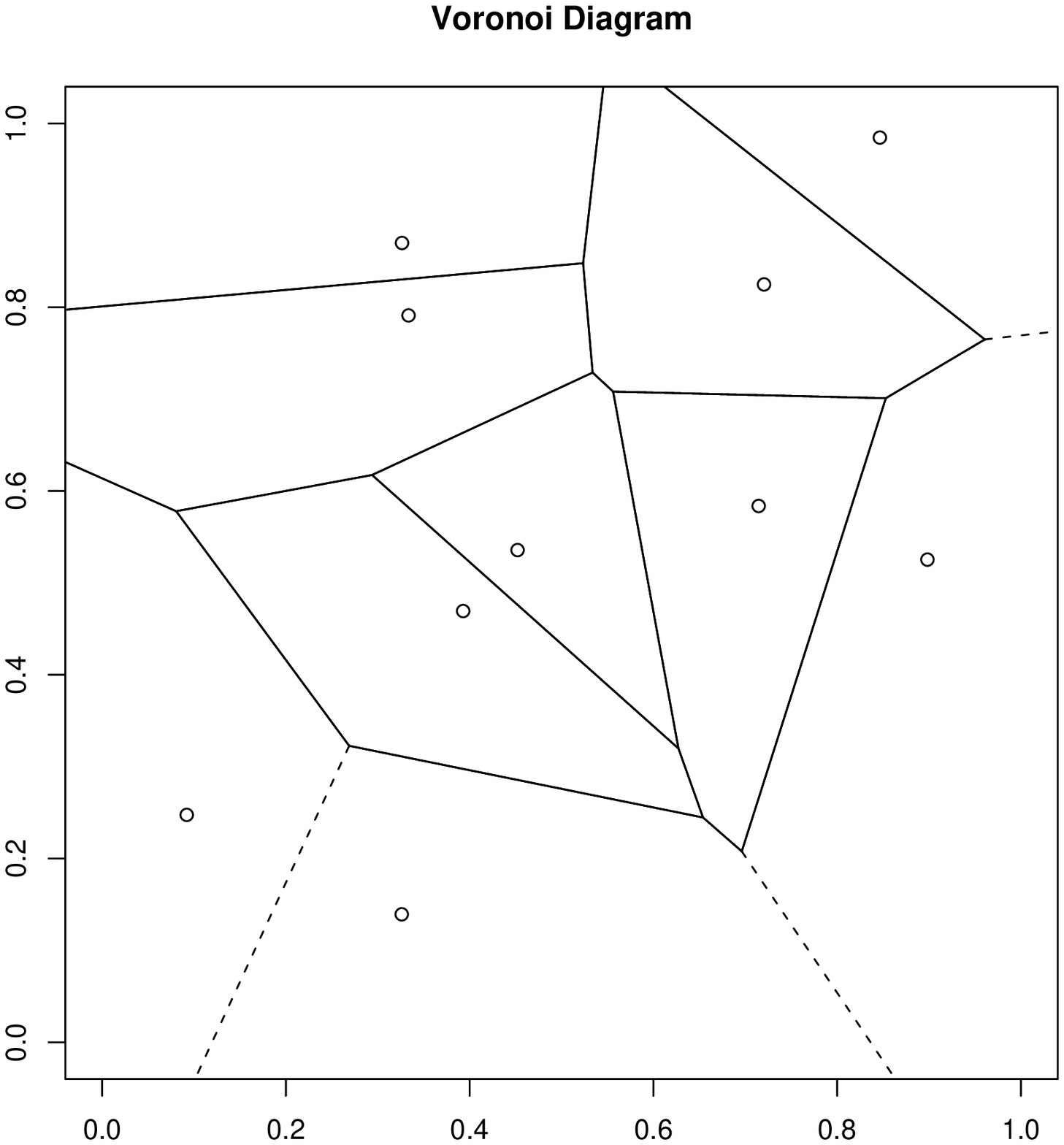}}}
\rotatebox{-90}{ \resizebox{2.2 in}{!}{ \includegraphics{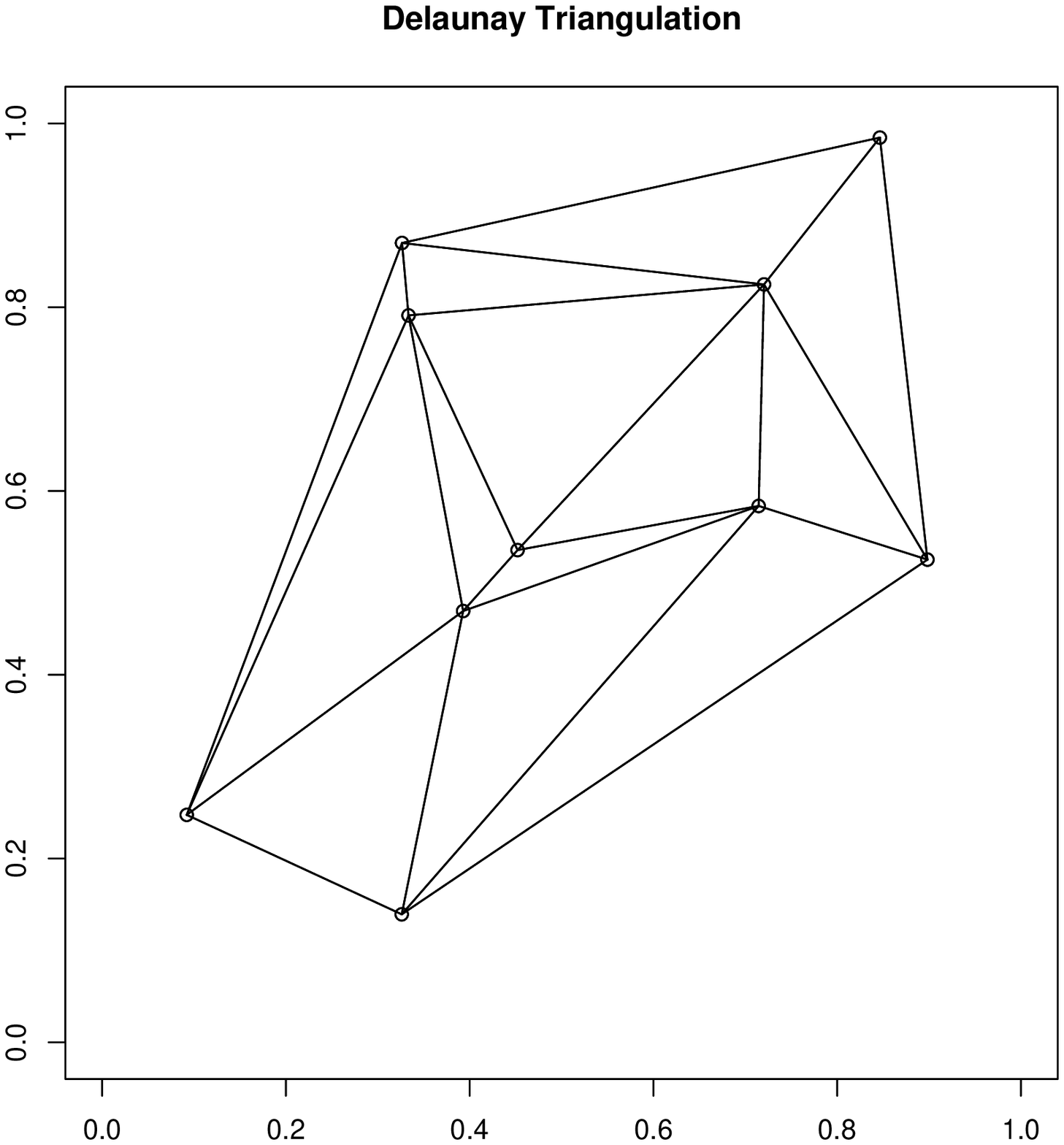}}}
\caption{
\label{fig:Vor-Del-tri}
Depicted are 10 $\Y$ points generated iid $\U(0,1)\times(0,1)$ (top),
the corresponding Voronoi diagram (bottom left) and the
Delaunay triangulation (bottom right).
}
\end{figure}

A \emph{tessellation} is a partition of a space into convex
polytopes; tessellation of the plane (into convex polygons) is the
most frequently studied  case (\cite{schoenberg:2002}).
Given $2\le n <\infty$ distinct points in $\R^2$, we associate all
points in the space with the closest member(s) of the point set
with respect to the Euclidean distance.  The result is a tessellation
of the plane into a set of regions associated with the $n$ points.
We call this tessellation the \emph{planar ordinary Voronoi diagram}
generated by the point set and the regions
\emph{ordinary Voronoi polygons}.
See Figure \ref{fig:Vor-Del-tri} for an example with $m=|\Y_m|=10$
$\Y$ points iid from $\U \bigl((0,1)\times (0,1)\bigr)$.

In general, let $P=\{p_1,p_2,\ldots,p_n\}$ be $n$ points in $\R^d$ where
$2 \le n <\infty$ and $p_i \not= p_j$ for $ i\not= j$,
$i,j \in [n]:=\bigl\{ 1,2,\ldots,n \bigr\}$ and
let $|| \cdot ||$ denote the norm functional.
We call the region
$\displaystyle \V_C(p_i)=\bigl\{x \in \R^d:\; ||x-p_i|| \le ||x-p_j|| \text{ for $j \not= i, \; j \in [n]$} \bigr\}$
the (ordinary) \emph{Voronoi polygon} or \emph{cell} associated with $p_i$
and the set $\mathfrak V=\bigl\{ \V_C(p_1),\ldots,\V_C(p_n) \bigr\}$
the \emph{Voronoi diagram} or \emph{Dirichlet tessellation} generated by $P$.
We call $p_i$ the \emph{generator} of $\V_C(p_i)$.
The Voronoi diagram partitions the space into disjoint regions
(which are also called \emph{tiles} or \emph{Thiessen polygons} in $\R^2$).
Notice that we still say $\mathfrak V$ \emph{partitions} the space $\R^d$,
although $\V_C(p_j)$ are not necessarily disjoint,
but if nonempty the intersection lies
in a lower dimension, or equivalently, has zero $\R^{d}$-Lebesgue measure.
We stick to this convention throughout the article.
The intersection of two Voronoi cells, if nonempty, i.e., for $i \not= j$,
$\displaystyle \V_C(p_i) \cap \V_C(p_j) \not= \emptyset$, is called a \emph{Voronoi edge}.
If a Voronoi edge is not a point, then $\V_C(p_i)$ and
$\V_C(p_j)$ are said to be \emph{adjacent}.
An end point of a Voronoi edge is called a \emph{Voronoi vertex}.
$\mathfrak V$ is called \emph{degenerate}
if at a Voronoi vertex, more than three Voronoi polygons intersect,
and \emph{non-degenerate} otherwise.
A detailed discussion including the history of
Voronoi diagrams is available in  \cite{okabe:2000}.

\begin{figure}[ht]
\centering
\rotatebox{-90}{ \resizebox{3.0 in}{!}{ \includegraphics{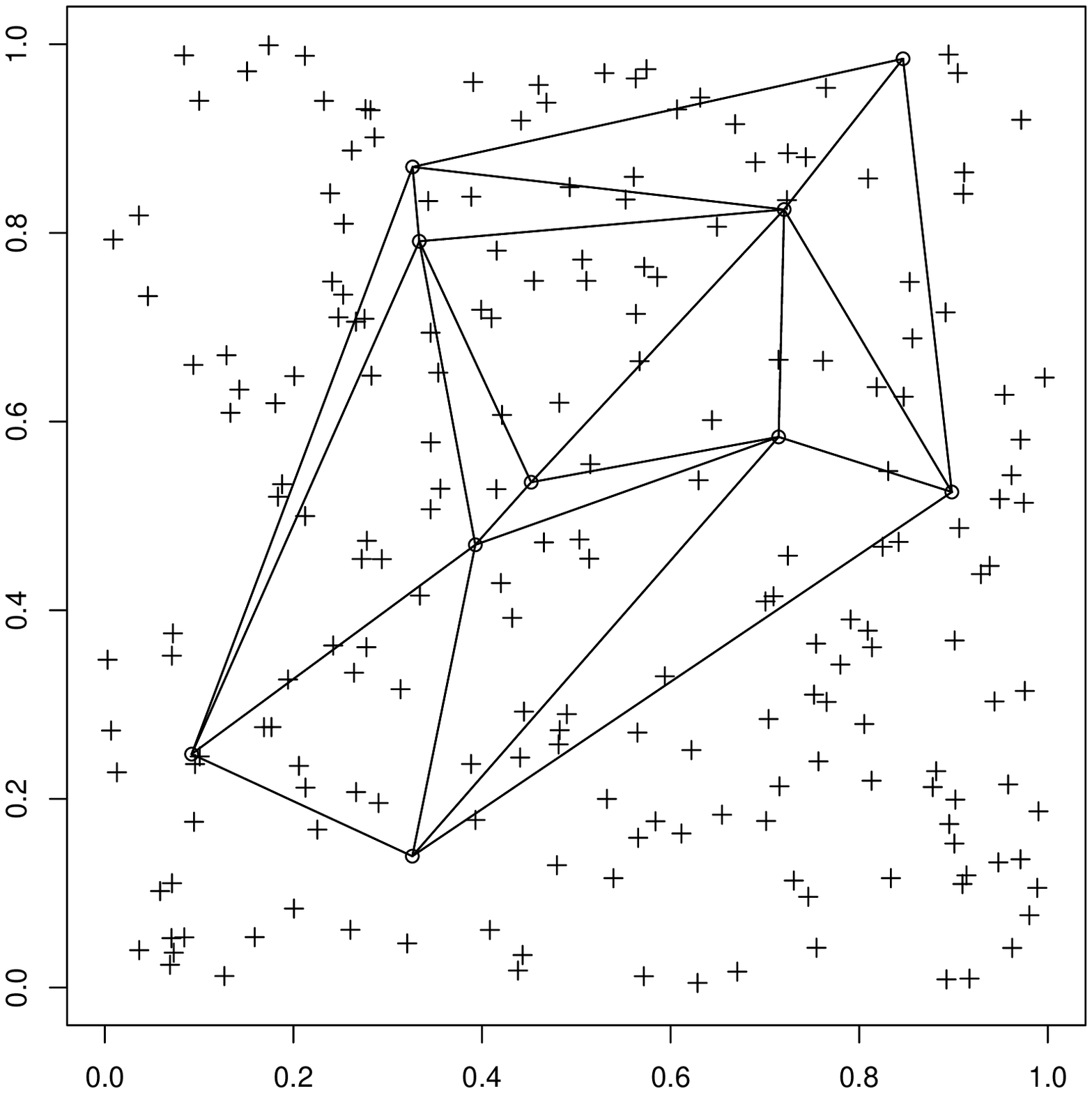}}}
\rotatebox{-90}{ \resizebox{3.0 in}{!}{\includegraphics{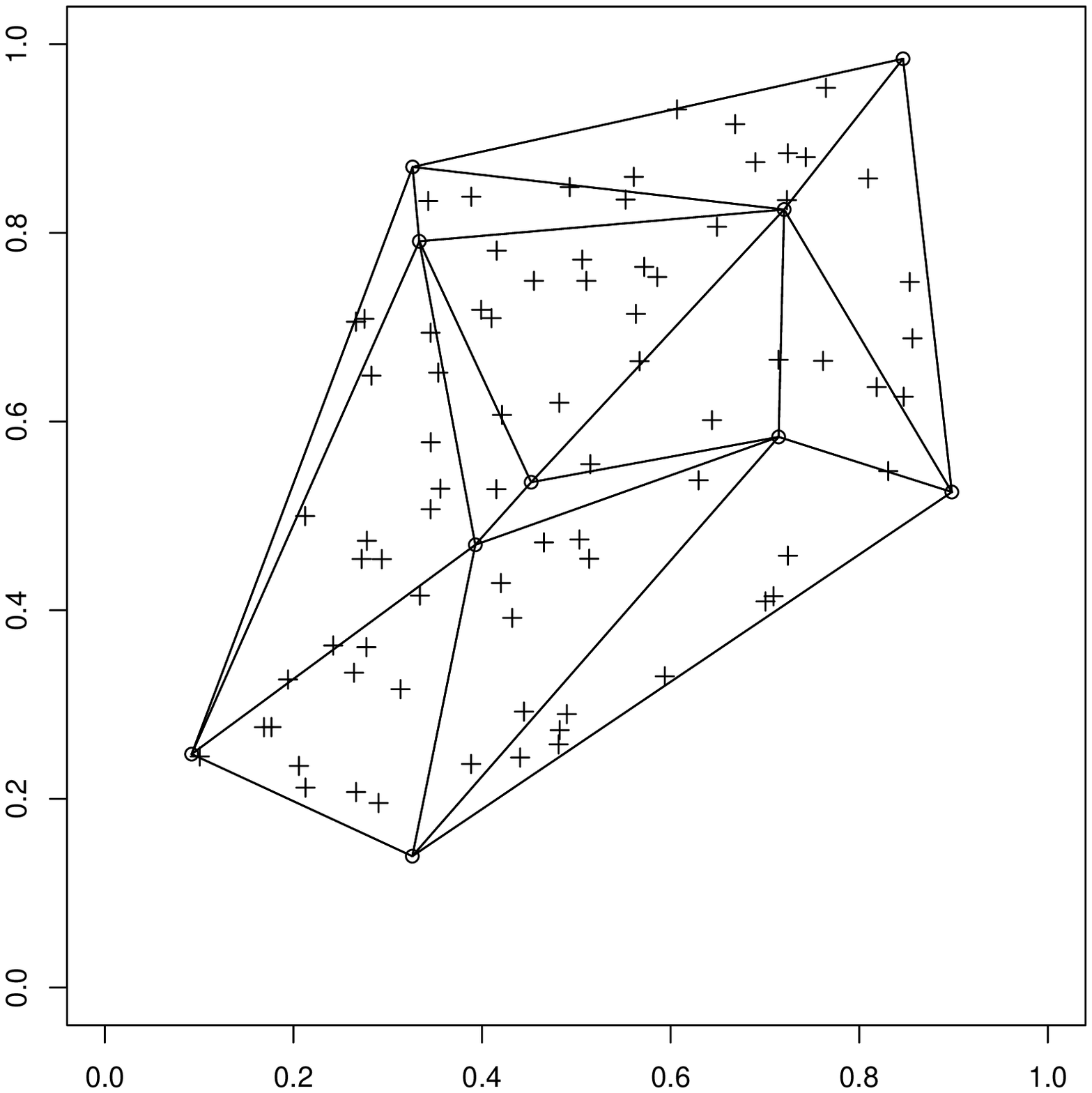} }}\\
\epsfig{figure=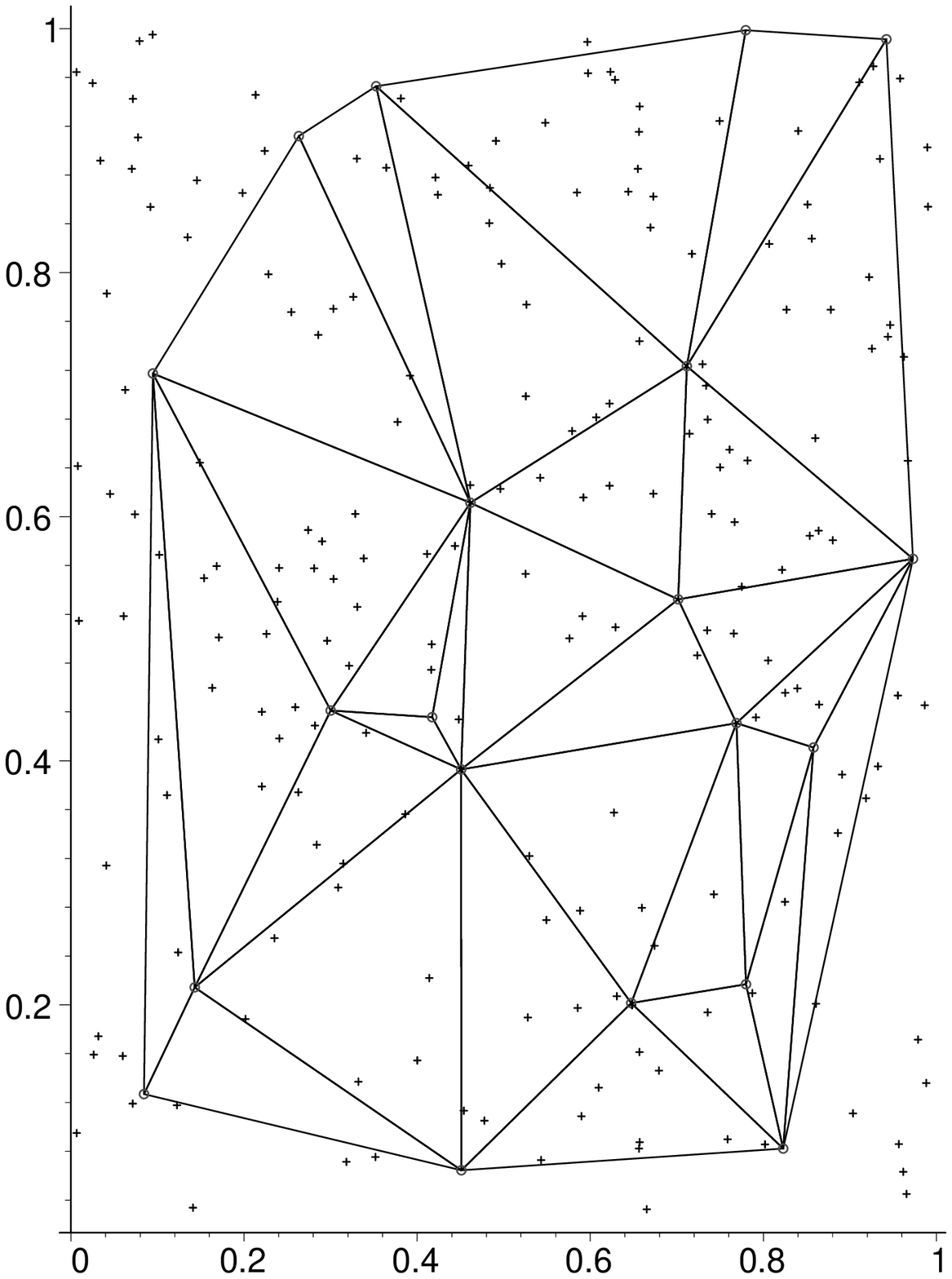, height=160pt , width=160pt}
\caption{\label{fig:deltri}
A realization of $200$ $X$ points and the Delaunay triangulation
based on 10 $\Y$ points in Figure \ref{fig:Vor-Del-tri} (top).
Another realization of 200 $\X$ and 20 $\Y$ points and the Delaunay triangulation based on $\Y_m$
(bottom).
}
\end{figure}

Given a Voronoi diagram with $n \ge d+1$ non-coplanar
(i.e., not all the points lie on a $(d-1)$-dimensional hyperplane)
generators, $P$, in $\R^d$, we join all pairs of generators whose Voronoi cells have
a common Voronoi edge.
The resulting tessellation is called the
\emph{Delaunay tessellation} of $P$.
See Figure \ref{fig:Vor-Del-tri} for the Delaunay triangulation
associated with the Voronoi diagram based on 10 points iid from $\U \bigl((0,1)\times (0,1)\bigr)$.
By definition a Delaunay tessellation
of a finite set of points, $P$, is the dual of the \emph{Voronoi diagram}
based on the same set.
The tessellation yields a (unique)
polytopization provided that no more than $(d+1)$ points in $\R^d$ are
cospherical (i.e., no more than $(d+1)$ points lie on the
boundary of a (hyper)sphere in $\R^d$).
Moreover, the circumsphere of each
Delaunay polytope
(i.e., the sphere that contains the vertices of the Delaunay polytope
on its boundary) is pure from the set $P$;
i.e., the interior of the circumsphere of
the Delaunay polytope does
not contain any points from $P$.
The Delaunay tessellation partitions $\C_H(P)$.
In particular, in $\mathbb R^2$, the
tessellation is a \emph{triangulation} that yields triangles $T_j$,
$j=1,\ldots,J$ (see, e.g., \cite{okabe:2000}) provided that no more than
three points are cocircular (i.e., no more than three points lie on
the boundary of some circle in $\R^2$).
In this article we adopt the convention that
a triangle refers to the closed region bounded by
its edges.
See Figure \ref{fig:deltri} for two examples:
an example with $n=200$ $\X$ points
$\stackrel{iid}{\sim}\U \bigl((0,1)\times (0,1)\bigr)$,
the uniform distribution on the unit square
and the Delaunay triangulation based on the 10 $\Y$ points in Figure \ref{fig:Vor-Del-tri};
and
an example with $n=200$ $\X$ points
$\stackrel{iid}{\sim}\U \bigl((0,1)\times (0,1)\bigr)$,
the uniform distribution on the unit square
and the Delaunay triangulation is based on 20 $\Y$ points also $\stackrel{iid}{\sim}\U \bigl((0,1)\times (0,1)\bigr)$.

\subsection{Poisson Delaunay Triangles}
\label{sec:poisson-delaunay-triangles}
The Delaunay triangles are based on a given set of points $\Y_m$.
The set $\Y_m$ can be assumed to come from a Poisson point process
on a finite region, and in the application of PCDs,
to remove the conditioning on $\Y_m$,
it is suggested that $\Y_m$ comes from a Poisson
point process for prospective research directions.
We briefly describe the Poisson point processes and Poisson Delaunay triangles.

A \emph{stochastic point process} on $\R^d$ is defined to be a process
in which points are generated according to a probability
distribution $P\bigl(|\mathcal Y \cap B|=k\bigr)$,
$k=0,1,2, \ldots$, over any $B \subseteq \R^d$.
For example, a \emph{binomial point process} is a stochastic
point process in which $n$ points are generated over a
bounded set $S \subsetneq \R^d$ according to the uniform distribution.
In particular, if $d=2$ the process is called a
\emph{planar stochastic point process}.  If two
points coincide with probability zero, then it
is a \emph{simple stochastic point process}.
A stochastic point process is said to be \emph{locally finite}
if any finite region $B \subset \R^d$ contains a finite number of
points with probability 1 under the process (see \cite{okabe:2000}).

We have built the Delaunay tessellation using $\Y_m$ with finite sample size.
Suppose $\Y$ is from a stochastic point process.
One of the most fundamental locally finite stochastic point processes is the
\emph{Poisson point process},
which is defined as the process that satisfies
$$P\bigl( |\Y \cap B|=k \bigr)=\frac{\lambda\,\cdot \text{V}(B)\cdot e^{-\lambda\,\cdot \text{V}(B)\,\cdot k}}{k!},\;\;k=0,1,2,\ldots$$
for any $B \subsetneq \R^d$, where $\text{V}(\cdot)$ denotes the $d$-dimensional
volume functional and $\lambda>0$ is the intensity
(number of points per unit volume) of the process.
 We can also define the Poisson point process as
the limit of the binomial point process in the sense of expanding
the finite region $S$ to an infinite region while keeping
$\lambda=n/\text{V}(S)$ constant.
We call the Delaunay tessellation based on a finite data set from a
Poisson point process \emph{Poisson Delaunay tessellation}
and denote it $\D_P$.
The associated Voronoi diagram is
called the \emph{Poisson Voronoi diagram} and is denoted by $\V_P$.
%Note that the cells of $\V_P$ are convex also.
For more detail on the properties of $\V_P$, see, \cite{okabe:2000}.

A \emph{simplex} in $\R^d$ is the convex hull of any $(d+1)$ points in general
position, i.e., no $(d+1)$ of the points lie in a $(d-1)$-dimensional
hyperplane in $\R^d$.
The simplex is the point itself for $d=0$, the line segment joining
the two points for $d=1$, a triangle for $d=2$, a tetrahedron for $d=3$, and so on.
Each Poisson Delaunay cell of $\D_P$ is a
$(d+1)$-dimensional simplex whose vertices are $x_0,x_1,\ldots,x_m$
from the Poisson point process.
Any $s$-face of $\D_P$
is an $(s+1)$-dimensional simplex with vertices $x_0, \ldots, x_s$,
also points from the Poisson point process.
There are $\displaystyle {d+s \choose s+1}$ many
$s$-faces contained in a Poisson Delaunay cell for $0 \le s \le d$.

%Let $c$ and $r$ be the circumcenter and circumradius, respectively,
%of a $(d+1)$-dimensional Poisson Delaunay cell in $\R^d$.
%Then the $(d+1)$ vertices
%of the cell are the points $\{c+r\,u_i\}$ where $\{u_i\}$ are the
%unit vectors for $i=0,1,\ldots,d$.
%The ergodic joint probability density function (pdf) of $\D_P$,
%the pdf of $r$, and $k^{th}$ moment of the volume of a typical Poisson Delaunay cell
%are provided in (\cite{okabe:2000}).

Let $c$ and $r$ be the circumcenter and circumradius, respectively,
of a $(d+1)$-dimensional Poisson Delaunay cell in $\R^d$.
Then the $(d+1)$ vertices
of the cell are the points $\{c+r\,u_i\}$ where $\{u_i\}$ are the
unit vectors for $i=0,1,\ldots,d$.
The ergodic joint probability density function (pdf) of $\D_P$ is completely specified as
$$f(r,u_0,\ldots,u_d)=a(\lambda,d)\,\Delta_d\,r^{d^2-1}\,\exp{\left(-\lambda\,w_d\,r^d\right)}$$
where $\Delta_d$ is the volume of the $(d+1)$-simplex with vertices $u_0,\ldots,u_d$, $w_d=\pi^{d/2}/\G(d/2+1)$ and
$$a(\lambda,d)=\frac{\pi^{(d^2+1)/2}\,\G(m^2/2)\,\left[2\,\lambda\, \Gamma\bigl((d+1)/2\bigr)\right]^d}
{d^{d-2}\,\G(d/2)^{2\,d+1}\,\Gamma\bigl((d^2+1)/2\bigr)}.$$

The circumradius $r$ may be viewed as a measure of size of the
$(d+1)$-simplex, and is independent of $\{u_0,\ldots,u_d\}$.
The pdf of $r$ is a generalized gamma function
with $t=d,\;q=d^2$ and $b=\lambda_d$, where a 3-parameter
generalized gamma function is
$$f(x)=r\,b^{q/t}\,x^{q-1}\,\exp{(-b\,x^t)}/\G\bigl(q(t)\bigr)\;\text{ for } x,t,p,q>0.$$
Let $V_d$ denote the volume of a typical
Poisson Delaunay cell. Then for $d=2$, the expected value of the
$k^{th}$ moment of volume of a typical Poisson Delaunay triangle is
$$\E\left[ V_2^k \right]=\frac{\G\bigl((3\,k+5)/2\bigr)\,\G(k/2+1)}{3\,\G\bigl((k+3)/2\bigr)^2\,2^k\,\pi^{k-1/2}\,\lambda^k},\,\;\;k=1,2,\ldots$$

In $\R^2$, the joint pdf of a pair of inner angles arbitrarily selected from an
arbitrary triangle in $\D_P$ is given by
$$f(x,y)=\frac{8}{3 \pi}(\sin x)(\sin y)\sin (x+y), \text{  for $x,y >0$
and $x+y <\pi$.} $$  Notice that the mode of this density is at $x=y=\pi/3$,
which implies that the most frequent triangles in a $\D_P$ are nearly
equilateral triangles.
By integrating over $y$, we obtain the pdf of
a randomly selected inner angle of an arbitrary triangle from $\D_P$:
$$f(x)=\left[\frac{4}{3 \pi} \left( (\pi-x)\cos x+\sin x \right) \sin x \right]\; \I(0<x<\pi)$$
where $\I(\cdot)$ is the indicator function.
Then the expected value of $X$ is $\E[X]=\pi/3$ and $\E\left[ X^2 \right]=2\,\pi^2/9-5/6$. %See \cite{Col1968}.
%The pdf  of the minimum angle and the
%pdf  of the maximum angle, and the
%distribution of the length of an arbitrary edge of an arbitrary
%triangle from $\D_P$ are also provided in \cite{okabe:2000}
%with relevant references.

\cite{mardia:1977} derived the pdf of the minimum angle as
$$f_1(x)=\left(\frac{2}{\pi} (\pi-3 x) \sin 2x + \cos 2x-\cos 4x\right) \I\left( 0<x<\frac{\pi}{3} \right),$$
and \cite{boots:1986} gave the pdf  of maximum angle,
\begin{multline*}
f_3(x)=\left[\frac{2}{\pi}(3x (\sin 2x)-\cos 2x +
\cos 4x -\pi \sin 2x)\right]\,\I\left(\frac{\pi}{3}<x<\frac{\pi}{2}\right)+ \\
\left[\frac{1}{\pi}\left(4 \pi (\cos x) (\sin x) +
3 \sin x^2- \cos x^2-4x(\cos x)(\sin x) +1\right)\right]\, \I\left(\frac{\pi}{2}<x<\pi\right).
\end{multline*}
The distribution of the length of an arbitrary edge $e$ of an arbitrary triangle from $\D_P$ is
$$f_L(x)=\left(\pi \lambda \frac{x}{3}\right)\left(\sqrt{\lambda}x \exp\left(-\pi \lambda \frac{x^2}{4}\right)+
\text{erfc}\left(\sqrt{\pi \lambda}\frac{x}{2}\right)\right)$$
where $\text{erfc}$ is the complimentary error function
$\displaystyle \text{erfc}(z)=\frac{2}{\sqrt{\pi}}\int_z^{\infty}
e^{-t^2}dt=\pi^{-\frac{1}{2}}\G\left(\frac{1}{2},z^2\right)$ (\cite{okabe:2000}).

\section{The Appealing Properties of Spherical Proximity Maps in $\R$}
\label{sec:prox-maps}
The CCCDs have desirable properties such as
having the finite sample and asymptotic distributions of the domination number available.
In this section, we determine some appealing properties of the
proximity map associated with CCCD for uniform data in a compact interval
in $\R$ and use them as guidelines for defining new proximity maps
in higher dimensions.
We believe these properties cause the CCCD to behave so ``nicely"
in $\R$ and the more they are satisfied by the new PCDs in higher dimensions,
the more likely the new PCDs to have similar behaviour.
Furthermore, we introduce the auxiliary tools used for the
construction of the new proximity maps, as well as some
related concepts that will be used in the investigation and
comparison of the proximity maps.

Let $\Y_m=\left \{\y_1,\ldots,\y_m \right\} \subset \R$.
Then the proximity map associated with CCCD is defined
as the open ball $\NS(x):=B(x,r(x))$ for all $x \in \R \setminus \Y_m$,
where $r(x)=\min_{\y\in \Y_m}d(x,\y)$ (see Section \ref{sec:PCD} and \cite{priebe:2001})
with $d(x,y)$ being the Euclidean distance between $x$ and $y$.
For $x \in \Y_m$, define $\NS(x)=\{x\}$.
Notice that a ball is a sphere in higher dimensions,
hence the name \emph{spherical proximity map} and the notation $N_S$.
Furthermore, dependence on $\Y_m$ is through $r(x)$.
Note that, this proximity map is based on
the intervals $I_i=\left(\y_{(i-1):m},\y_{i:m} \right)$
for $i=0,\ldots,(m+1)$ with $\y_{0:m}=-\infty$ and $\y_{(m+1):m}=\infty$
where $y_{i:m}$ is the $i^{th}$ order statistic of $\Y_m$.

For $X_i \stackrel{iid}{\sim} \U(I_i)$,
without loss of generality we can assume $I_i=(0,1)$.
Then the arc probability $\mu(N_S)=P(X_2 \in N_S(X_1))=1/2$,
since $\displaystyle P(X_2 \in N_S(X_1))=\int_0^{1/2}\int_0^{2x_1} dx_2dx_1+
\int_{1/2}^1\int_{2-2x_1}^1 dx_2dx_1=1/2$.

A natural extension of the proximity region $\NS(x)$
to multiple dimensions (i.e., to $\R^d$ with $d>1$)
is obtained by the same definition as above;
that is, $N_S(x):=B(x,r(x))$ where $r(x):=\min_{\y \in \Y_m} d(x,\y)$.
The spherical proximity map $N_S(x)$ is well-defined for all
$x \in \R^d$ provided that $\Y_m \not= \emptyset$.
Extensions to
$\R^2$ and higher dimensions with the spherical proximity
map --- with applications in classification --- are
investigated in  \cite{devinney:2002a}, \cite{marchette:2003},
\cite{priebe:2003b}, \cite{priebe:2003a}, and \cite{devinney:2006}.
However, finding the minimum
dominating set of the PCD associated with $N_S(\cdot)$ is an NP-hard problem (\cite{devinney:Phd-thesis}) and
the distribution of the domination number
is not analytically tractable for $d>1$ (\cite{ceyhan:Phd-thesis}).
This drawback has motivated the definition of new types of proximity maps
in higher dimensions.
Note that for $d=1$, such problems do not exist.
We state some appealing properties of the proximity map
$N_S(x)=B(x,r(x))$ in $\R$ and use them as guidelines for
our definition of new proximity maps:
\begin{itemize}
\item[\textbf{P1}]
$\NS(x)$ is well-defined for all $x \in \C_H\left( \Y_m \right)=[\y_{1:m},\y_{m:m}]$.
\item[\textbf{P2}]
$x \in \NS(x)$ for all $x \in \C_H\left( \Y_m \right)$.
\item[\textbf{P3}]
The point $x$ is at the \emph{center} of $\NS(x)$ for all $x \in \C_H\left( \Y_m \right)$.
\item[\textbf{P4}]
For $x \in I_i \subseteq \C_H\left( \Y_m \right)$, $\NS(x)$ and $I_i$
are of the \emph{same type}; i.e., they are both intervals.
\item[\textbf{P5}]
For $x \in I_i \subseteq \C_H\left( \Y_m \right)$,
$\NS(x)$ mimics the shape of $I_i$; i.e.,  it is (geometrically) \emph{similar} to $I_i$.
\item[\textbf{P6}]
For $x \in I_i$, $\NS(x)$ is a proper subset
of $I_i$ for all $x \in I_i \setminus \{(\y_{(i-1):m}+\y_{i:m})/2\}$
(or almost everywhere in $I_i$).
\item[\textbf{P7}]
For $x \in I_i$ and $y \in I_j$ with $i \not= j$,
$\NS(x)$ and $\NS(y)$ are disjoint regions.
\item[\textbf{P8}]
The size (i.e., measure) of $\NS(x)$ is continuous in $x$; that is,
for each $\ve >0$ there exists a $\delta(\ve)>0$
such that $\bigl||\NS(y)|-|\NS(x)|\bigr|<\ve$
whenever $|d(x,y)|<\delta(\ve)$.
%\item[\textbf{P9}]
%The size (i.e., measure) of $\NS(x)$ increases as $d(x,\Y_m)$ increases.
\item[\textbf{P9}]
The arc probability $\mu(N_S)$ does not depend on the support interval for uniform data in $\R$.
\end{itemize}

Notice that properties \textbf{P1}, \textbf{P2}, and \textbf{P3}
also hold for all $x \in \R$.
\textbf{P9} implies that not only the arc probability but also
the distribution of the relative arc density and domination
number do not depend on the support interval either.
This independence of the support set is called
\emph{geometry invariance} in higher dimensions
(see Section \ref{sec:transformations}).
For $N_S$ it suffices to work with $\U(0,1)$ data,
and in higher dimensions we will be able to consider only uniform data
in an equilateral triangle for PCDs based on
proximity maps that satisfy \textbf{P9}.

Suppose we partition the convex hull of $\Y_m$, $\C_H\left( \Y_m \right)$ %, in $\R$ to intervals
%can be viewed as a tessellation.
%For $d>1$, a natural tessellation that partitions $\C_H\left( \Y_m \right)$ is
%the Delaunay tessellation, where each Delaunay cell is a $(d+1)$-simplex
%(e.g., a Delaunay cell is an interval for $d=1$, a triangle for $d=2$,
%a tetrahedron for $d=3$, etc.)
%Furthermore, Delaunay tessellation is the graph theoretic dual
%of Voronoi diagrams.
by Delaunay tessellation.
Let $\T_j$ be the $j^{th}$ Delaunay cell in the
Delaunay tessellation of $\Y_m$ for $j=1,\ldots,J$,
where $J$ is the total number of Delaunay cells.
See Figure \ref{fig:deltri} for two $\Y$ sets
of sizes 10 and 20 and the corresponding Delaunay triangulations.
%In $\R$, the Delaunay cell that contains $x$ is implicitly used
%to define the proximity map for $x \in \C_H\left( \Y_m \right)$.

%\begin{figure}[ht]
%\begin{center}
%\epsfig{figure=Ypoints.eps, height=111pt , width=148pt}
%%\epsfig{figure=Voronoi.ps, height=120pt , width=160pt}
%\epsfig{figure=DTofY.eps, height=111pt , width=148pt}
%\end{center}
%\caption{
%\label{fig:delaunay}
%Depicted are 10 $\Y$ points generated iid $\U(0,1)\times(0,1)$ (left),
%the corresponding
%%Voronoi diagram (middle) and the
%Delaunay triangulation (right).}
%\end{figure}

Note that \textbf{P4} and \textbf{P5} are equivalent when $d=1$
for $x \in \C_H\left( \Y_m \right)$, since any two (compact) intervals in $\R$
are (geometrically) similar.
For $d>1$, \textbf{P5} implies \textbf{P4} only, since, for example,
for $d=2$, any two triangles are not necessarily similar,
but similar triangles are always of the same type;
they are triangles.

Notice that $N_S(\cdot)$ satisfies only \textbf{P1}, \textbf{P2},
\textbf{P3}, and \textbf{P8} in $\R^d$ with $d>1$.
\textbf{P4} and \textbf{P5} fail since $N_S(x)$ is a sphere for $x \not\in \Y_m$,
but $\T_j$ is a $(d+1)$-simplex.
For any $x \in \T_j \subset \R^d$,
$B(x,r(x)) \not \subset \T_j$, so \textbf{P6} also fails,
furthermore this also implies that
$N_S(x)$ and $N_S(y)$ might overlap
for $x,y$ from two distinct cells,
hence \textbf{P7} is violated.
The arc probability $\mu(N_S)$ depends on the support set $\T_i$ for $d>1$
so \textbf{P9} is violated.

The appealing properties mentioned above can be extended to more general measurable spaces.
Let $(\Omega,\M)$ be a measurable space, and let $\Omega_j$, $j\in \{1,2,\ldots,J\}$
partition $\Omega$, and $\mu$ be the associated measure on $\Omega$.
Then the appealing properties are
\begin{itemize}
\item[\textbf{P1}]
$N(x)$ is well-defined for all $x \in \Omega$.
\item[\textbf{P2}]
$x \in N(x)$ for all $x \in \Omega$.
\item[\textbf{P3}]
$x$ is at the \emph{center} of $N(x)$ for all $x \in \Omega$.
\item[\textbf{P4}]
For $x \in \Omega_j$, $N(x)$ and $\Omega_j$
are of the \emph{same type}; they have the same functional form.
\item[\textbf{P5}]
For $x \in \Omega_j \subseteq \Omega$,
$N(x)$ mimics the shape of $\Omega_j$; i.e., it is \emph{similar} to $\Omega_j$.
\item[\textbf{P6}]
For $X \in \Omega_j$, $N(X)$ is a proper subset of $\Omega_j$ a.s.
\item[\textbf{P7}]
For $x \in \Omega_j$ and $y \in \Omega_k$
with $j\not=k$, $N(x)$ and $N(y)$ are disjoint. %a.s.
\item[\textbf{P8}]
The measure of $N(x)$ is continuous in $x$;
that is, for each $\ve >0$ there exists a $\delta(\ve)>0$ such that
$|\mu(N(y))-\mu(N(x))|<\delta(\ve)$ whenever $||y-x|| <\ve$.
%\item[\textbf{P9}]
%The measure of $\NS(x)$ increases as $d(x,\Y_m)$ increases.
\item[\textbf{P9}]
The arc probability $\mu(N_S)$ does not depend on
the support set for uniform data in $\Omega$.
\end{itemize}

Property \textbf{P6} suggests a new concept.
\begin{definition}
The \emph{superset region} for any proximity map $N(\cdot)$
in $\Omega$ is defined to be $\RS(N):=\bigl\{ x \in \Omega: N(x) = \Omega \bigr\}$. $\square$
\end{definition}
For example, for $\Omega=I_i \subsetneq \R$,
$\RS(\NS):=\{x \in I_i: \NS(x) = I_i\}=\left\{ \left(\y_{(i-1):m}+\y_{i:m} \right)/2 \right\}$,
and for $\Omega=\T_i \subsetneq \R^d$,
$\RS(\NS):=\{x \in \T_i: \NS(x) = \T_i\}$.
Note that for $x \in I_i$, $\lambda(\NS(x)) \le \lambda(I_i)$
and $\lambda(\NS(x)) = \lambda(I_i)$ iff $x \in \RS(\NS)$
where $\lambda(\cdot)$ is the Lebesgue measure on $\R$
(also called $\R$-Lebesgue measure).
So the proximity region of a point in $\RS(\NS)$ has the largest
$\R$-Lebesgue measure.
Note that for $\Y_m=\left\{ \y_1,\ldots,\y_m \right\} \subset \R$,
$\displaystyle \RS(\NS)=\left\{ \frac{\y_{1:m}+\y_{2:m}}{2},\ldots,\frac{\y_{(m-1):m}+\y_{m:m}}{2}\right \}$.
Note also that given $\Y_m$, $\RS(\NS)$ is not a random set,
but $\I(X\in \RS(\NS))$ is a random variable.
Furthermore, \textbf{P6} is equivalent to $\RS(\NS)$ having
zero $\R$-Lebesgue measure.
On the other hand, for $x\in \partial(I_i)=\{\y_{(i-1):m},\y_{i}\}$,
the proximity region $\NS(x)=\{x\}$ has zero $\R$-Lebesgue measure.
This suggests the following concept.

\begin{definition}
Let $(\Omega,\mu)$ be a measurable space.
The \emph{$\Lambda_0$-region} for any proximity map $N(\cdot)$ is
defined to be $\Lambda_0(N):=\bigl\{ x \in \Omega: \mu(N(x)) = 0 \bigr\}$. $\square$
\end{definition}
For $\Omega=\R^d$, $\Lambda_0(\NS):=\bigl\{ x \in \R^d: \lambda(\NS(x)) = 0 \bigr\}$.
For example, for $\Omega=\C_H\left( \Y_m \right) \subsetneq \R$, $\Lambda_0(\NS)=\Y_m$,
since $\lambda(\NS(x))=0$ iff $x \in \Y_m$.

Furthermore, given a set $B$ of size $n$ in $[\y_{1:m},\y_{m:m}] \setminus \Y_m$,
\textbf{P7} implies that the number of
disconnected components in the PCD based on $N_S(\cdot)$
is at least the cardinality of $\{i \in [m]: B \cap I_i \not=\emptyset \}$,
which is the set of indices of the intervals that contain
some point(s) from $B$ and $[m]:=\bigl\{ 0,1,\ldots,m-1 \bigr\}$.

\section{Transformations Preserving Uniformity on Triangles in $\R^2$}
\label{sec:transformations}
The property \textbf{P9} suggests that in higher dimensions
the arc probability of PCDs based on uniform data is geometry invariant,
i.e., does not depend on the geometry of the support set.
The set $\X_n$ is assumed to be a set of iid uniform random variables on
the convex hull of $\Y_m$; i.e., a random sample from $\U\left(\C_H\left( \Y_m \right)\right)$.
In particular, conditional on $|\X_n \cap T_j|>0$
being fixed, $\X_n \cap T_j$ will
also be a set of iid uniform random variables on $T_j$ for $j\in \{1,2,\ldots,J\}$.
The geometry invariance property will reduce the triangle $T_j$
as much as possible while preserving uniformity and the
probabilities related to PCDs will simplify in notation and calculations.
Below, we present such a transformation that reduces a single triangle
to the standard equilateral triangle $T_e=T\left((0,0),(1,0),\left( 1/2,\sqrt{3}/2 \right)\right)$.

Let $\Y_3 =\{\y_1,\y_2,\y_3\} \subset \R^2$ be three non-collinear points
and $\TY$ be the triangle with vertices $\y_1,\y_2,\y_3$.
Let $X_i \stackrel{iid}{\sim} \; \UT$,
the uniform distribution on $\TY$, for $i=1, \ldots, n$.
The pdf of $\UT$ is $$f(u)=\frac{1}{A(\TY)}\I(u \in \TY),$$
where $A(\cdot)$ is the area functional.

\begin{figure}[ht]
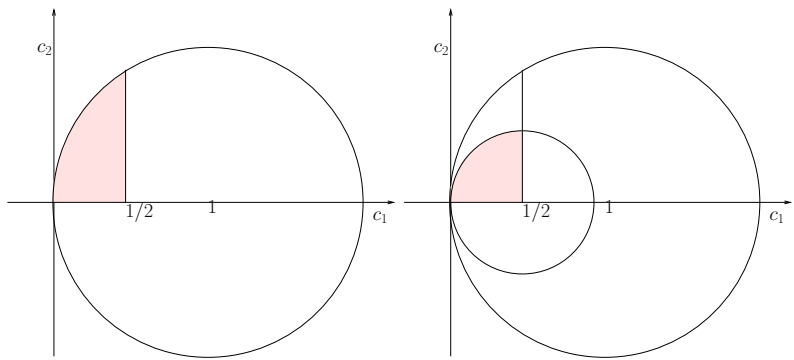

\centering
\scalebox{.3}{\input{domainc1c2.pstex_t}}
\scalebox{.3}{\input{obtdomc1c2.pstex_t}}
\caption{
The shaded regions are the domain of $(c_1,c_2)$ values for the basic triangle $T_b$ (left)
the domain where $T_b$ is a non-acute triangle (right).
%and the domain where the two centers (circumcenter and orthocenter) are inside the triangle (right).
}
\label{fig:domain-c1c2}
\end{figure}

The triangle $\TY$ can be carried into the first quadrant by a
composition of transformations in such a way that the largest
edge has unit length and lies on the $x$-axis, and the
$x$-coordinate of the vertex nonadjacent to largest edge is less than $1/2$.
We call the resultant triangle the \emph{basic triangle} and denote it as $T_b$
where $T_b=\bigl((0,0),(1,0),(c_1,c_2)\bigr)$ with $0< c_1 \le 1/2$,
and $c_2 > 0$ and $(1-c_1)^2+c_2^2 \le 1$.
We will describe such transformations below:
Let $e_i$ be the edge opposite to the vertex $\y_i$ for $i \in \{1,2,3\}$.
Find the lengths of the edges;
say $e_3$ is of maximum length.
Then scale the triangle so that $e_3$ is of unit length.
Next translate $\y_1$ to $(0,0)$,
and rotate  (if necessary) the triangle so that $\y_2=(1,0)$.
If the $y$-coordinate of $\y_3$ is negative reflect the triangle around the $x$-axis,
then if $x$-coordinate of $\y_3$ is greater than $1/2$,
reflect it around $x=1/2$, then the associated basic triangle $T_b$ is obtained.
So the basic triangle $T_b$ can be obtained by a transformation denoted by $\phi_b$
which is a composition of rigid motion
transformations (namely translation, rotation, and reflection) and scaling.
Hence if $\TY$ is transformed into $T_b$, then $\TY$ is similar to $T_b$
and $\phi_b\left(\TY\right)=T_b$.
Thus the random variables $X_i \stackrel{iid}{\sim}\UT$ transformed
along with $\TY$ in the described fashion by $\phi_b$ satisfy $\phi_b(X_i) \stackrel{iid}{\sim}\U(T_b)$.
So, without loss of generality, we can assume $\TY$ to be the basic triangle,
The functional form of $T_b$ is
$$T_b=\left\{ (x,y) \in \R^2 : y \ge 0;\; y \le (c_2\,x)/c_1;\; y \le c_2\,(1-x)/(1-c_1) \right\}.$$
If $c_1=1/2$ and $c_2=\sqrt{3}/2$, then $T_b$ is an equilateral triangle;
if $c_2 < \sqrt{c_1-c_1^2}$, then $T_b$ is an obtuse triangle;
if $c_2=\sqrt{c_1-c_1^2}$, then $T_b$ is a right triangle;
and if $c_2 > \sqrt{c_1-c_1^2}$, then $T_b$ is an acute triangle.
If $c_2=0$, then the $T_b$ reduces to the unit interval $(0,1)$.
See Figure \ref{fig:domain-c1c2} for the domain of $(c_1,c_2)$ for $T_b$
and the part on which $T_b$ is a non-acute triangle.

\begin{lemma}
\label{lem:scale-inv-arc-prob}
The arc probability $\mu(\NY)$ of the PCD based on $\NY$ for uniform data on $\TY$
is rigid-motion and scale invariant; i.e.,
$\mu(\NY)$ does not change under rigid motion transformations and does not
depend on the scale of the support triangle $\TY$.
\end{lemma}
\noindent \textbf{Proof:}
We have shown that for $X_i \stackrel{iid}{\sim}\UT$,
we have $\phi_b(X_i) \stackrel{iid}{\sim}\U(T_b)$
since $\TY$ is similar to $T_b$.
For uniform data,
the set probabilities are calculated as the ratio
of the area of the set to the total area.
So $P(X \in S \subseteq \TY)=A(S)/A(\TY)$
and $P(\phi_b(X) \in \phi_b(S) \subseteq \phi_b(\TY))=
P(\phi_b(X) \in \phi_b(S) \subseteq T_b)=A(\phi_b(S))/A(T_b)
=[k A(S)]/[k A(\TY)]=A(S)/A(\TY)$ where $k$ is the scaling factor.
Letting $X=X_j$ and $S=\NY(X_i)$, the desired result follows.
$\blacksquare$

Based on Lemma \ref{lem:scale-inv-arc-prob},
for uniform data we can, without loss of generality,
assume $\TY$ to be the basic triangle $T_b$.

\subsection{Transformation of $T_b$ to $T_e$}
\label{sec:transform-into-Te}
There are also transformations that preserve uniformity
of the random variable, but not similarity of the triangles.
We only describe the transformation that maps $\TY$ to the
standard equilateral triangle,
$T_e=T\left((0,0),(1,0),\left( 1/2,\sqrt{3}/2 \right)\right)$
for exploiting the symmetry in calculations using $T_e$.

\begin{figure}[ht]
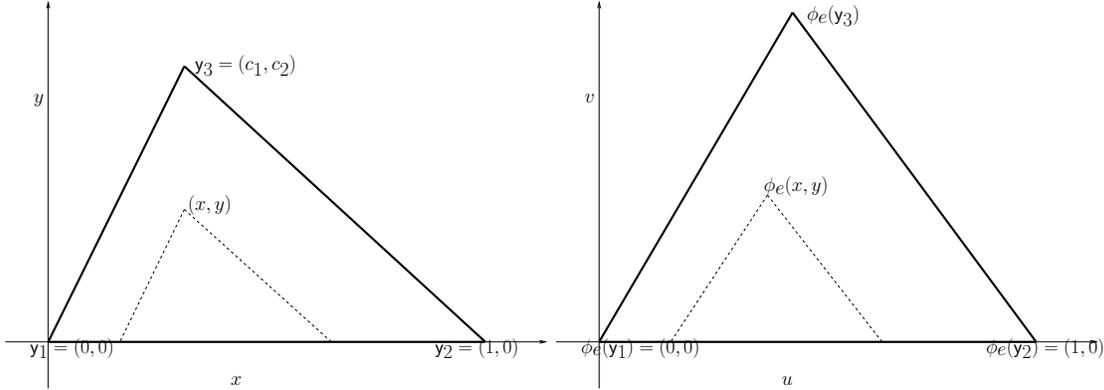

\centering
\scalebox{.3}{\input{Basic_Tri.pstex_t}}
\scalebox{.3}{\input{Equal_Tri_phi.pstex_t}}
\caption{The description of $\phi_e(x,y)$ for $(x,y) \in T_b$ (left) and the equilateral triangle $\phi_e(T_b)=T_e$ (right).}
\label{fig:transform-equal}
\end{figure}

Let $\phi_e:\,(x,y) \rightarrow (u,v)$, where
$\displaystyle u(x,y)=x+\frac{1-2\,c_1}{\sqrt{3}}\,y$ and
$\displaystyle v(x,y)=\frac{\sqrt{3}}{2\,c_2}\,y$.
Then $\y_1$ is mapped to $(0,0)$, $\y_2$ is mapped to $(1,0)$,
and $\y_3$ is mapped to $\left( 1/2,\sqrt{3}/2 \right)$.
See also Figure \ref{fig:transform-equal}.
Note that the inverse transformation is
$\phi_e^{-1}(u,v)=\bigl(x(u,v),y(u,v)\bigr)$ where
$\displaystyle x(u,v)=u-\frac{(1-2\,c_1)}{\sqrt{3}}\,v$ and
$\displaystyle y(u,v)=\frac{2\,c_2}{\sqrt{3}}\,u$.
Then the Jacobian is given by
\begin{eqnarray*}
J(x,y)&=&\left|\begin{array}{cc}
\frac{\partial{x}}{\partial{u}} & \frac{\partial{x}}{\partial{v}}
\vspace{.1 in}\\
\frac{\partial{y}}{\partial{u}} & \frac{\partial{y}}{\partial{v}}
 \end{array}\right|=\left|\begin{array}{cc}
1 & \frac{2\,c_1-1}{\sqrt{3}}
\vspace{.1 in}\\
0 & \frac{2\,c_2}{\sqrt{3}}
 \end{array}\right|=\frac{2\,c_2}{\sqrt{3}}.
\end{eqnarray*}
So $f_{U,V}(u,v)=f_{X,Y}(\phi_e^{-1}(u,v))\,|J| = \frac{4}{\sqrt{3}}\, \I\bigl((u,v)\in T_e\bigr)$. Hence uniformity is preserved.

%\clearpage

\begin{theorem}
\label{thm:geo-inv-arc-prob}
The arc probability $\mu(\NY)$ of the PCD based on $\NY$ for uniform data on $T_b$
is geometry invariant iff $A(\phi_e(\NY(x)))=A(N_{\phi_e\left(\Y_3\right)}(\phi_e(x)))$ for all $x \in T_b$.
\end{theorem}
\noindent \textbf{Proof:}
By Lemma \ref{lem:scale-inv-arc-prob},
the PCD based on $\NY$ for uniform data on $\TY$ is rigid-motion
and scale invariant.
So $\TY$ can be transformed to $T_b$ preserving the uniformity of the data
and the arc probability for the associated PCD.
For uniform data,
the set probabilities are calculated as the ratio
of the area of the set to the total area.
Suppose the arc probability is geometry invariant.
Then $\mu(\NY)=P(X \in \NY(x))=P(\phi_e(X) \in N_{\phi_e\left(\Y_3\right)}(\phi_e(x)))$.
But $P(X \in \NY(x))=A(\NY(x))/A(T_b)$ and
$P(\phi_e(X) \in N_{\phi_e\left(\Y_3\right)}(\phi_e(x)))=A(N_{\phi_e\left(\Y_3\right)}(\phi_e(x)))/A(T_e)$.
Moreover $A(\NY(x))/A(T_b)=A(\phi_e(\NY(x)))/A(\phi_e(T_b))=A(\phi_e(\NY(x)))/A(T_e)$
since the Jacobian cancels out and $\phi_e(T_b)=T_e$.
Hence $A(N_{\phi_e\left(\Y_3\right)}(\phi_e(x)))/A(T_e)=A(\phi_e(\NY(x)))/A(T_e)$
implies $A(\phi_e(\NY(x)))=A(N_{\phi_e\left(\Y_3\right)}(\phi_e(x)))$ for all $x \in T_b$.
The converse can be proved similarly.
%Hence the geometry invariance for arc probability follows
%if the area ratio $A(\NY(x))/A(T_b)$ for all $x \in T_b$ is preserved under $\phi_e$.
%But $A(\NY(x))/A(T_b)=A(\phi_e(\NY(x)))/A(\phi_e(T_b))=A(\phi_e(\NY(x)))/A(T_e)$
%since the Jacobian cancels out.
%The arc probability for the transformed data as $\phi_e(X)$ is calculated with the area ratio
%as $A(N_{\phi_e\left(\Y_3\right)}(\phi_e(x)))/A(T_e)$.
%But these two ratios (for the original and the transformed data)
%are equal iff $A(\phi_e(\NY(x)))=A(N_{\phi_e\left(\Y_3\right)}(\phi_e(x)))$.
$\blacksquare$

\begin{corollary}
\label{cor:geo-inv-arc-prob}
If $\phi_e(\NY(x))=N_{\phi_e\left(\Y_3\right)}(\phi_e(x))$ for all $x \in T_b$,
then the arc probability $\mu(\NY)$ of the PCD based on $\NY$ for uniform data on $T_b$
is geometry invariant.
\end{corollary}
\noindent \textbf{Proof:}
Let $x \in T_b$.
Then $\phi_e(\NY(x))=N_{\phi_e\left(\Y_3\right)}(\phi_e(x))$ implies
$A(\phi_e(\NY(x)))=A(N_{\phi_e\left(\Y_3\right)}(\phi_e(x)))$.
Hence the result follows by Theorem \ref{thm:geo-inv-arc-prob}.
$\blacksquare$

\begin{definition}
\label{def:geo-inv-arc-prob}
The $M$-edge regions are said to be geometry invariant if
$\phi_e(R_M(e_i))=R_{\phi_e(M)}(\phi_e(e_i))$ for $i=1,2,3$.
The $M$-vertex regions are said to be geometry invariant if
$\phi_e(R_M(\y_i))=R_{\phi_e(M)}(\phi_e(\y_i))$ for $i=1,2,3$.
$\square$
\end{definition}
%The geometry invariance for edge and vertex regions
%follow by the preservation of the area ratios under $\phi_e$
%as in Theorem \ref{thm:geo-inv-arc-prob}.

\begin{corollary}
\label{cor:geo-inv-arc-prob-1}
Suppose $\NY$ is based on geometry invariant edge or vertex regions.
If the proximity regions are based on boundary of $\TY$
and parallel lines to edges,
then geometry invariance of the arc probability for uniform data follows.
\end{corollary}
\noindent \textbf{Proof:}
Such proximity maps with geometry invariant edge or vertex regions,
satisfy $\phi_e(\NY(x))=N_{\phi_e\left(\Y_3\right)}(\phi_e(x))$.
Hence the desired result follows by Corollary \ref{cor:geo-inv-arc-prob}.
$\blacksquare$

\begin{corollary}
\label{cor:geo-inv-arc-prob-2}
If the edge or vertex regions are based on specific angles in $T_b$
in the sense that their vertices have specific angular values,
then these regions are not geometry invariant.
Similarly if the proximity regions are based on specific angles in $T_b$
then they are not geometry invariant either.
\end{corollary}
\noindent \textbf{Proof:}
The transformation $\phi_e$ clearly does not preserve the angles in $T_b$.
Hence the regions dependent on (inner) angles of $T_b$ fail to be preserved.
$\blacksquare$

\section{Triangle Centers}
\label{sec:triangle-centers}
The PCDs will be defined using the vertex and edge regions, which will
be constructed using a point, preferably, in the interior of the triangle,
e.g., a \emph{triangle center}.
Let $\Y_3=\{\y_1,\y_2,\y_3\} \subset \R^2$ be non-collinear and $\TY$ be
the corresponding triangle.
The \emph{trilinear coordinates} of a point
$P$ with respect to $\TY$ are an ordered triple of numbers,
which are proportional to the distances from $P$ to the edges.
Trilinear coordinates are denoted as $(\alpha:\beta:\gamma)$ and also are known as
\emph{homogeneous coordinates} or \emph{trilinears}.
Trilinear coordinates were introduced by Pl$\ddot{\text{u}}$cker in 1835
(see \cite{weisstein.1}).
The triplet of trilinear coordinates obtained by multiplying
a given triplet by any positive constant $k$ describes the same point;
i.e., $(\alpha:\beta:\gamma)=(k\alpha:k\beta:k\gamma)$, for any $k  > 0$.
By convention, the three vertices $\y_1,\; \y_2,$ and $\y_3$ of
$\TY$ are commonly written as $(1:0:0)$, $(0:1:0)$, and $(0:0:1)$,
respectively (see \cite{weisstein.1}).

\begin{definition}
A \emph{triangle center} is a point whose trilinear coordinates
are defined in terms of the edge lengths and (inner) angles of a triangle.
The function giving the coordinates $(\alpha:\beta:\gamma)$ is
called the \emph{triangle center function}. $\square$
\end{definition}
\cite{kimberling.1} enumerates 360 triangle centers,
among which four have been widely known since the ancient times; namely,
\emph{circumcenter} ($CC$), \emph{incenter} ($IC$), \emph{center of mass} or \emph{centroid} ($CM$),
and \emph{orthocenter} ($OC$).
The point where the center is located in $\TY$
will be labeled accordingly; e.g., $M_{CC}$ will denote the circumcenter of $\TY$.

\begin{figure}[ht]
\begin{center}
\psfrag{A}{\scriptsize{$\y_1$}}
\psfrag{B}{\scriptsize{$\y_2$}}
\psfrag{C}{\scriptsize{$\y_3$}}
\psfrag{CC}{\scriptsize{$M_{CC}$}}
\psfrag{r}{\scriptsize{$r_{cc}$}}
\psfrag{M1}{\scriptsize{$M_1$}}
\psfrag{M2}{\scriptsize{$M_2$}}
\psfrag{M3}{\scriptsize{$M_3$}}
\psfrag{x}{}
\psfrag{D}{}
\psfrag{E}{}
\epsfig{figure=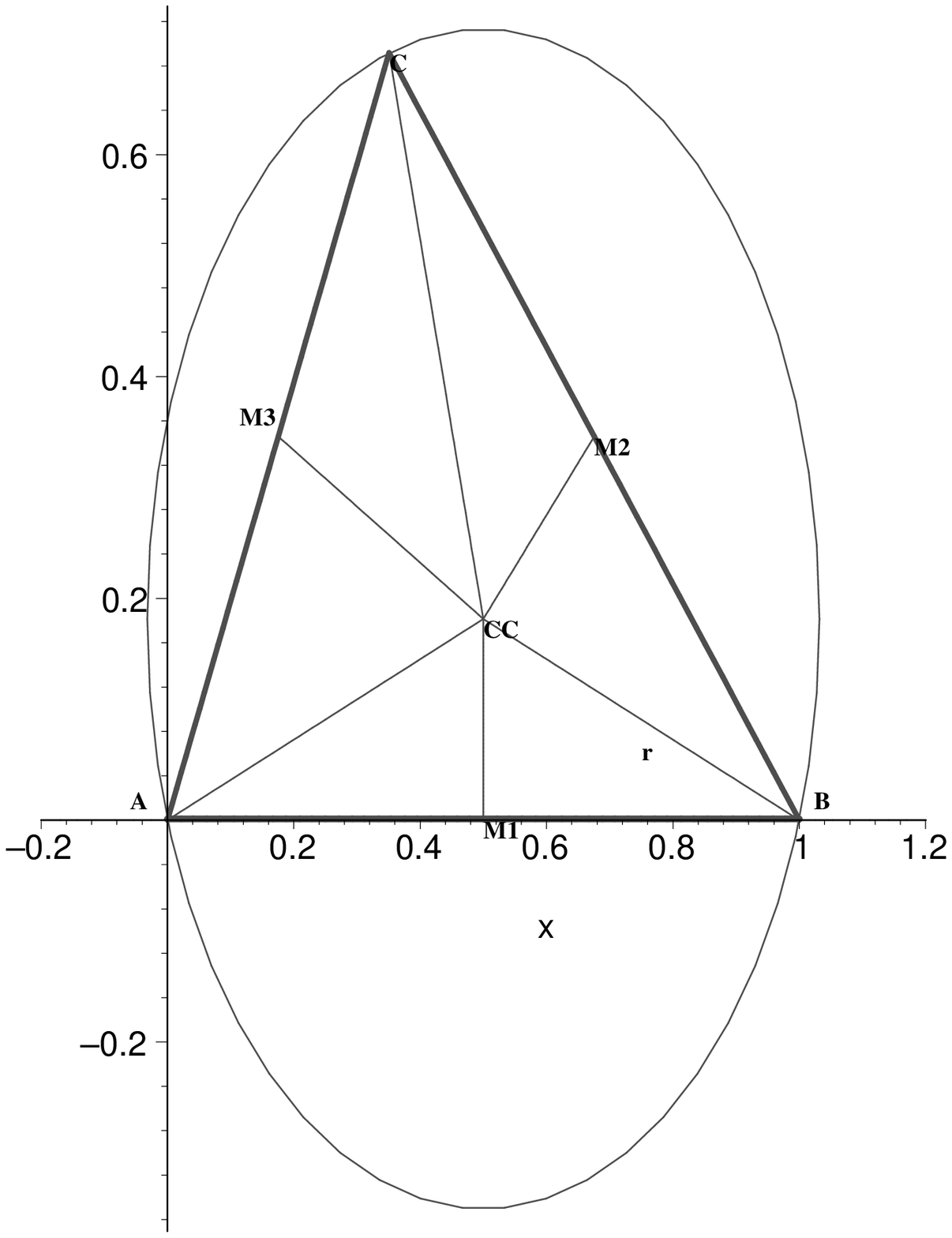, height=150pt , width=200pt}
\epsfig{figure=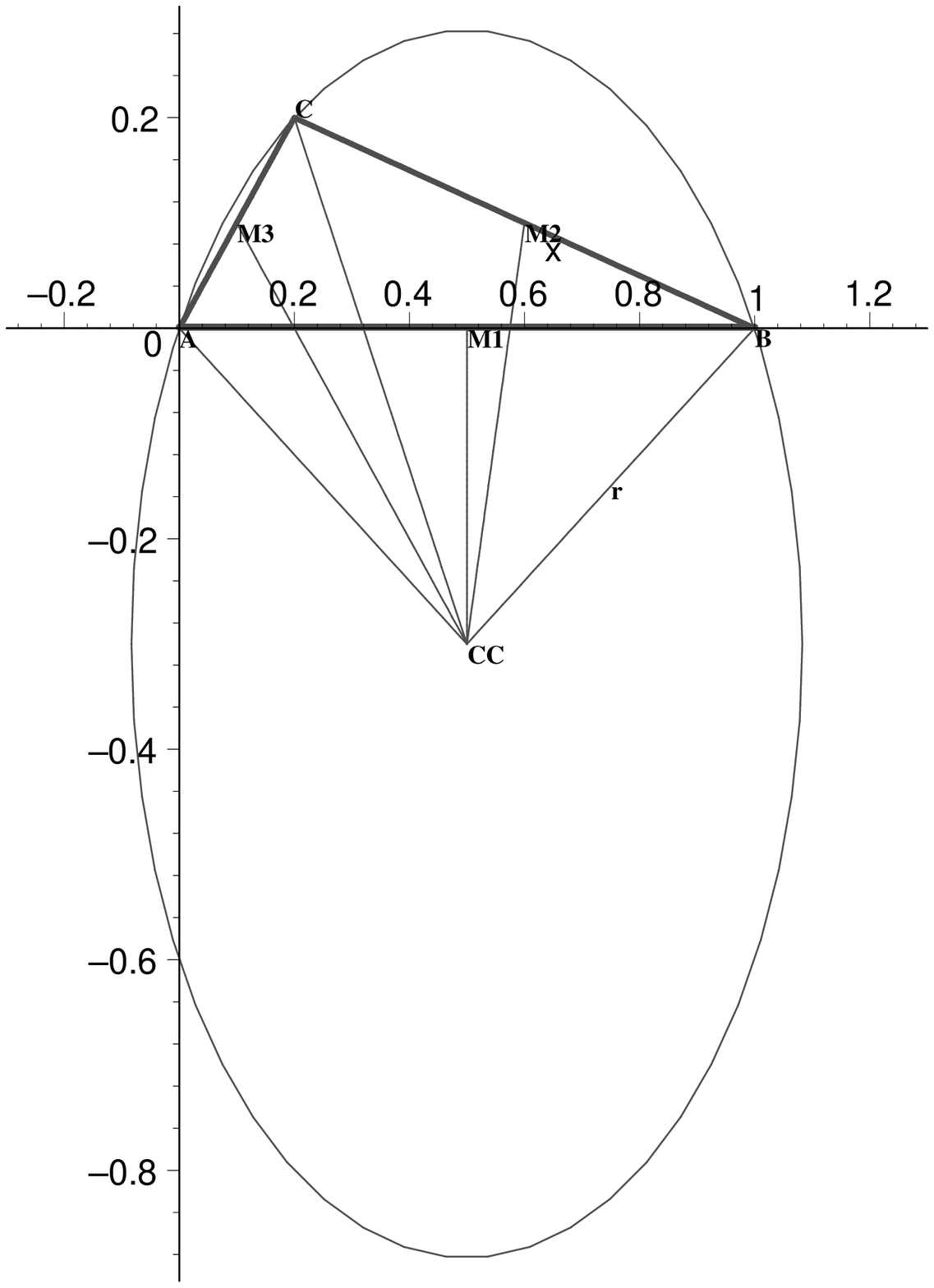, height=150pt , width=200pt}
\end{center}
\caption{The circumcenter, circumcircle, and circumradius of an acute triangle (left),
an obtuse triangle (right).}
\label{fig:circumcenter}
\end{figure}

The \emph{circumcircle} is a triangle's circumscribed circle; i.e.,
the unique circle that passes through each of the triangle's three
vertices, $\y_1,\, \y_2,\,\y_3$.
The center of the circumcircle is
called the \emph{circumcenter}, denoted as $M_{CC}$, and the circle's
radius is called the \emph{circumradius}, denoted as $r_{cc}$.
By construction, the distances from circumcenter to the vertices are equal (to $r_{cc}$).
Furthermore, the triangle's three edge bisectors
perpendicular to edges $e_i$ at $M_i$ for $i \in \{1,2,3\}$
intersect at ${M_{CC}}$.
See Figure \ref{fig:circumcenter}.
The trilinear coordinates of $M_{CC}$ are
$\left( \cos \theta_1:\cos \theta_2:\cos \theta_3 \right)$ where
$\theta_i$ is the inner angle of $\TY$ at vertex $\y_i$ for $i \in \{1,2,3\}$
and the trilinears for $M_{CC}$ can also be written as
$(r_{cc} \cos \theta_1: r_{cc} \cos \theta_2: r_{cc} \cos \theta_3)$.

The circumcenter of a triangle is in the interior,
at the midpoint of the hypotenuse, or in the exterior of
the triangle, if the triangle is acute, right,
or obtuse, respectively. See
Figure \ref{fig:circumcenter} where an acute and an obtuse triangle are depicted.
Using the pdf  of
an arbitrary angle of a triangle $T_j$ from Poisson Delaunay triangulation $\D_P$ (\cite{mardia:1977}),
we see that,
$$P\left( \text{$T_j$ is a right triangle} \right)=P(\theta=\pi/2)=0, $$
hence $P(M_{CC} \text{ is the midpoint of the hypothenuse})=0$.
Furthermore,
\begin{eqnarray}
\label{eqn:prob-obtuse}
P\left( \text{$T_j$ is an obtuse triangle} \right)
&=&P(M_{CC} \not\in T_j)=P(\theta_{\max} > \pi/2)=\int_{\pi/2}^{\pi}f_3(x)dx \nonumber\\
&=& \frac{\left(3\,f_S\left( \sqrt{2\,\pi } \right) -f_C \left( \sqrt{2\,\pi } \right)
-3\,f_S\left( \sqrt{\frac{\pi}{2}} \right) +f_C\left( \sqrt{\frac{\pi}{2}}\right) \right) }{\sqrt{2\,\pi }}\nonumber\\
&\approx& .03726 \nonumber
\end{eqnarray}
where
\begin{multline*}
f_3(x)=\left[\frac{2}{\pi}(3x (\sin 2x)-\cos 2x + \cos 4x -\pi \sin 2x)\right]\,
\I(\pi/3<x<\pi/2))+ \\
\left[\frac{1}{\pi}(4 \pi (\cos x) (\sin x) +
3 \sin x^2- \cos x^2-4x(\cos x)(\sin x) +1)\right]\, \I(\pi/2<x<\pi)
\end{multline*}
is the pdf of the maximum angle,
$\displaystyle f_C(x)=\int _{0}^{x}\!\cos (\pi \,t^2/2) {dt}$, and
$\displaystyle f_S(x)=\int _{0}^{x}\!\sin (\pi \,t^2/2) {dt}$ are
the Fresnel cosine and sine functions, respectively.
The coordinates of $M_{CC}$ in the basic triangle
$T_b$ are $\displaystyle \left( \frac{1}{2},\frac{c_1^2-c_1+c_2^2}{2\,c_2} \right)$.

\begin{figure}[ht]
\begin{center}
\psfrag{A}{\scriptsize{$\y_1$}}
\psfrag{B}{\scriptsize{$\y_2$}}
\psfrag{C}{\scriptsize{$\y_3$}}
\psfrag{IC}{\scriptsize{$M_{I}$}}
\psfrag{r}{\scriptsize{$r_{ic}$}}
\psfrag{P1}{\scriptsize{$P_1$}}
\psfrag{P2}{\scriptsize{$P_2$}}
\psfrag{P3}{\scriptsize{$P_3$}}
\psfrag{x}{}
\epsfig{figure=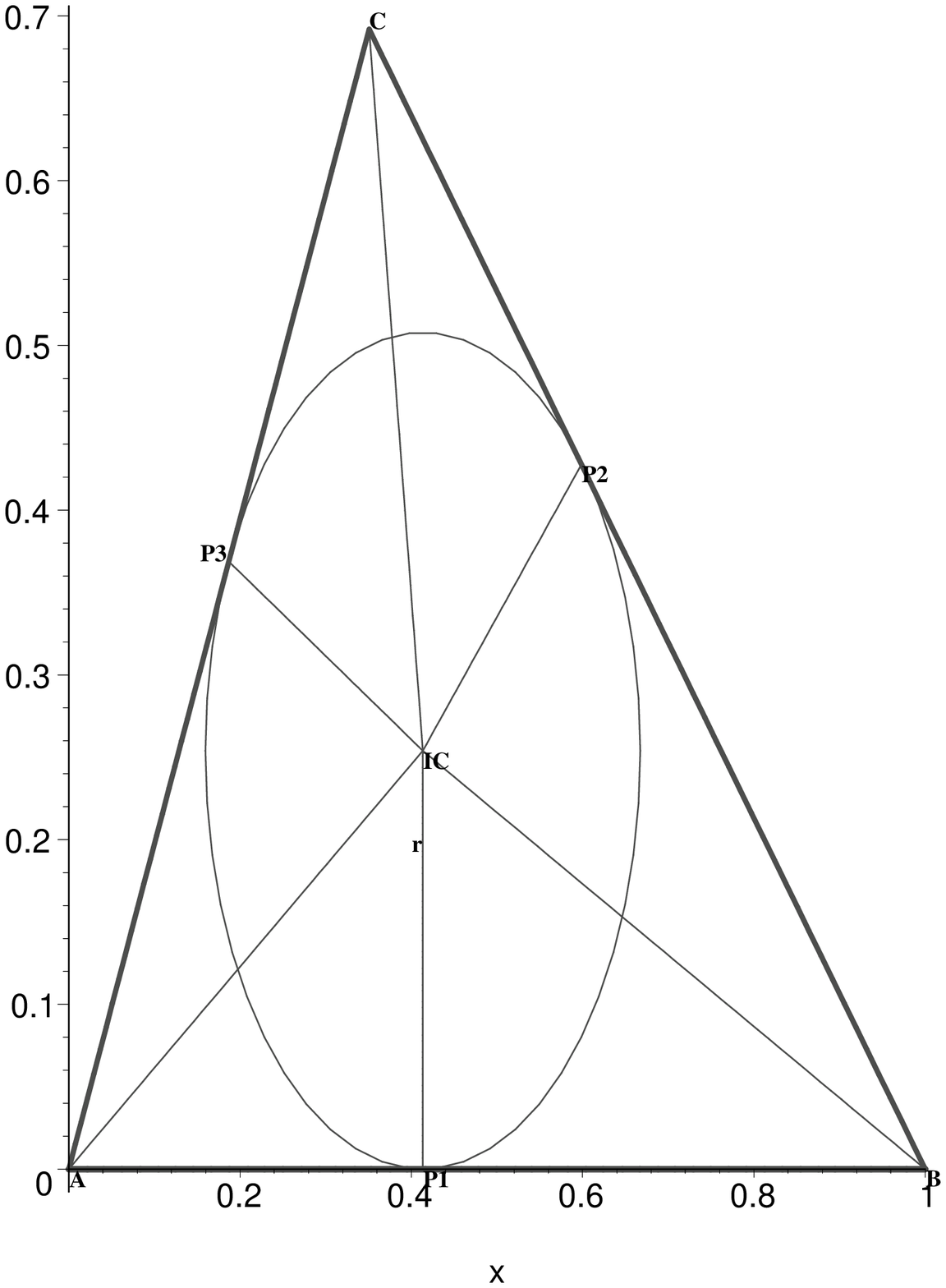, height=150pt, width=200pt}
\psfrag{CM}{\scriptsize{$M_C$}}
\psfrag{r}{\scriptsize{$r_{ic}$}}
\psfrag{M1}{\scriptsize{$M_1$}}
\psfrag{M2}{\scriptsize{$M_2$}}
\psfrag{M3}{\scriptsize{$M_3$}}
\epsfig{figure=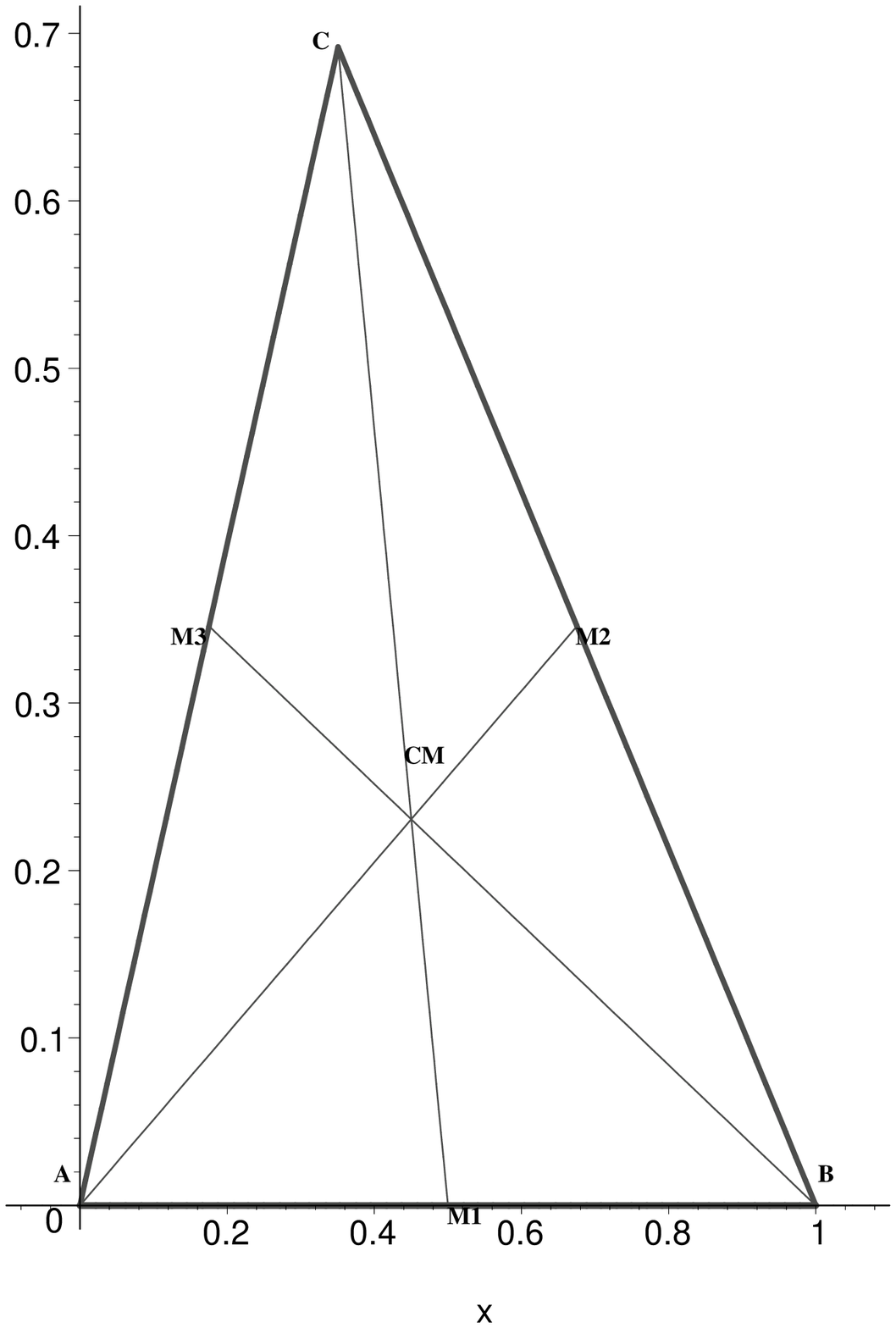, height=150pt, width=200pt}
\end{center}
\caption{The incircle, incenter, inradius of a triangle (left),
and the centroid or center of mass of a triangle (right).}
\label{fig:in-masscenter}
\end{figure}

The \emph{incircle} is the inscribed circle of a triangle, i.e.,
the unique circle that is tangent to the edges $e_i$ at $P_i$ for $i \in \{1,2,3\}$.
The center of the incircle is called the \emph{incenter}, denoted as $M_I$,
and the radius of the incircle is called the \emph{inradius}, denoted as $r_{ic}$.
Incenter has trilinear coordinates $(1:1:1)$.
The incenter is the point where the triangle's inner angle bisectors meet.
See Figure \ref{fig:in-masscenter} (left).

The coordinates of $M_I$ for the basic triangle $T_b$ are $(x_{ic},y_{ic})$, where
%\begin{align*}
%x_{ic}&=\tan\left(\frac{1}{2}\,\arctan\left({\frac{c_2}{1-{c_1}}}\right)\right)
%\left(\tan\left(\frac{1}{2}\,\arctan\left({\frac{c_2}{c_1}}\right)\right)+
%\tan\left(\frac{1}{2}\,\arctan\left({\frac{c_2}{1-{c_1}}}\right)\right)\right)^{-1},\\
%y_{ic}&=\tan\left(\frac{1}{2}\,\arctan\left({\frac{c_2}{c_1}}\right)\right)
%\tan\left(\frac{1}{2}\,\arctan\left({\frac{c_2}{1-{c_1}}}\right)\right)
%\left(\tan\left(\frac{1}{2}\,\arctan\left({\frac{c_2}{c_1}}\right)\right)+
%\tan\left(\frac{1}{2}\,\arctan\left({\frac{c_2}{1-{c_1}}}\right)\right)\right)^{-1},
%\end{align*}
%or
$$x_{ic}={\frac{c_1-\sqrt{c_1^2+c_2^2}}{1+\sqrt{c_1^2+c_2^2}+\sqrt{(1-c_1)^2+c_2^2}}}, \;\; y_{ic}={\frac{c_2}{1+\sqrt{c_1^2+c_2^2}+\sqrt{(1-c_1)^2+c_2^2}}}.$$
%We will use the latter more frequently for simplicity.

Note that, $M_{CC}$ and $M_I$ do not necessarily concur.
The distance between $M_{CC}$ and $M_I$ is $d(M_{CC},M_I)=\sqrt{r_{cc}(r_{cc}-2\,r_{ic})}$.
Unlike the circumcenter, the incenter is guaranteed to be inside the triangle.

%\begin{definition}
The \emph{median line} of a triangle is the line from one of its vertices to
the midpoint of the opposite edge.
%\end{definition}
The three median lines of any triangle intersect at the triangle's
\emph{centroid}, denoted as $M_C$.
The centroid is the \emph{center of mass}
of the vertices of a triangle.
Since $M_C$ is also the intersection of the
triangle's three median lines, it is sometimes called the \emph{median point}.
It has trilinear coordinates $\left( 1/|e_1|:1/|e_2|:1/|e_3| \right)$ or
$\left( \csc \theta_1:\csc \theta_1: \csc \theta_1 \right)$ where $e_i$ denotes
the edge opposite to the vertex $\y_i$ for $i\in \{1,2,3\}$.
The centroid is also
guaranteed to be in the
interior of the triangle.  See Figure \ref{fig:in-masscenter} (right).
The coordinates of $M_C$ for the basic triangle are
$\bigl( (1+c_1)/3,c_2/3 \bigr)$.

The intersection of the three altitudes of a triangle is called the
\emph{orthocenter}, $M_O$, which has trilinear coordinates
$\left( \cos \theta_2 \cos \theta_3:\cos \theta_1 \cos \theta_3: \cos \theta_1 \cos \theta_2 \right)$.
The orthocenter of a triangle is in the interior, at vertex $\y_3$,
or in the exterior of the basic triangle, $T_b$,
if $T_b$ is acute, right, or obtuse, respectively.
The functional form of $M_O$ in the basic triangle is
$\left( c_1,c_1\,(1-c_1)/c_2 \right)$.

Note that in an equilateral triangle, $M_I=M_{CC}=M_O=M_C$
(i.e., all four centers we have described coincide).

\section{Vertex and Edge Regions}
\label{sec:vertex-edge-regions}
The new proximity maps will be based on the Delaunay cell
$\T_j$ that contains $x$.
The region $\NY(x)$ will also depend on the location of $x$
in $\T_j$ with respect to
the vertices or faces (edges in $\R^2$) of $\T_j$.
Hence for $\NY(x)$ to be well-defined, the vertex or face of
$\T_j$ associated with $x$
should be uniquely determined.
This will give rise to two new concepts: \emph{vertex regions} and
\emph{face regions} (\emph{edge regions} in $\R^2$).

\subsection{Vertex Regions}
\label{sec:vertex-regions}
Let $\Y_3=\{\y_1,\y_2,\y_3\}$ be three non-collinear points
in $\R^2$ and $\TY=T(\y_1,\y_2,\y_3)$ be the triangle
with vertices $\Y_3$.
Then for $x \in \TY$,
$N_S(x)=B(x,r(x))$ where $r(x)=\min_{\y \in \Y_3}d(x,\y)$.
That is, $r(x)=d(x,\y_i)$ iff $x \in \mathcal V_C(\y_i) \cap \TY$
for $i \in \{1,2,3\}$, where $\mathcal V_C(\y_i)$ is the Voronoi cell
generated by $\y_i$ in the Voronoi tessellation based on $\Y_3$.
Notice that these cells partition the triangle $\TY$
and each $\mathcal V_C(\y_i) \cap \TY$ is adjacent only to
vertex $\y_i$ and their intersection is the point $M$ which
is equidistant to the vertices,
so $M$ is in fact the circumcenter, $M_{CC}$, of $\TY$.
See Figure \ref{fig:cc-vertex-regions}.
To define new proximity regions based on some sort of distance
or dissimilarity relative to the vertices $\Y_3$, we
associate each point in $\TY$ to a vertex of $\TY$ as in the spherical case.
This gives rise to the concept of \emph{vertex regions}.
Note that $N_S(x)$ is constructed
using the vertex region based on the closest vertex,
$\argmin_{\y \in \Y_3}d(x,\y)$.
If two vertices were equidistant from $x$
(i.e., $\argmin_{\y \in \Y_3}d(x,\y)$ were not unique),
$x$ is arbitrarily assigned to a vertex region.
In fact, for $N_S$, by construction,
it would not matter which vertex to pick when the vertices are equidistant to $x$,
the region $N_S(x)$ will be the same.

\begin{definition}
The connected regions that partition the triangle, $\TY$,
(in the sense that the intersections of the regions have zero $\R^2$-Lebesgue measure)
 such that each region has one and only one vertex of $\TY$ on its boundary
are called \emph{vertex regions}. $\square$
\end{definition}
This definition implies that there are three vertex regions.
In fact, the vertex regions can be described
starting with a point $M \in \R^2 \setminus \Y_3$.
Join the point $M$ to a point on each edge by a curve such
that the resultant regions satisfy the above definition.
Such regions are called \emph{$M$-vertex regions} and
we denote the vertex region associated with
vertex $\y$ as $R_M(\y)$ for $\y \in \Y_3$.
%Vertex regions can be defined using any point $M \in \R^2 \setminus \Y_3$
%by joining $M$ to a point on each edge.
In particular, one can use a \emph{center} of the triangle $\TY$ as
the starting point $M$ for vertex regions.
See the discussion of triangle centers in Section \ref{sec:triangle-centers}.
The points in $R_M(\y)$ can be thought as
being ``closer" to $\y$ than to the other vertices.

It is reasonable to require that the area of the region $R_M(\y)$ gets larger
as $d(M,\y)$ increases.
Usually the curves will be taken to be
lines or even the orthogonal projections to the edges.
But these lines do not necessarily yield three vertex regions for
$M$ in the exterior of $\TY$.
Unless stated otherwise, $M$-vertex regions will refer to
regions constructed by joining $M$ to the edges
with \emph{straight line segments}.

\subsubsection{$M$-Vertex Regions}
For $M \in \TY^o$, $M$-vertex regions are defined by two ways:

\textbf{(I)} Geometrically, one can construct $M$-vertex regions by drawing the \textbf{orthogonal projections to the edges},
denoted as $R^{\perp}_{M}(\y)$. For instance see Figure \ref{fig:cc-vertex-regions} with $M=M_{CC}$.

The functional forms of $R^{\perp}_{M}(\y)$ for $M=(m_1,m_2)$ in the basic triangle are:
\begin{align*}
R^{\perp}_{M}(\y_1)&=\left \{ (x,y) \in \TY : x \le m_1;\; y \le m_2-(x-m_1)\,c_1/c_2  \right \},\\
R^{\perp}_{M}(\y_2)&=\left \{(x,y) \in \TY : x \ge m_1;\; y \le m_2+(1-c_1)\,(x-m_1)/c_2 \right \}, \text{ and } \\
R^{\perp}_{M}(\y_3)&=\left \{(x,y) \in \TY : y \ge m_2-c_1\,(x-m_1)/c_2 ;\; y \ge m_2+(1-c_1)\,(x-m_1)/c_2 \right\}.
\end{align*}
However, the orthogonal projections from $M$ to the edges does not necessarily fall on the boundary of $\TY$.
For example, letting $P_2^{M}$ be the orthogonal projection of $M$ to edge $e_2$,
it is easy to see that $P_2^{M}$ might fall outside $\TY$ which
contradicts the definition of vertex regions.
In fact $P_2^{M} \in e_2$ iff
$\displaystyle \frac {c_2\,(m_2\,c_2+c_1\,m_1)}{c_1^2+c_2^2} \le c_2$ iff
$c_2\,(c_2-m_2)+c_1(c_1-m_1) \ge 0$.

\textbf{(II)} One can also construct $M$-vertex regions with $M \in \TY^o$ by
\textbf{using the extensions of the line segments joining $\y$ to $M$}
for all $\y \in \Y_3$.
See Figure \ref{fig:cm-vertex-regions} with $M=M_C$.
The functional forms of $R_{M}(\y_i)$ for $i \in \{1,2,3\}$ with
$M=(m_1,m_2)$ and $m_1 > c_1$ in the basic triangle, $T_b$, are given by
\begin{align*}
R_{M}(\y_1)&=\left \{(x,y) \in \TY : y \le \frac{m_2\,(x-1)}{m_1-1};\; y \le \frac{m_2(c_1-x)+c_2(x-m_1)}{c_1-m_1} \right\},\\
R_{M}(\y_2)&=\left \{(x,y) \in \TY : y \le \frac{m_2\,x}{m_1};\; y \ge \frac{m_2(c_1-x)+c_2(x-m_1)}{c_1-m_1} \right\}, \text{ and } \\
R_{M}(\y_3)&=\left \{(x,y) \in \TY : y \ge \frac{m_2\,x}{m_1} ;\; y \ge \frac{m_2\,(x-1)}{m_1-1} \right\}.
\end{align*}
For $m_1<c_1$, $R_M(\y_i)$ for $i \in \{1,2,3\}$ are defined similarly.

If $x$ falls on the boundary of two $M$-vertex regions, then
$x$ is arbitrarily assigned to one of the $M$-vertex regions.

To distinguish between these two types, the vertex regions constructed
by using orthogonal projections are denoted as $R^{\perp}_M(\y)$ and
the vertex regions constructed by using the lines joining vertices to $M$ are denoted as $R_M(\y)$.
By definition, $R_{OC}(\y)$ and $R^{\perp}_{OC}(\y)$ are identical.
But, for $M=M_C,M_{CC},M_I$, $R_{M}(\y)$ can have both versions,
so the above distinction is necessary for them.

\subsubsection{$CC$-Vertex Regions}
\begin{figure}[ht]
\begin{center}
\psfrag{A}{\scriptsize{$\y_1$}}
\psfrag{B}{\scriptsize{$\y_2$}}
\psfrag{C}{\scriptsize{$\y_3$}}
\psfrag{R(A)}{\scriptsize{$R_{CC}(\y_1)$}}
\psfrag{R(B)}{\scriptsize{$R_{CC}(\y_2)$}}
\psfrag{R(C)}{\scriptsize{$R_{CC}(\y_3)$}}
\psfrag{CC}{\scriptsize{$M_{CC}$}}
\psfrag{M1}{\scriptsize{$M_3$}}
\psfrag{M2}{\scriptsize{$M_1$}}
\psfrag{M3}{\scriptsize{$M_2$}}
\psfrag{x}{}
\epsfig{figure=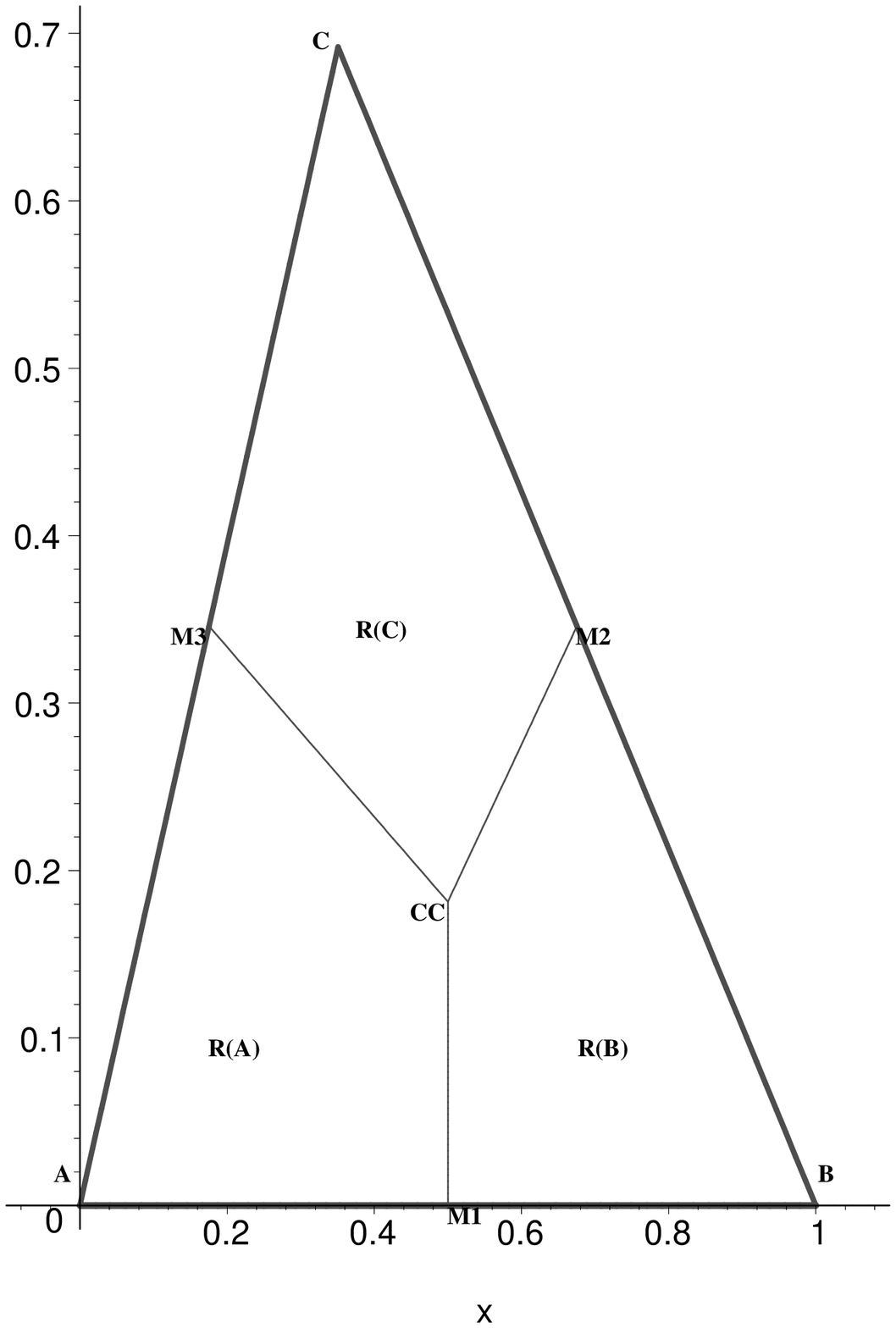, height=150pt, width=200pt}
\epsfig{figure=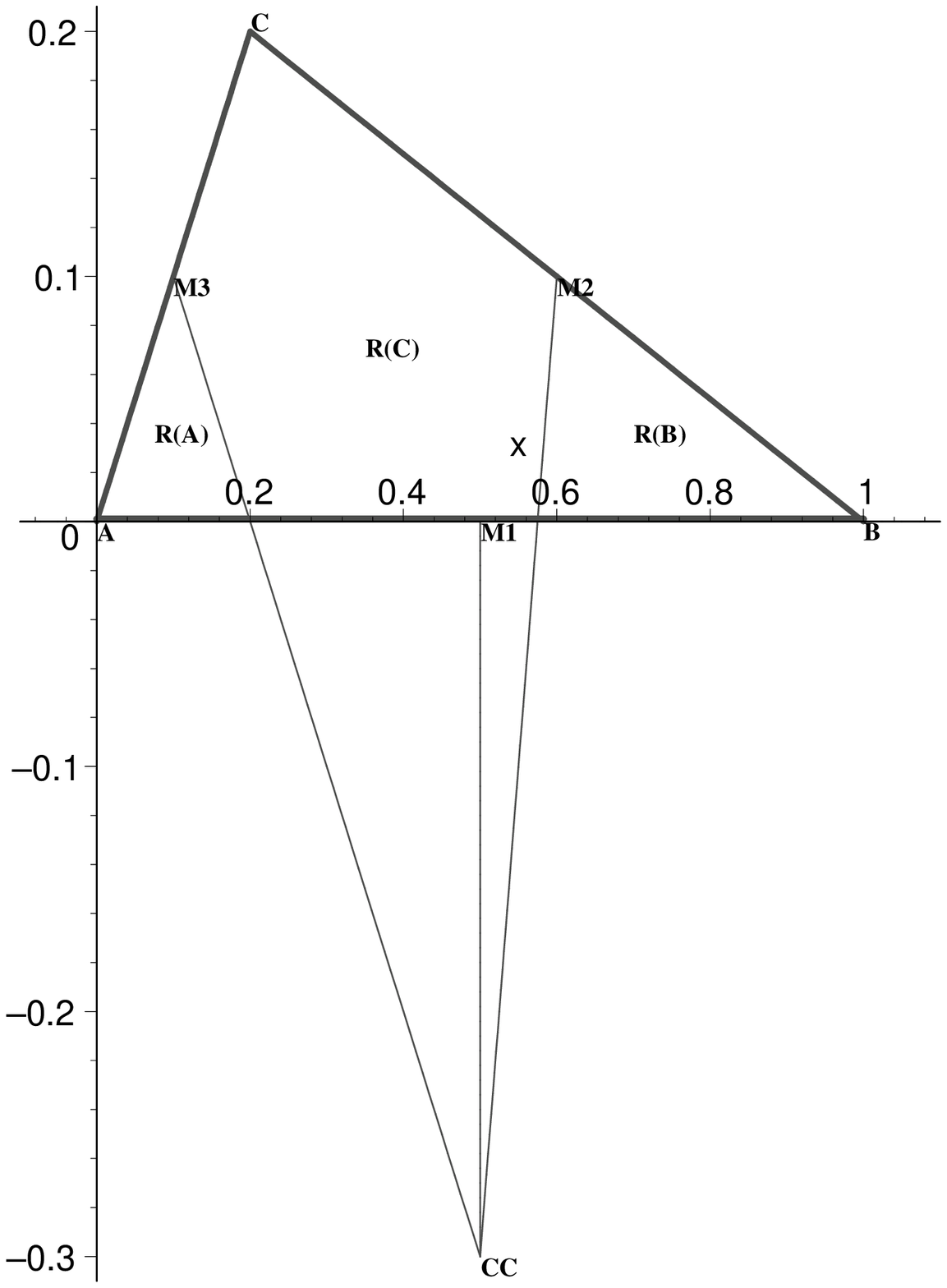, height=100pt, width=200pt}
\end{center}
\caption{
\label{fig:cc-vertex-regions}
The $CC$-vertex regions in an acute triangle (left) and in an obtuse triangle (right).}
\end{figure}

The region $\V_C(\y) \cap \TY$ is a special
type of vertex regions, which can also be obtained
geometrically by starting at
$M_{CC}$ and drawing the orthogonal projections to the edges.
Hence these regions are called \emph{$CC$-vertex regions}.
One can also construct $CC$-vertex regions
by drawing the perpendicular (mid)edge bisectors
or by finding the circumcenter and drawing the orthogonal projections to the edges.
See Figure \ref{fig:cc-vertex-regions},
where $M_i$ are the midpoints of the edges.

The functional forms of $R_{CC}(\y_i)$ for $i \in \{1,2,3\},$ in the basic triangle
$T_b=\bigl( (0,0),(1,0),(c_1,c_2) \bigr)$ (see Section \ref{sec:transformations}) are given by
\begin{align*}
R_{CC}(\y_1)&=\left \{(x,y) \in T\left( \Y_3 \right):
x \le \frac{1}{2};\; y \le \frac{c_1^2+c_2^2-2\,c_1\,x}{2\,c_2} \right\},\\
R_{CC}(\y_2)&=\left \{(x,y) \in T\left( \Y_3 \right):
x \ge \frac{1}{2};\; y \le \frac{c_1^2+c_2^2-1+2\,(1+c_1)\,x}{2\,c_2} \right\}, \\
R_{CC}(\y_3)&=\left \{(x,y) \in T\left( \Y_3 \right):
y \ge \frac{c_1^2+c_2^2-2\,c_1\,x}{2\,c_2};\; y \ge \frac{c_1^2+c_2^2-1+2\,(1+c_1)\,x}{2\,c_2} \right\}.
\end{align*}

One can also define $CC$-vertex regions by using the line segments
which join $M_{CC}$ to edge $e_i$
and are extensions of the lines joining $M_{CC}$ to the vertex $v_i$ for $i \in \{1,2,3\}$,
but this definition only works for acute triangles,
since $M_{CC} \not\in \TY^o$ for non-acute triangles.

%If $x$ falls on the boundary of two $CC$-vertex regions, then
%$x$ is arbitrarily assigned to one of the $CC$-vertex regions.

\subsubsection{$CM$-Vertex Regions}
\begin{figure}[ht]
\centering
\psfrag{c1}{\scriptsize{$c_1$}}
\psfrag{c2}{\scriptsize{$c_2$}}
\epsfig{figure=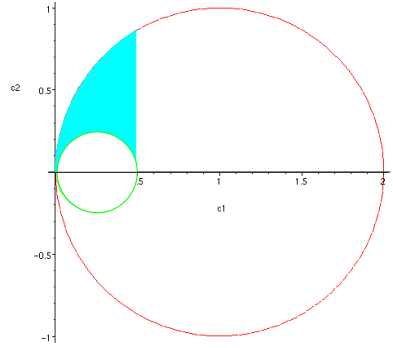, height=145pt , width=145pt}
\caption{
\label{fig:domain-c1c2-P3CMinT}
Depicted is the domain of $(c_1,c_2)$ where $P_2^{CM} \in e_2$.
}
\end{figure}

\begin{figure}[ht]
\begin{center}
\psfrag{A}{\scriptsize{$\y_1$}}
\psfrag{B}{\scriptsize{$\y_2$}}
\psfrag{C}{\scriptsize{$\y_3$}}
\psfrag{R(A)}{\scriptsize{$R^{\perp}_{CM}(\y_1)$}}
\psfrag{R(B)}{\scriptsize{$R^{\perp}_{CM}(\y_2)$}}
\psfrag{R(C)}{\scriptsize{$R^{\perp}_{CM}(\y_3)$}}
\psfrag{CM}{\scriptsize{$M_{C}$}}
\psfrag{P1}{\scriptsize{$P^{CM}_3$}}
\psfrag{P2}{\scriptsize{$P^{CM}_1$}}
\psfrag{P3}{\scriptsize{$P^{CM}_2$}}
\psfrag{x}{}
\epsfig{figure=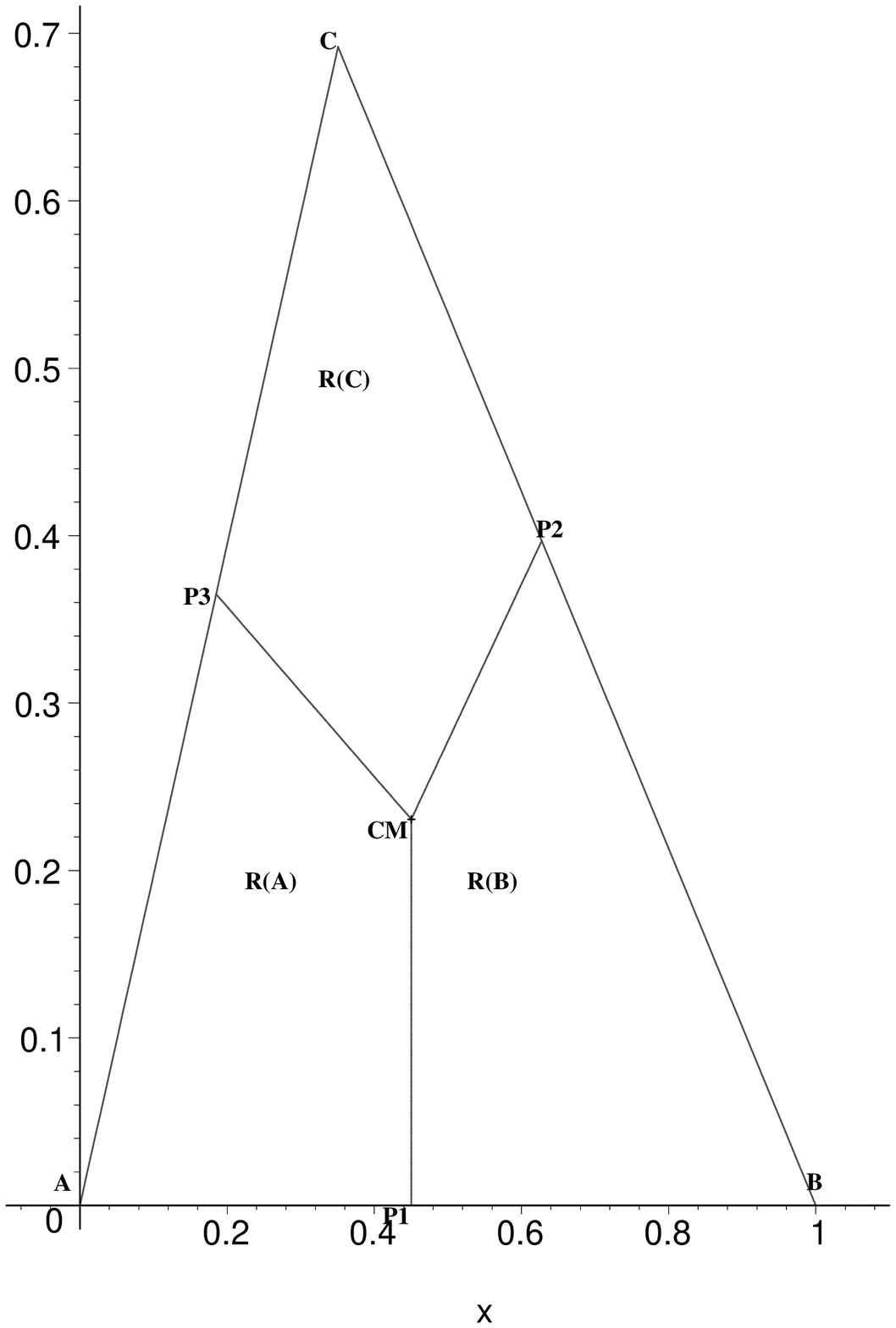, height=140pt , width=200pt}
\psfrag{R(A)}{\scriptsize{$R_{CM}(\y_1)$}}
\psfrag{R(B)}{\scriptsize{$R_{CM}(\y_2)$}}
\psfrag{R(C)}{\scriptsize{$R_{CM}(\y_3)$}}
\psfrag{P1}{$M_3$}
\psfrag{P2}{$M_1$}
\psfrag{P3}{$M_2$}
\epsfig{figure=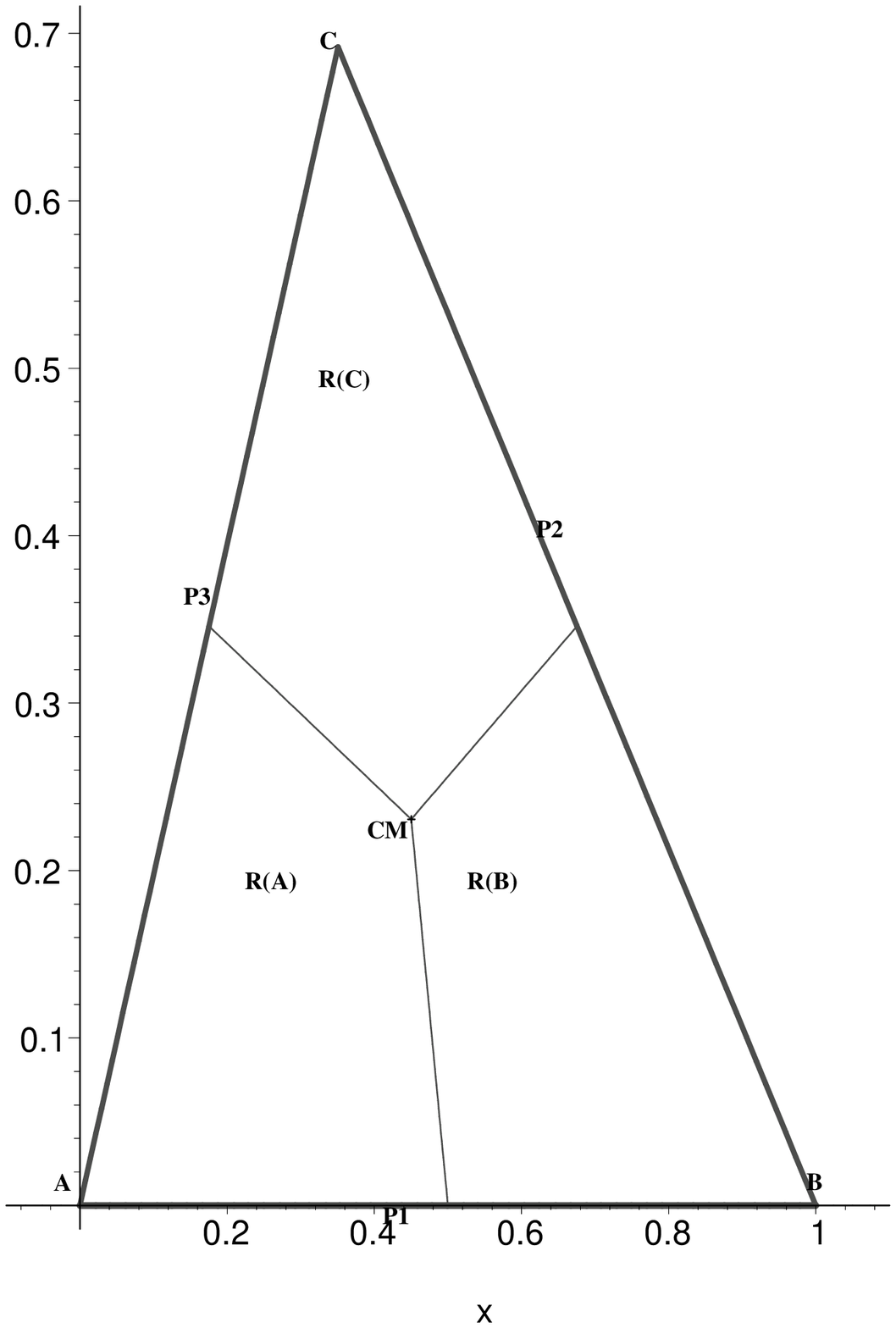, height=140pt , width=200pt}
\end{center}
\caption{
\label{fig:cm-vertex-regions}
The $CM$-vertex regions with orthogonal projections (left) and
with median lines (right).}
\end{figure}

The motivation behind $CM$-vertex regions is that
unlike the circumcenter, center of mass is guaranteed to be
inside the triangle.
We define the $CM$-vertex regions using the median lines and
denote the regions as $R_{CM}(\y_i)$
for $i \in \{1,2,3\}$ (see Figure \ref{fig:cm-vertex-regions}).
\noindent
However, the method with orthogonal projections of $M_C$ to
the edges does not always work.
Let $P_i^{CM}$ be the point at which orthogonal projection
of $M_C$ on $e_i$ crosses $e_i$ for $i \in \{1,2,3\}$.
Then, $P_2^{CM}$ might fall outside $\TY$ in which case
$R^{\perp}_{CM}(\y_1)$ is adjacent to
two vertices $\y_1$ and $\y_3$,
while $R^{\perp}_{CM}(\y_3)$ is not adjacent to any of the vertices.
Hence the definition of the vertex regions is violated.
In fact $P_2^{CM} \in e_2$ iff $\displaystyle \frac{c_2\,(c_2^2+c_1^2+c_1)}{3\,(c_1^2+c_2^2)} \le c_2$
iff $2\,c_1^2+2\,c_2^2-c_1 \ge 0$.
See Figure \ref{fig:domain-c1c2-P3CMinT} for the domain of $(c_1,c_2)$
in $T_b$ for $P_2^{CM} \in e_2$.

The functional forms of $R_{CM}(\y_i)$ for $i \in \{1,2,3\}$
in the basic triangle $T_b$ are given by
\begin{align*}
R_{CM}(\y_1)&=\left \{(x,y) \in T\left( \Y_3 \right):
y \le \frac{c_2\,(2\,x-1)}{2\,c_1-1};\; y \le \frac{c_2\,(x-1)}{c_1-2} \right\},\\
R_{CM}(\y_2)&=\left \{(x,y) \in T\left( \Y_3 \right):
y \le \frac{c_2\,(2\,x-1)}{2\,c_1-1};\; y \le \frac{c_2\,x}{1+c_1} \right\}, \\
R_{CM}(\y_3)&=\left \{(x,y) \in T\left( \Y_3 \right):
y \ge \frac{c_2\,(x-1)}{c_1-2};\; y \ge \frac{c_2\,x}{1+c_1} \right\}.
\end{align*}

%If $x$ falls on the boundary of two $CM$-vertex regions, then
%$x$ is arbitrarily assigned to one of the $CM$-vertex regions.

\subsubsection{$IC$-Vertex Regions}
\begin{figure}[ht]
\begin{center}
\psfrag{A}{\scriptsize{$\y_1$}}
\psfrag{B}{\scriptsize{$\y_2$}}
\psfrag{C}{\scriptsize{$\y_3$}}
\psfrag{R(A)}{\scriptsize{$R^{\perp}_{IC}(\y_1)$}}
\psfrag{R(B)}{\scriptsize{$R^{\perp}_{IC}(\y_2)$}}
\psfrag{R(C)}{\scriptsize{$R^{\perp}_{IC}(\y_3)$}}
\psfrag{IC}{\scriptsize{$M_{IC}$}}
\psfrag{P1}{\scriptsize{$P^{IC}_3$}}
\psfrag{P2}{\scriptsize{$P^{IC}_1$}}
\psfrag{P3}{\scriptsize{$P^{IC}_2$}}
\psfrag{x}{}
\epsfig{figure=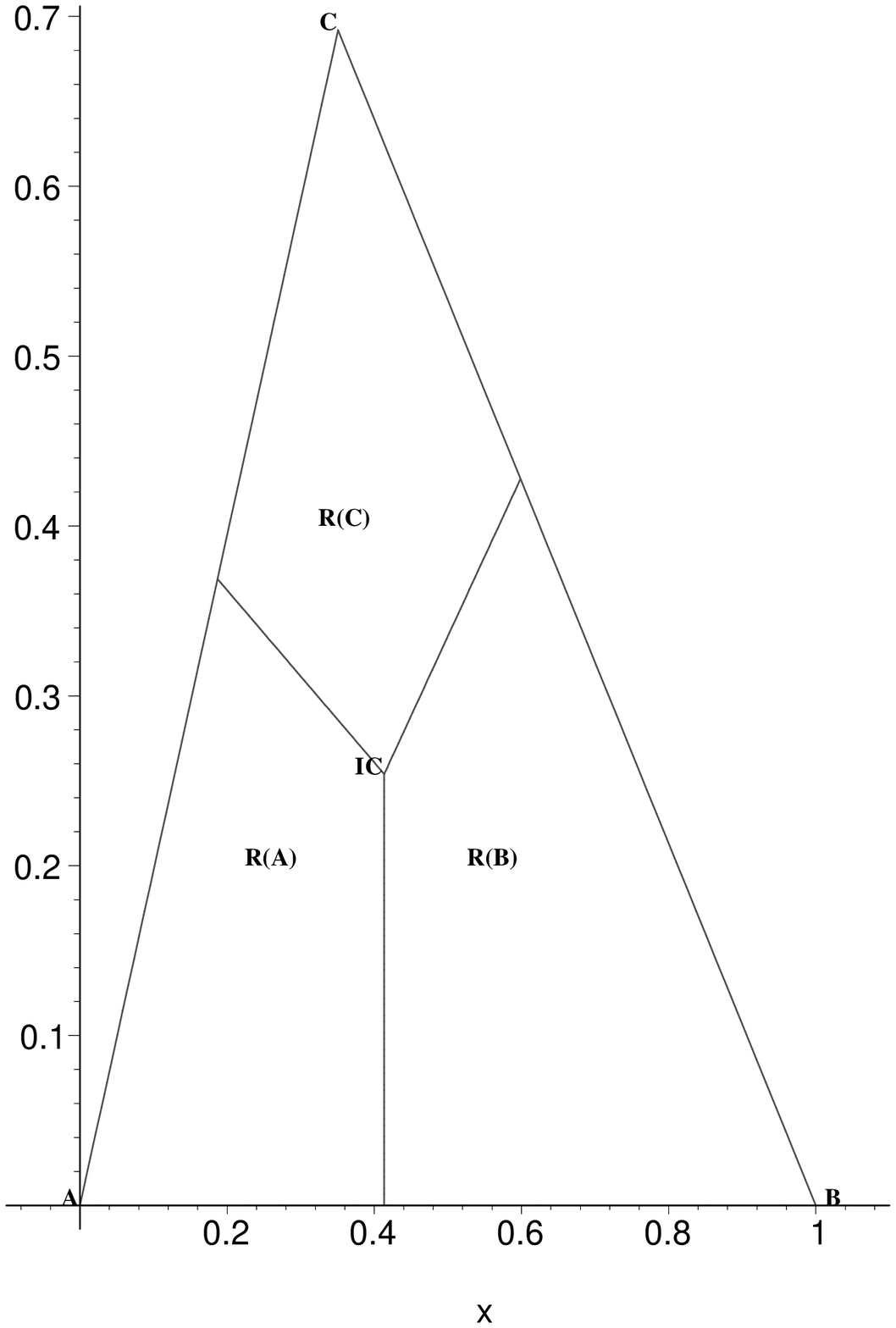, height=140pt , width=200pt}
\psfrag{R(A)}{\scriptsize{$R_{IC}(\y_1)$}}
\psfrag{R(B)}{\scriptsize{$R_{IC}(\y_2)$}}
\psfrag{R(C)}{\scriptsize{$R_{IC}(\y_3)$}}
\psfrag{P1}{$M_3$}
\psfrag{P2}{$M_1$}
\psfrag{P3}{$M_2$}
\epsfig{figure=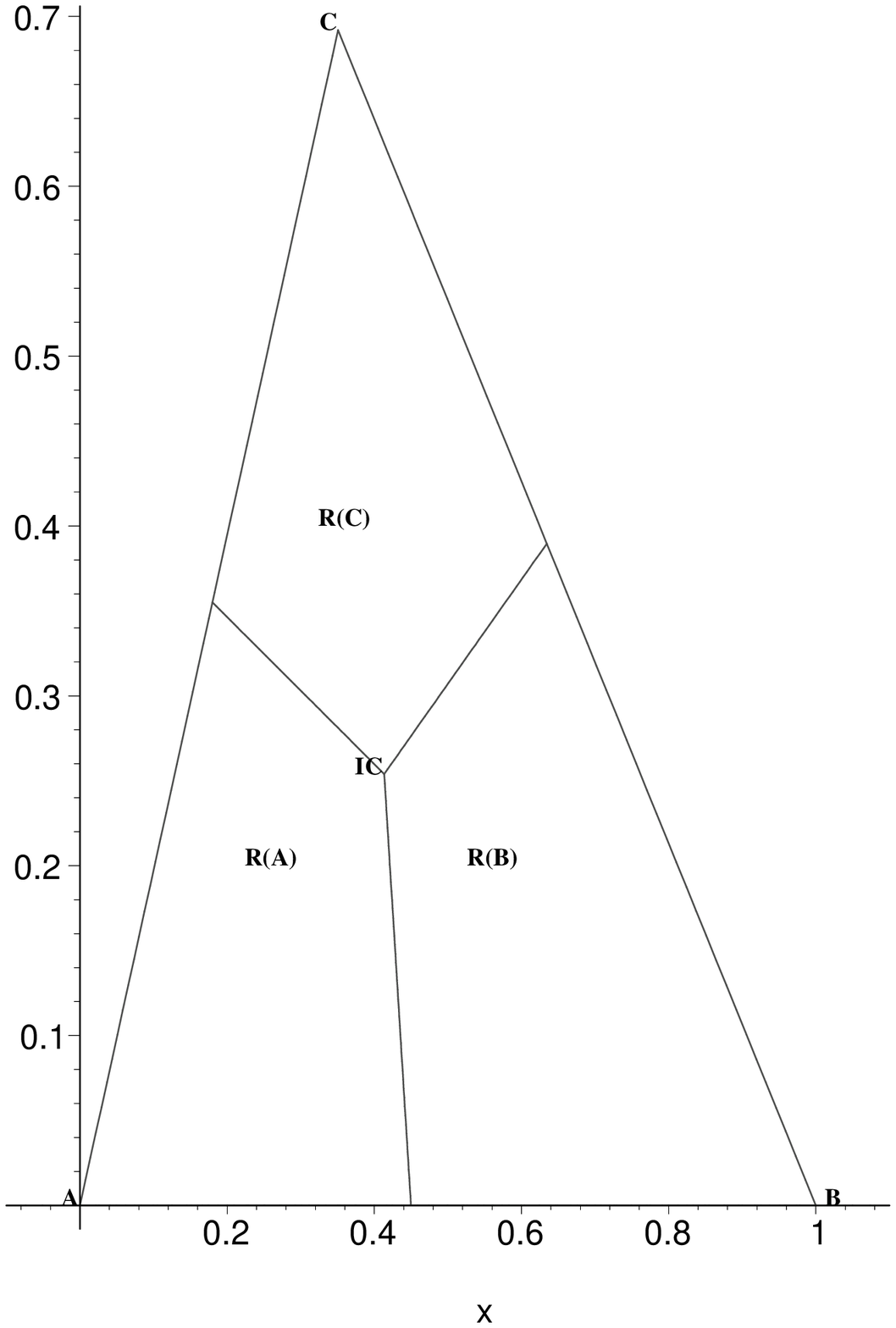, height=140pt , width=200pt}
\end{center}
\caption{
\label{fig:ic-vertex-regions}
The $IC$-vertex regions with orthogonal projections (left) and
with extension of the line segments joining the vertices to $M_{IC}$ (right).}
\end{figure}

One can also define the incenter vertex regions by using the inner angle bisectors.
With orthogonal projections, $IC$-vertex regions are bounded by the edges of $\TY$ and
the inradii crossing the tangential points of the incircle on the edges.
These three regions are denoted as $R^{\perp}_{IC}(\y)$ for $\y \in \Y_3$.
With the inner angle bisectors,
the incenter $M_I$ is used
and the parts of the inner angle bisectors that join $M_I$ to the edges.
These vertex regions are denoted as $R_{IC}(\y)$.
See Figure \ref{fig:ic-vertex-regions} for both versions of the vertex regions.

Note that one might also use the orthocenter, $M_O$, to define the vertex regions.
However, for non-acute triangles $OC$-vertex regions cannot naturally be defined.

\subsection{Edge Regions}
\label{sec:edge-regions}
The spherical proximity region seen earlier is constructed by using the vertex
region based on the closest vertex, $\argmin_{\y \in \Y_3}d(x,\y)$.
One can also use the closest edge, $\argmin_{i \in \{1,2,3\}}d(x,e_i)$,
in defining a proximity region,
which suggests the concept of \emph {edge regions}.

While using the edge $\argmin_{i \in \{1,2,3\}}d(x,e_i)$,
the triangle is again partitioned into three regions
whose intersection is some point $M$ with Euclidean
distance to the edges $d(M,e_1)=d(M,e_2)=d(M,e_3)$,
so $M$ is in fact the incenter of $\TY$ and $d(M,e)=r_{ic}$ is the inradius
(see Section \ref{sec:triangle-centers} for incenter and inradius).

\begin{definition}
The connected regions that partition the triangle, $\TY$, in such a way that
each region has one and only one edge of $\TY$ on its boundary,
are called \emph{edge regions}. $\square$
\end{definition}

This definition implies that there are exactly three edge regions which
intersect at only one point, $M$ in $\TY^o$, the interior of $\TY$.
In fact, one can describe the edge regions starting with $M$.
Join the point $M$ to the vertices by curves such that
the resultant regions satisfy the above definition.
Such regions are called \emph{$M$-edge regions} and
the edge region for edge $e$ is denoted as $R_M(e)$ for $e \in \{e_1,e_2,e_3\}$.
Unless stated otherwise, $M$-edge regions will refer to the regions
constructed by joining $M$ to the vertices by straight lines.
In particular, one can use a \emph{center} of $\TY$
for the starting point $M$.
One can also consider the points in $R_M(e)$ to be
``closer" to $e$ than to the other edges.
Furthermore, it is reasonable to require that the area of
the region $R_M(e)$ get larger as $d(M,e)$ increases.
Moreover, in higher dimensions,
the corresponding regions are called ``face regions".

The functional forms of $R_{M}(e_i)$ for $i \in \{1,2,3\},$ for $M=(m_1,m_2) \in \TY^o$
and $m_1 > c_1$ in the basic triangle are given by
\begin{align*}
R_{M}(e_1)&=\left \{(x,y) \in \TY :
y \ge \frac{m_2\,(x-1)}{m_1-1};\; y \ge \frac{c_2\,(x-m_1)-m_2\,(x-c_1)}{c_1-m_1} \right\},\\
R_{M}(e_2)&=\left \{(x,y) \in \TY :
y \ge \frac{m_2\,x}{m_1};\; y \le \frac{c_2\,(x-m_1)-m_2\,(x-c_1)}{c_1-m_1} \right\}, \text{ and } \\
R_{M}(e_3)&=\left \{(x,y) \in \TY :
\frac{m_2\,x}{m_1};\; y \le \frac{m_2\,(x-1)}{m_1-1} \right\}.
\end{align*}

If $x$ falls on the boundary of two $M$-edge regions, then
$x$ is arbitrarily assigned to one of the $M$-edge regions.

The center of mass edge regions ($CM$-edge regions) are described in detail,
as we will use them in defining a new class of proximity maps.

\subsubsection{$CM$-Edge Regions}
\begin{figure}[ht]
\begin{center}
\psfrag{A}{\scriptsize{$\y_1$}}
\psfrag{B}{\scriptsize{$\y_2$}}
\psfrag{C}{\scriptsize{$\y_3$}}
\psfrag{R(AB)}{\scriptsize{$R_{CM}(e_3)$}}
\psfrag{R(BC)}{\scriptsize{$R_{CM}(e_1)$}}
\psfrag{R(AC)}{\scriptsize{$R_{CM}(e_2)$}}
\psfrag{CM}{\scriptsize{$M_{C}$}}
\psfrag{x}{}
\epsfig{figure=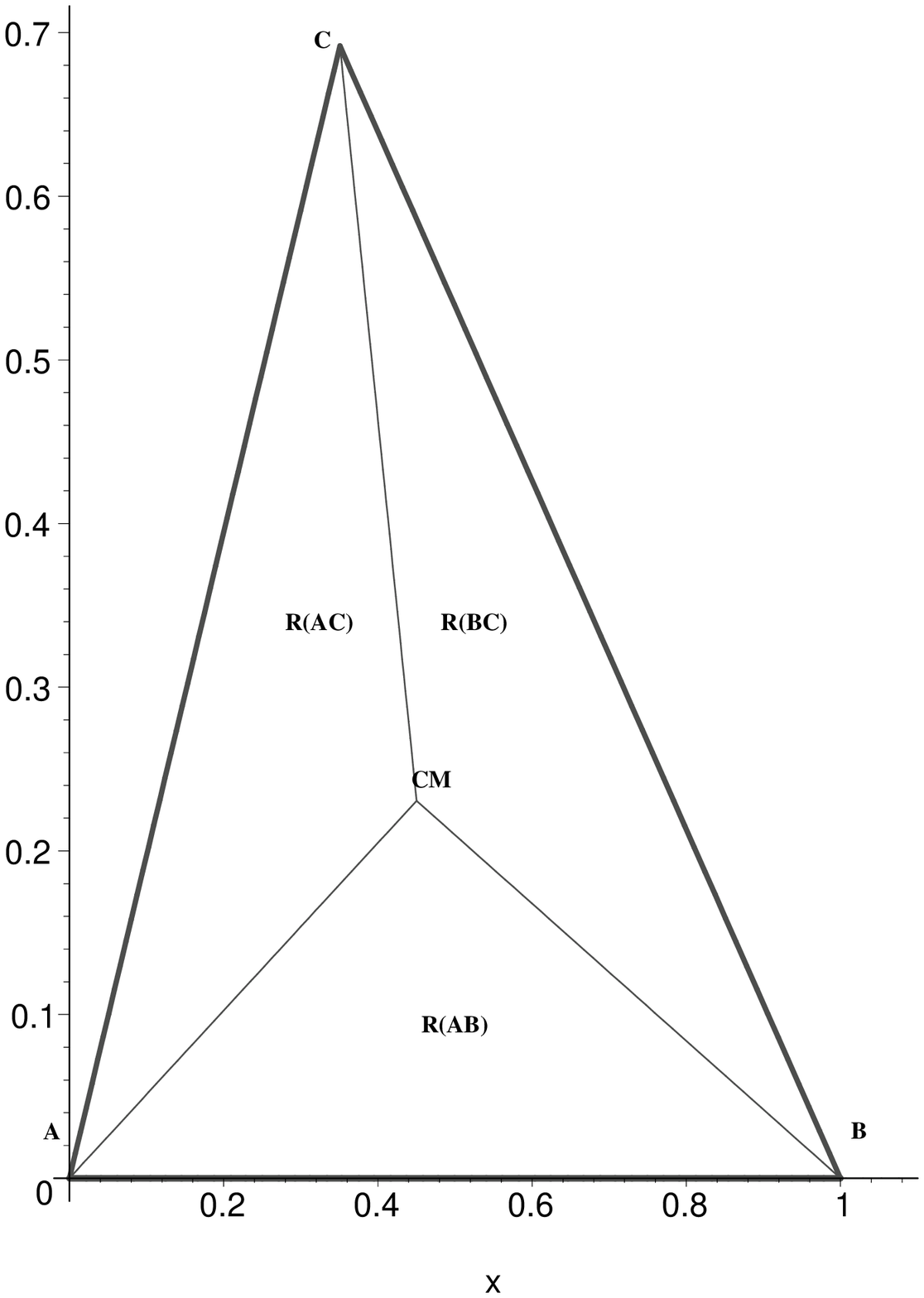, height=140pt , width=200pt}
\psfrag{R(AB)}{\scriptsize{$R_{IC}(e_3)$}}
\psfrag{R(BC)}{\scriptsize{$R_{IC}(e_1)$}}
\psfrag{R(AC)}{\scriptsize{$R_{IC}(e_2)$}}
\psfrag{IC}{\scriptsize{$M_I$}}
\psfrag{x}{}
\epsfig{figure=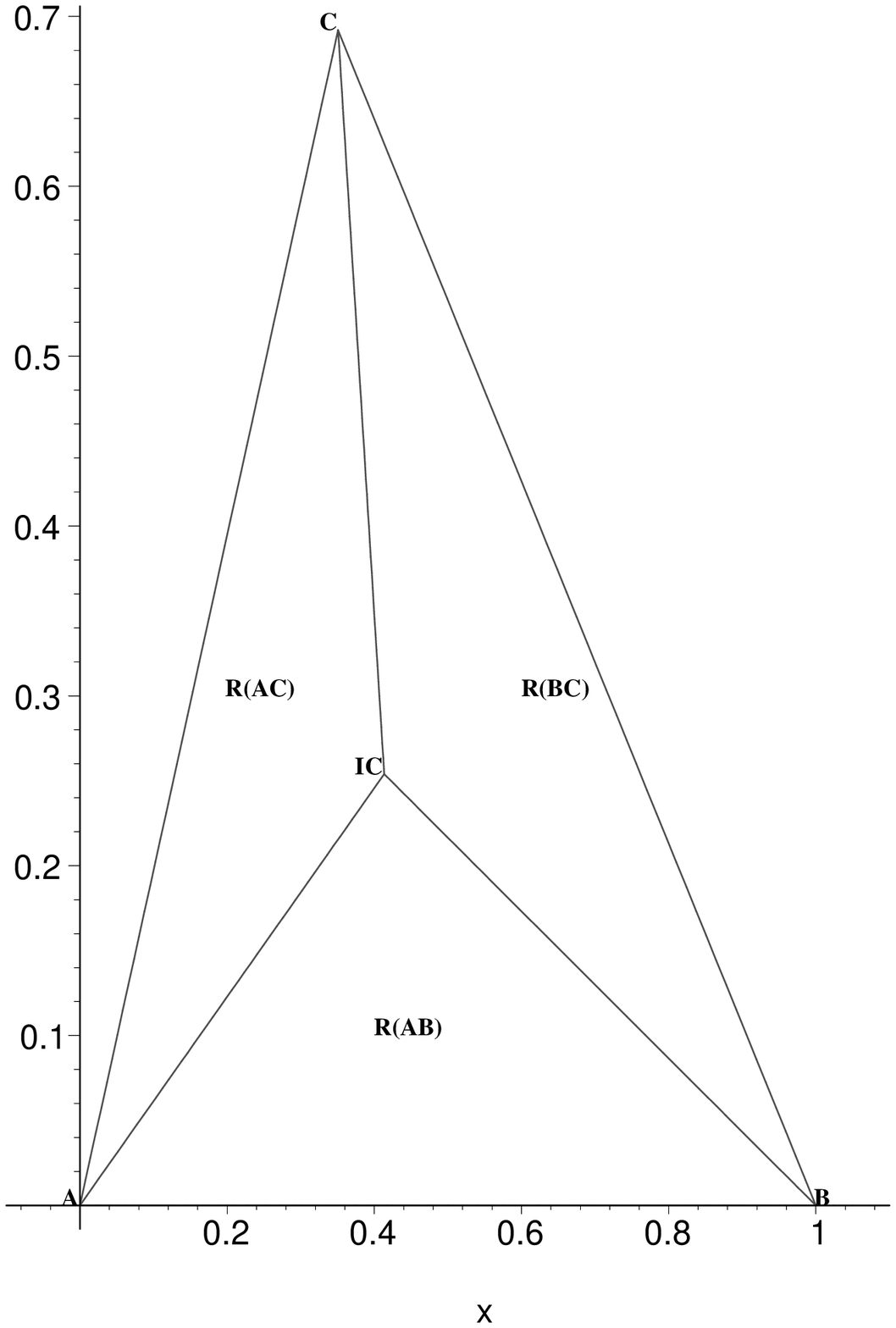, height=140pt , width=200pt}
\end{center}
\caption{$CM$-edge regions $R_{CM}(e_i)$ (right) and $R_{IC}(e_i)$ (left) for $i \in \{1,2,3\}$.}
\label{fig:edge-regions}
\end{figure}

One can divide $\TY$ into three regions by using the median
lines which intersect at the centroid,
or equivalently, joining the centroid $M_C$ to the vertices by
straight lines will yield the $CM$-edge regions.
Let $R_{CM}(e)$ be the region for edge $e \in \left\{ e_1,e_2,e_3 \right \}$.
See Figure \ref{fig:edge-regions} (left).

The functional forms of $R_{CM}(e_i)$ for $i \in \{1,2,3\},$ in the basic triangle,
$T_b$, are given by
\begin{align*}
&R_{CM}(e_1)=\left \{(x,y) \in T\left( \Y_3 \right):
y \le \frac{c_2\,(1-2\,x)}{1-2\,c_1};\; y \ge \frac{c_2\,(1-x)}{2-c_1} \right\},\\
&R_{CM}(e_2)=\left \{(x,y) \in T\left( \Y_3 \right):
y \ge \frac{c_2\,(1-2\,x)}{1-2\,c_1};\; y \ge \frac{c_2\,x}{1+c_1} \right \},\\
&R_{CM}(e_3)=\left \{(x,y) \in T\left( \Y_3 \right):
y \le \frac{c_2\,x}{1+c_1} ;\; y \le \frac{c_2\,(1-x)}{2-c_1} \right\}.
\end{align*}

%If $x$ falls on the boundary of two $CM$-edge regions, then
%$x$ is arbitrarily assigned to one the $CM$-edge regions.

\begin{remark}
One can also divide $\TY$ into three regions by using the inner angle bisectors
which intersect at the incenter, yielding the incenter edge regions ($IC$-edge regions).
Let $R_{IC}(e)$ be the region for edge $e \in \{e_1,e_2,e_3\}$.
Notice that the closest edge to any point in $R_{IC}(e)$ is edge $e$, i.e.,
$x \in R_{IC}(e)$ iff $\argmin_{u \in \{e_1,e_2,e_3\}}d(x,u)=e$.
If two edges are equidistant from $x$,
$x$ is arbitrarily assigned to an edge region.
See Figure \ref{fig:edge-regions} (right).
$\square$
\end{remark}

\begin{remark} In $\R$, one can view the end points of $[0,1]$,
$\{0,1\}$, as vertices or edges. So $[0,1/2]$ and $[1/2,1]$
can be viewed as either vertex regions or edge regions.
$\square$
\end{remark}

\section{Proximity Regions in Delaunay Tessellations}
\label{sec:prox-regions-in-Deltri}
Let $\Y_m=\left \{\y_1,\ldots,\y_m \right\}$
be $m$ points in general position in $\R^d$ and $\T_j$ be
the $j^{th}$ Delaunay cell for $j=1,\ldots,J$.
Let also that $\X_n$ be a
random sample from $F$ with support $\mS(F) \subseteq \C_H(\Y_m)$.
That is, $\Omega=\C_H(\Y_m)$ and the $\Omega_j=\T_j$ with $\mu$ being the
Lebesgue measure. Then the appealing properties for proximity
regions in Section \ref{sec:prox-maps} become:
\begin{itemize}
\item[\textbf{P1}]
$N(x)$ is well defined for all $x \in \C_H(\Y_m)$.
\item[\textbf{P2}]
$x \in N(x)$ for all $x \in \C_H(\Y_m)$.
\item[\textbf{P3}]
$x$ is at the \emph{center} of $N(x)$ for all $x \in \C_H(\Y_m)$.
\item[\textbf{P4}]
For $x \in \T_j \subseteq \C_H(\Y_m)$, $N(x)$ and
$\T_j$ are of the \emph{same type}; they are both $(d+1)$-simplicies.
\item[\textbf{P5}]
For $x \in \T_j \subseteq \C_H(\Y_m)$, $N(x)$ mimics the shape of $\T_j$;
i.e., it is \emph{similar} to $\T_j$.
\item[\textbf{P6}]
Conditional on $X \in \T_j$, $N(X)$ is a proper subset of $\T_j$ a.s.
\item[\textbf{P7}]
For $x \in \T_j$ and $y \in \T_k$ with $j\not=k$,
$N(x)$ and $N(y)$ are disjoint a.s.
\item[\textbf{P8}]
The size of $N(x)$ is continuous in $x$;
that is, for each $\varepsilon >0$ there exists a $\delta(\varepsilon)>0$
such that $|\mu(N(y))-\mu(N(x))|<\delta(\varepsilon)$ whenever $||y-x|| <\varepsilon$.
%\item[\textbf{P9}]
%The size (i.e., measure) of $\NS(x)$ increases as $d(x,\Y_m)$ increases.
\item[\textbf{P9}]
The arc probability $\mu(N_S)$ does not depend on the support region for uniform data in $\R^d$.
\end{itemize}

In particular, for illustrative purposes, we focus on $\R^2$,
where a Delaunay tessellation is a triangulation,
provided that no more than three points of $\Y_m$ are cocircular.
Furthermore, for simplicity,
let $\Y_3=\{\y_1,\y_2,\y_3\}$ be three non-collinear points
in $\R^2$ and $\TY=T(\y_1,\y_2,\y_3)$ be the triangle
with vertices $\Y_3$.
Let $\X_n$ be a random sample from $F$ with
support $\mS(F) \subseteq \TY$.
The spherical proximity map is the first
proximity map defined in literature
(see \cite{devinney:2002a}, \cite{marchette:2003}, \cite{priebe:2003b}, \cite{priebe:2003a},
and \cite{devinney:2006})
where $CC$-vertex regions were implicitly used for points in $\C_H\left( \Y_m \right)$.
In the following sections, we will describe arc-slice proximity maps $\NAS(\cdot)$
and define two families of triangular proximity regions for which
 \textbf{P4} and \textbf{P5} will automatically hold.

%\subsection{Spherical Proximity Maps}

\subsection{Arc-Slice Proximity Maps}
Recall that for $N_S(\cdot)$ \textbf{P7} is violated,
since for any $x \in \T_j \subset \R^d$,
$B(x,r(x)) \not \subset \T_j$, which implies that two
proximity regions $N_S(x)$ and $N_S(y)$ might overlap
for $x,y$ in two distinct cells.
Such an overlap of the regions make the distribution of the
domination number of the PCD associated with $N_S(\cdot)$,
if not impossible, hard to calculate.
In order to avoid the overlap of
regions $B(x,r(x))$ and $B(y,r(y))$ for $x,\,y$ in different Delaunay cells,
the balls are restricted to the corresponding cells,
which leads to \emph{arc-slice proximity regions},
$\NAS(x):=\overline B(x,r(x)) \cap \T_j$, where $\overline B(x,r(x))$
is the closure of the ball $B(x,r(x))$.
The closed ball is used in the definition of the arc-slice
proximity map for consistency with the other proximity maps
that will be defined on Delaunay cells.
The arc-slice proximity map $\NAS(x)$ is well-defined only in $\C_H\left( \Y_m \right)$,
provided that $\Y_m$ is in general position and $m \ge d+1$ in $\R^d$.
%So \textbf{P1} holds. $\NAS(\cdot)$ also satisfies \textbf{P2}, \textbf{P6}, \textbf{P7}, and \textbf{P8}.
%We define new types of proximity maps that satisfy more
%of the properties \textbf{P1}-\textbf{P8} and
%introduce new explanatory concepts for investigating the
%properties of these proximity maps.
%We focus on possible
%extensions to multiple dimensions based on the Delaunay
%tessellation of $\Y_m$, hence the extensions will
%also be well-defined only in $\C_H\left( \Y_m \right)$.

\begin{figure}[ht]
\begin{center}
\psfrag{A}{\scriptsize{$\y_1$}}
\psfrag{B}{\scriptsize{$\y_2$}}
\psfrag{C}{\scriptsize{$\y_3$}}
\psfrag{CC}{\scriptsize{$M_{CC}$}}
\psfrag{M1}{\scriptsize{$M_3$}}
\psfrag{M2}{\scriptsize{$M_1$}}
\psfrag{M3}{\scriptsize{$M_2$}}
\psfrag{y}{}
\psfrag{x}{}
\psfrag{x1}{\scriptsize{$x$}}
\epsfig{figure=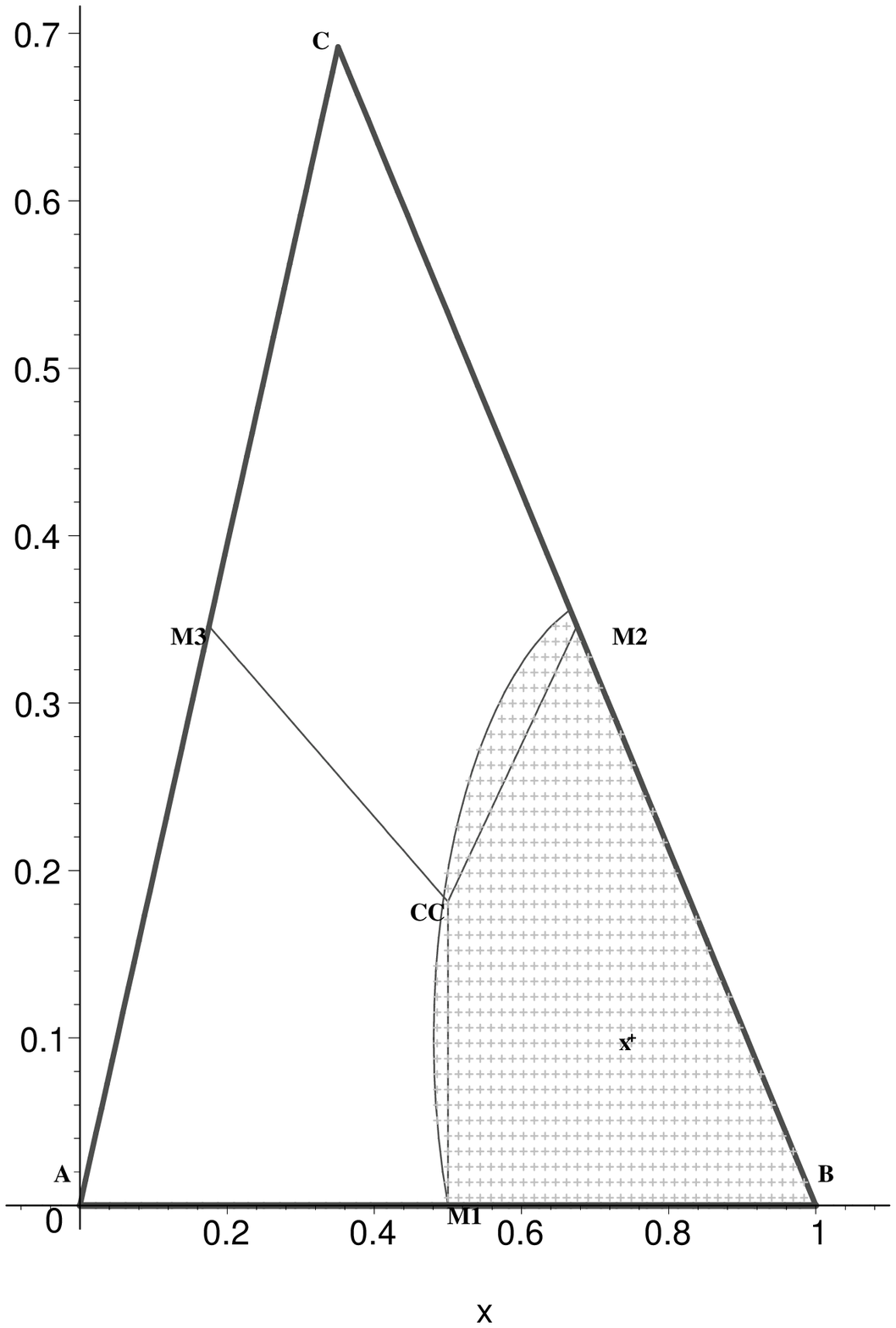,  height=140pt , width=190pt}\\
\psfrag{CM}{\scriptsize{$M_{C}$}}
\epsfig{figure=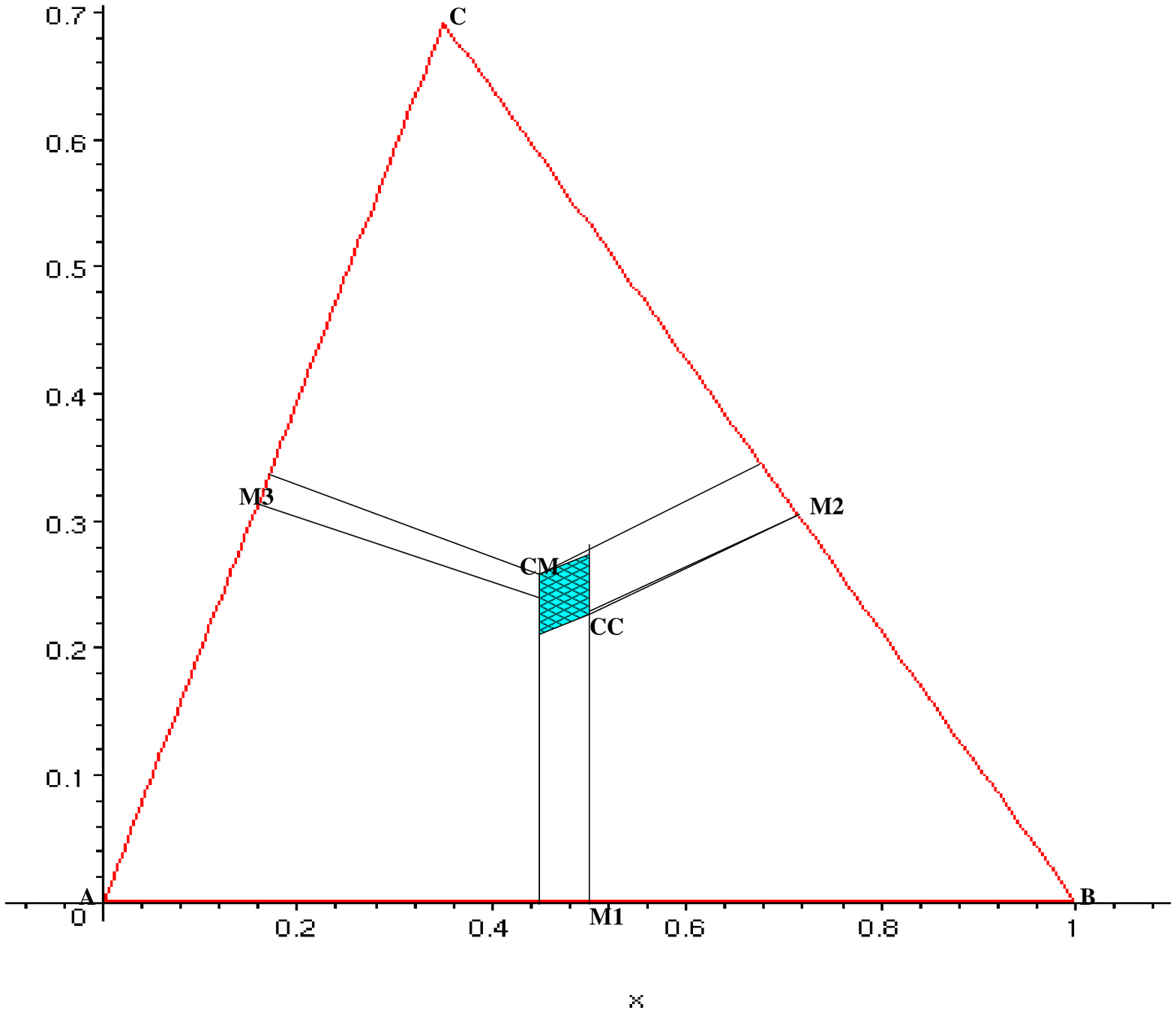, height=145pt , width=210pt}
\psfrag{IC}{\scriptsize{$M_{IC}$}}
\epsfig{figure=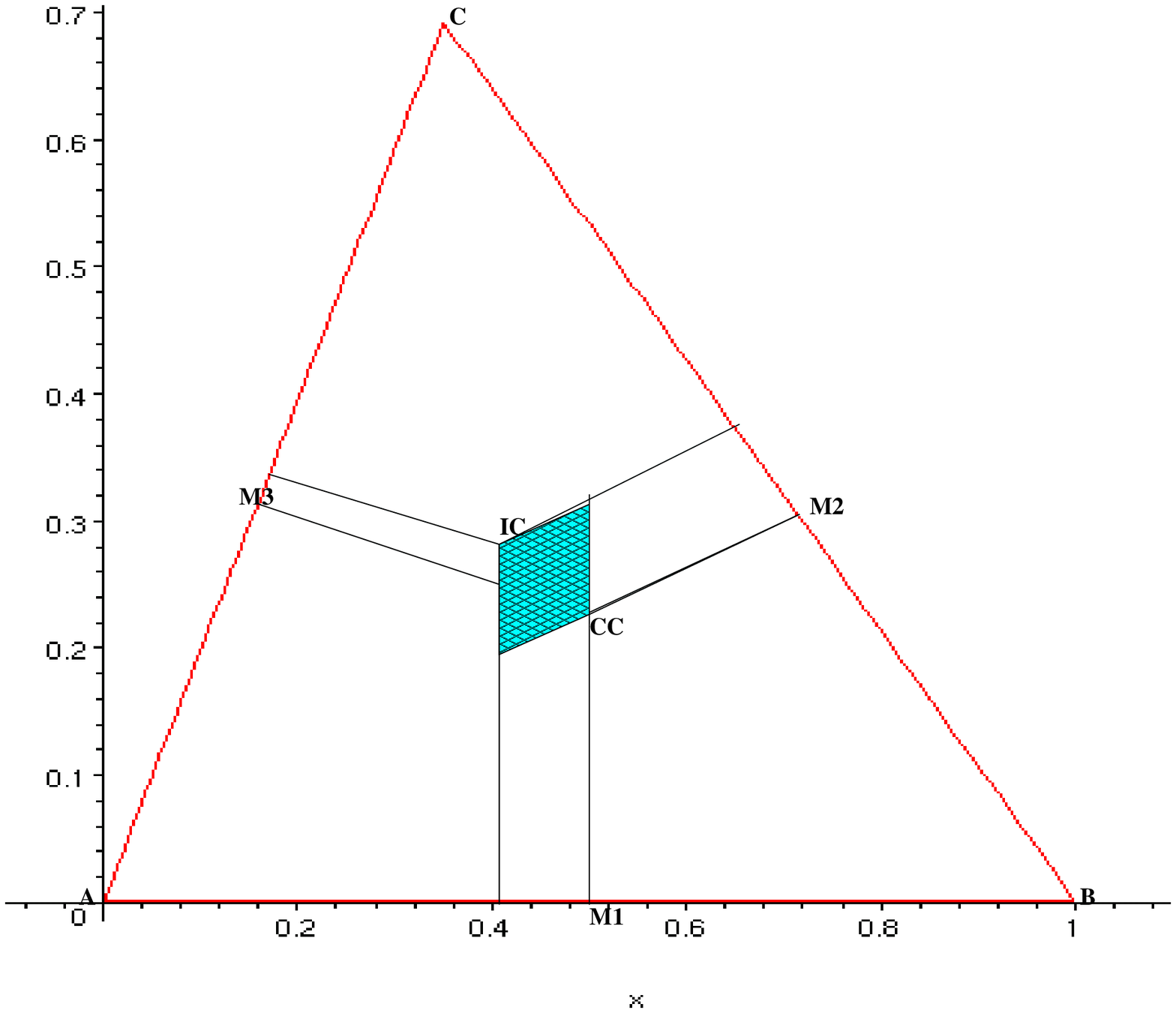, height=145pt , width=210pt}
\end{center}
\caption{$\NAS(x,M_{CC})$ with an $x \in R_{CC}(\y_2)$ (top)
and the superset regions $\RS(\NAS,M_C)$ (bottom left)
$\RS(\NAS,M_{IC})$ in $\TY$ (bottom right).}
\label{fig:AS-region}
\end{figure}

\begin{figure}[ht]
\centering
\rotatebox{-90}{ \resizebox{2.0 in}{!}{\includegraphics{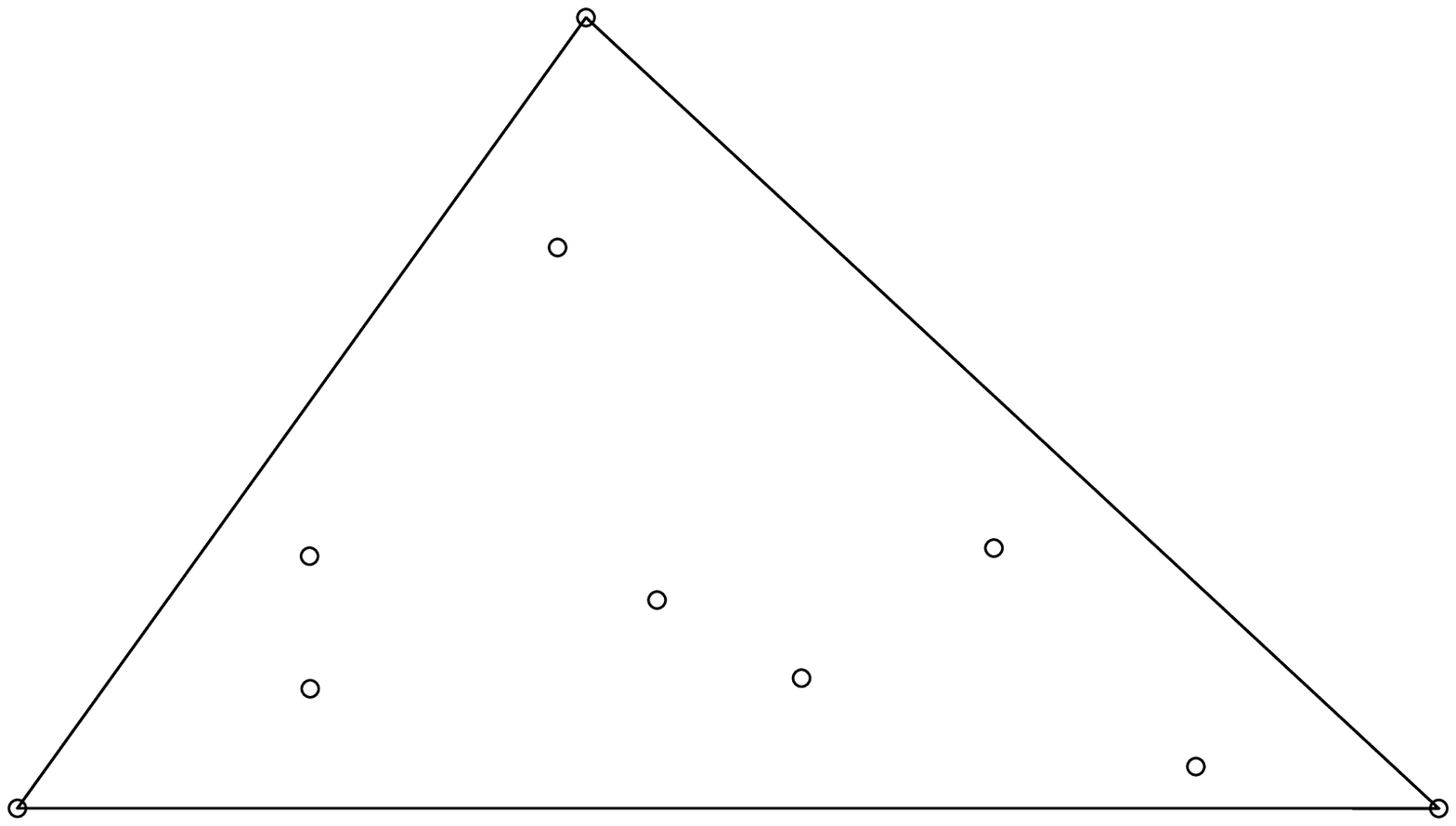} }}
\rotatebox{-90}{ \resizebox{2.0 in}{!}{\includegraphics{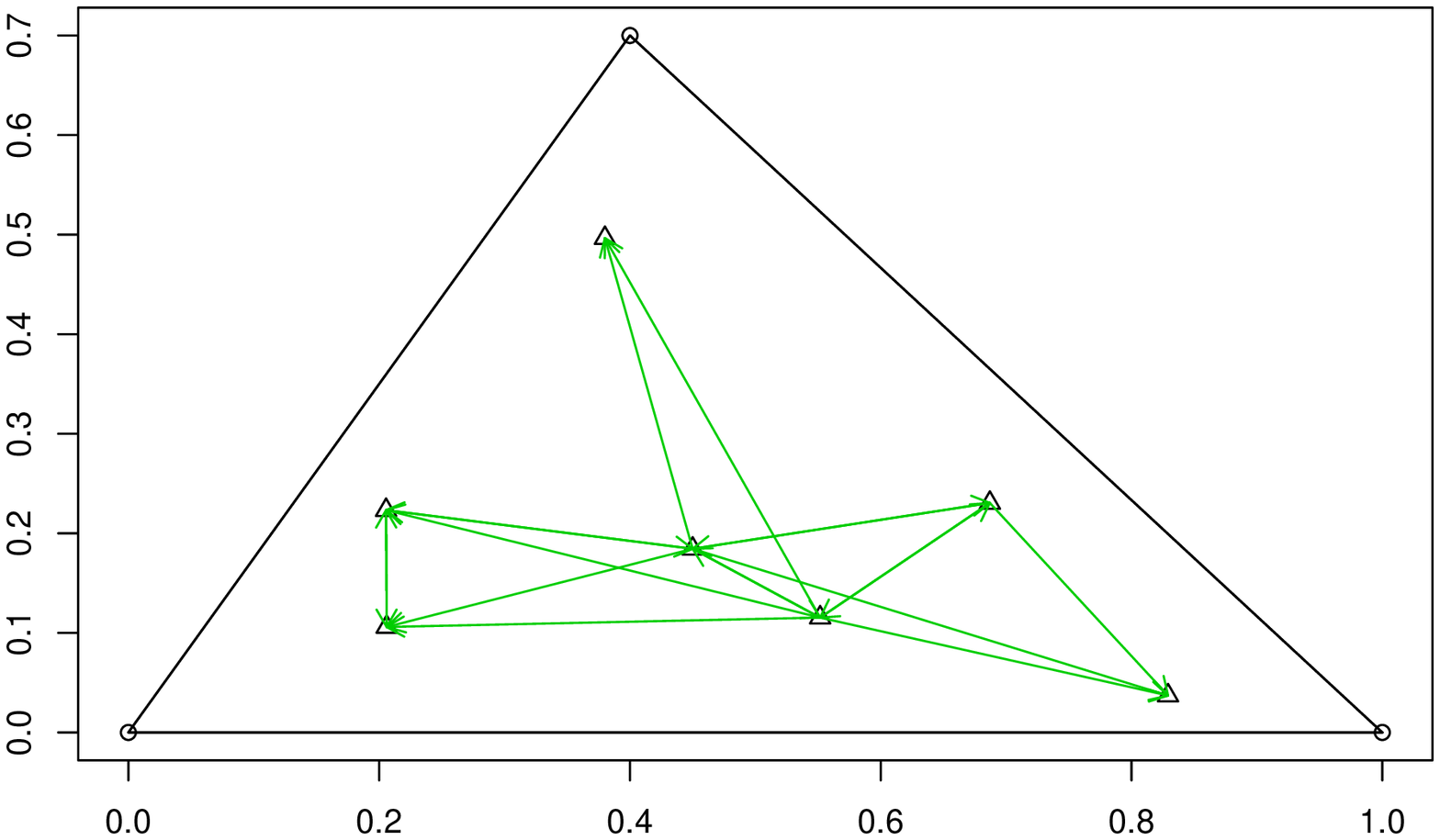} }}
\caption{
\label{fig:AS-arcs-1T}
A realization of 7 $\X$ points generated iid $\U(\TY)$ (left)
and the corresponding arcs for $\NAS^{r=2}(x,M_C)$ (right).
}
\end{figure}

\begin{figure}[ht]
\centering
%\rotatebox{-90}{ \resizebox{2.0 in}{!}{\includegraphics{DelaunayXY.ps} }}
%\rotatebox{-90}{ \resizebox{2.0 in}{!}{\includegraphics{DelaunayCH.ps} }}
\rotatebox{-90}{ \resizebox{2.0 in}{!}{\includegraphics{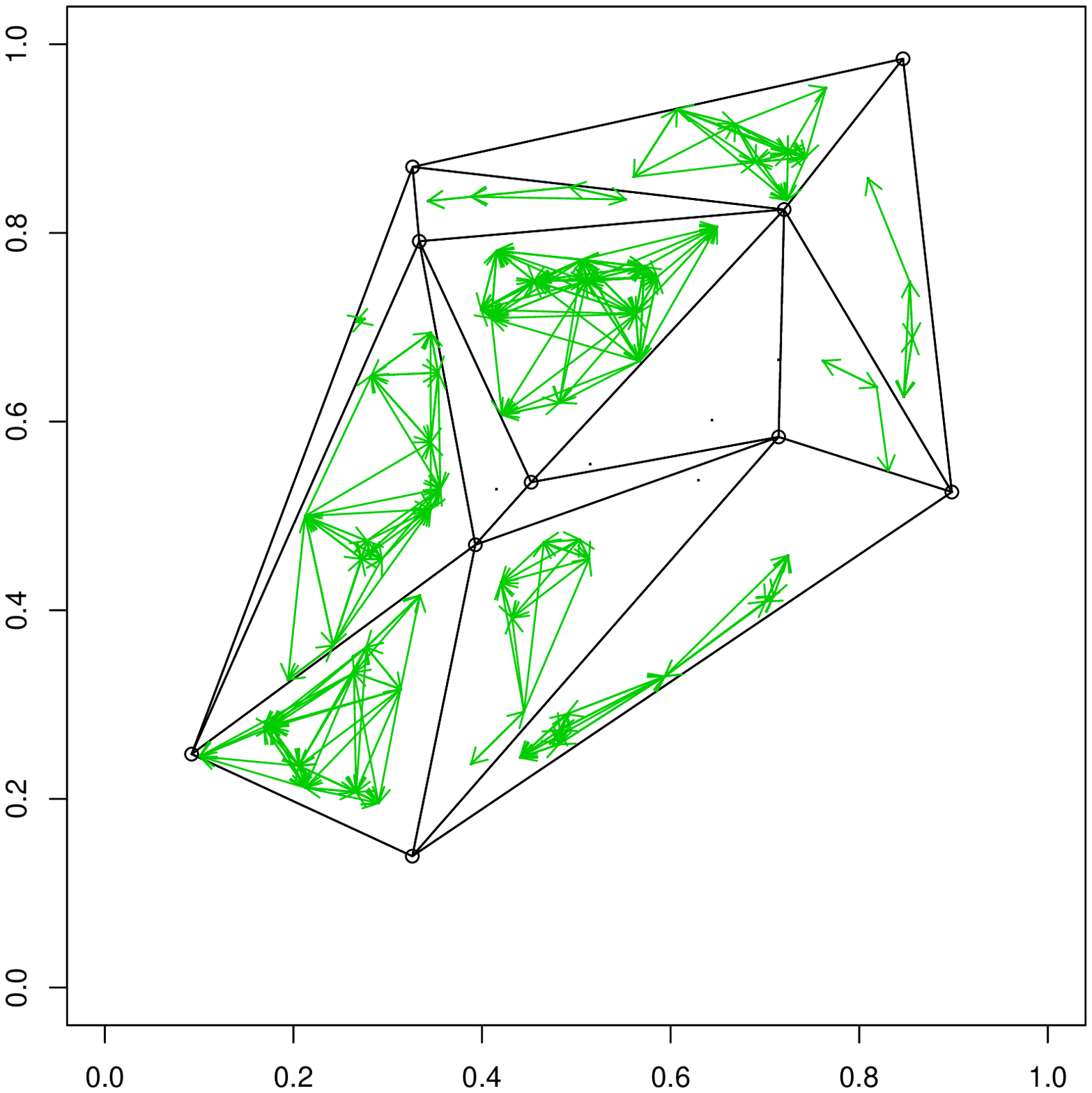} }}
\caption{
\label{fig:AS-arcs-multiT}
The arcs for arc-slice PCDs with $\NAS^{r=2}(x,M_C)$
for the 77 $\X$ points that lie in the $\C_H\left(\Y_{10}\right)$ (see Figure \ref{fig:deltri}).
}
\end{figure}

%Recall that the arc-slice proximity region of $x \in \TY$ is given by
% $\NAS(x):=\overline B(x,r(x)) \cap \TY$ where $\overline B(x,r(x))$ is the closed
%ball centered at $x$ with radius $r(x):=\min_{\y \in \Y_m}d(x,\y)$.
By construction, the $CC$-vertex regions are implicitly used,
since $x \in R_{CC}(\y)$ iff $\y=\argmin_{u \in \Y_m}d(x,u)$.
To make this dependence explicit, the notation $\NAS(\cdot,M_{CC})$ is used.
See Figure \ref{fig:AS-region} (top) for
$\NAS(x,M_{CC})$ for an $x \in R_{CC}(\y_2)$.
The functional form of $\NAS(x,M_{CC})$ for an $x\in R_{CC}(\y)$ is given by
$$\NAS(x,M_{CC}):=\bigl\{z \in \TY: d(z,x) \le r(x)=d(x,\y) \bigr\}.$$
Notice that, the region $\NAS(x,M_{CC})$ is a closed region,
unlike $N_S(x)=B(x,r(x))$.
The properties \textbf{P1}, \textbf{P2}, \textbf{P7} hold by definition.
Notice that $\NAS(x,M_{CC}) \subseteq \TY$ for all $x \in \TY$ and
$\NAS(x,M_{CC})=\TY$ iff $x=M_{CC}$, since $\overline B(x,r(x)) \supset \TY$
only when $x=M_{CC}$.
Hence the superset region for arc-slice proximity maps with
$CC$-vertex regions is $\RS(\NAS,M_{CC})=\{M_{CC}\}$.
So \textbf{P6} follows.
Furthermore, \textbf{P8} holds since the area $A\left(\NAS(x,M_{CC}) \right)$
is a continuous function of $r(x)=\min_{\y \in \Y_3}d(x,\y)$
which is a continuous function of $x$.
\textbf{P3}, \textbf{P4}, \textbf{P5}, and \textbf{P9} fail for $\NAS(x,M_{CC})$.
See Figure \ref{fig:AS-arcs-1T} for the arcs
based on $\NAS(x,M_C)$ for a realization
of 7 $\X$ points in the one triangle case,
and Figure \ref{fig:AS-arcs-multiT} for the arcs
for the realization of 77 $\X$ points in the
multi-triangle case in Figure \ref{fig:deltri} (top right).

One can define arc-slice proximity regions with
any type of $M$-vertex regions as
$$\NAS(x,M):=\overline B(x,r(x)) \cap \TY\text{ where }r(x):=d(x,\y) \text{ for }x \in R_M(\y).$$
But for $M \not= M_{CC}$, $\NAS(\cdot,M)$ satisfies only
\textbf{P1}, \textbf{P2}, and \textbf{P7}.
\textbf{P6} fails to hold since $\RS(\NAS,M)$ has
 positive area and \textbf{P8} fails since the size of $\NAS(x,M)$
is not continuous in $x$.
% for $M \not= M_{CC}$.
See, for example, Figure \ref{fig:AS-region} (right) for $\RS(\NAS,M_C)$.
In terms of the properties in Section \ref{sec:prox-maps},
$\NAS(\cdot,M_{CC})$ is the most appealing proximity map among the family
$\mathscr N_{AS}:=\{\NAS(\cdot,M): M \in \R^2 \setminus \Y_3 \}$.

Moreover, $\Lambda_0(\NAS,M)=\Y_3$ for all $M \in \R^2 \setminus \Y_3$
since $\lambda(\NAS(x,M))=0$ iff $x \in \Y_3$.

Next, we define \emph{triangular proximity regions},
which, by definition, will satisfy properties \textbf{P4} and \textbf{P5}.
These proximity regions are the building blocks of the PCDs for
which more rigorous mathematical analysis ---
compared to the PCDs based on spherical
and arc-slice proximity maps --- will be possible.
%See, for example, \cite{ceyhan:dom-num-NPE-SPL},
%\cite{ceyhan:arc-density-PE}, and \cite{ceyhan:arc-density-CS}.

\subsection{Proportional-Edge Proximity Maps}
\label{sec:r-factor}
The first type of triangular proximity map introduced
is the proportional-edge proximity map.
For this proximity map,
the asymptotic distribution of domination number and
the relative density of the corresponding
PCD will have mathematical tractability
(\cite{ceyhan:dom-num-NPE-SPL}, \cite{ceyhan:arc-density-PE} and \cite{ceyhan:dom-num-NPE-MASA}).

\begin{figure}[ht]
\centering
\scalebox{.45}{\input{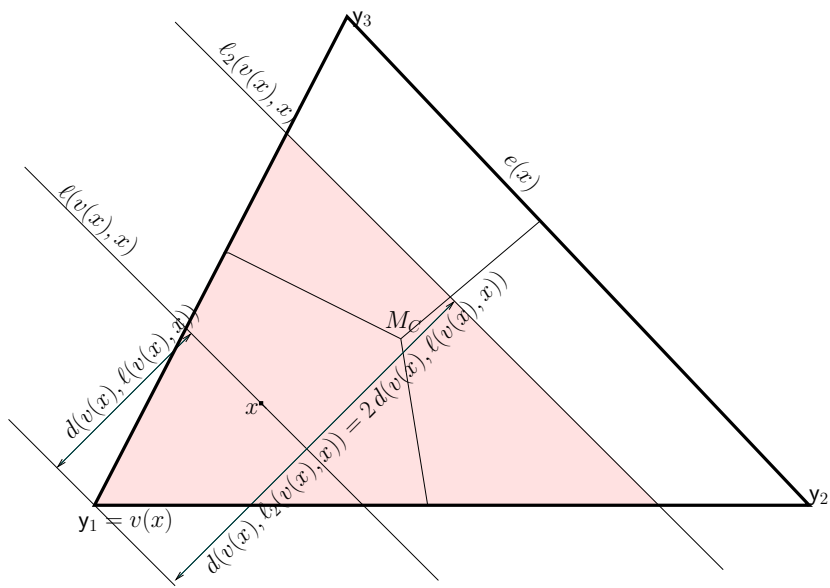}}
\caption{
\label{fig:ProxMapDef1}
Construction of proximity region, $\NPE^2(x)$ (shaded region)
for an $x \in R_{CM}(\y_1)$.}
\end{figure}

For $r \in [1,\infty]$,
define $\NPE^r(\cdot,M):=N\left(\cdot,M;r,\Y_3\right)$ to be the
\emph{proportional-edge proximity map}
with $M$-vertex regions as follows
(see also Figure  \ref{fig:ProxMapDef1} with $M=M_C$ and $r=2$).
For $x \in \TY \setminus \Y_3$, let $v(x) \in \Y_3$ be the
vertex whose region contains $x$; i.e., $x \in R_M(v(x))$.
If $x$ falls on the boundary of two $M$-vertex regions,
$v(x)$ arbitrarily assigned.
Let $e(x)$ be the edge of $\TY$ opposite $v(x)$.
Let $\ell(v(x),x)$ be the line parallel to $e(x)$ through $x$.
Let $d(v(x),\ell(v(x),x))$ be the Euclidean (perpendicular) distance
from $v(x)$ to $\ell(v(x),x)$.
For $r \in [1,\infty)$, let $\ell_r(v(x),x)$ be the line parallel to $e(x)$
such that
\begin{gather*}
d(v(x),\ell_r(v(x),x)) = r\,d(v(x),\ell(v(x),x))\\
\text{ and }\\
d(\ell(v(x),x),\ell_r(v(x),x)) < d(v(x),\ell_r(v(x),x)).
\end{gather*}
Let $T_r(x)$ be
the triangle similar to
and with the same orientation as $\TY$
having $v(x)$ as a vertex
and $\ell_r(v(x),x)$ as the opposite edge.
Then the {\emph r-factor proportional-edge proximity region}
$\NPE^r(x,M)$ is defined to be $T_r(x) \cap \TY$.
Notice that $\ell(v(x),x)$ divides the edges of $T_r(x)$
(other than $\ell_r(v(x),x)$) proportionally with the factor $r$.
Hence the name \emph{proportional edge proximity region} and the notation $\NPE^r(\cdot,M)$.

\begin{figure}[ht]
\centering
%\rotatebox{-90}{ \resizebox{2.0 in}{!}{\includegraphics{Pointsin1T.ps} }}
\rotatebox{-90}{ \resizebox{2.0 in}{!}{\includegraphics{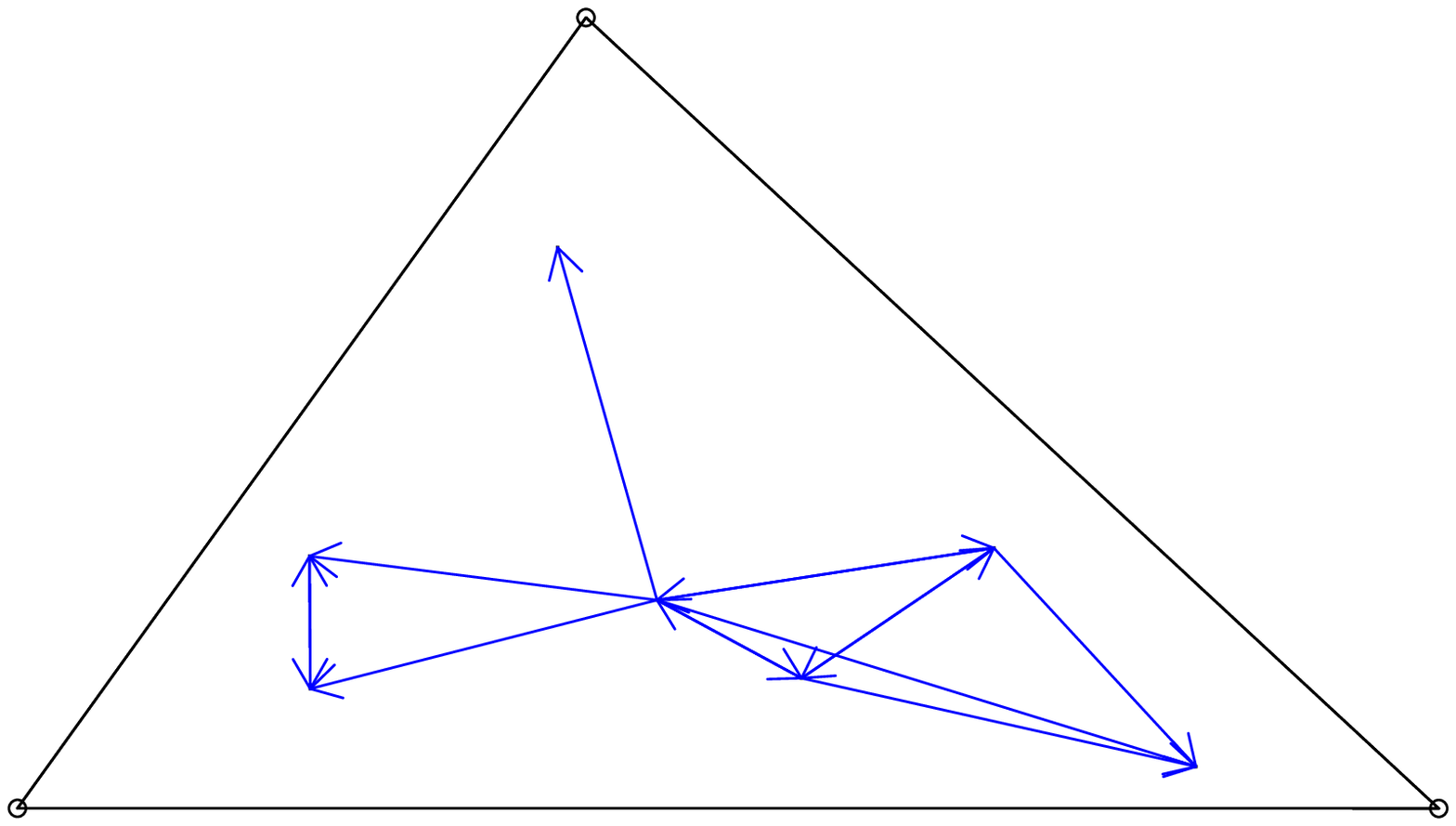} }}
\rotatebox{-90}{ \resizebox{2.0 in}{!}{\includegraphics{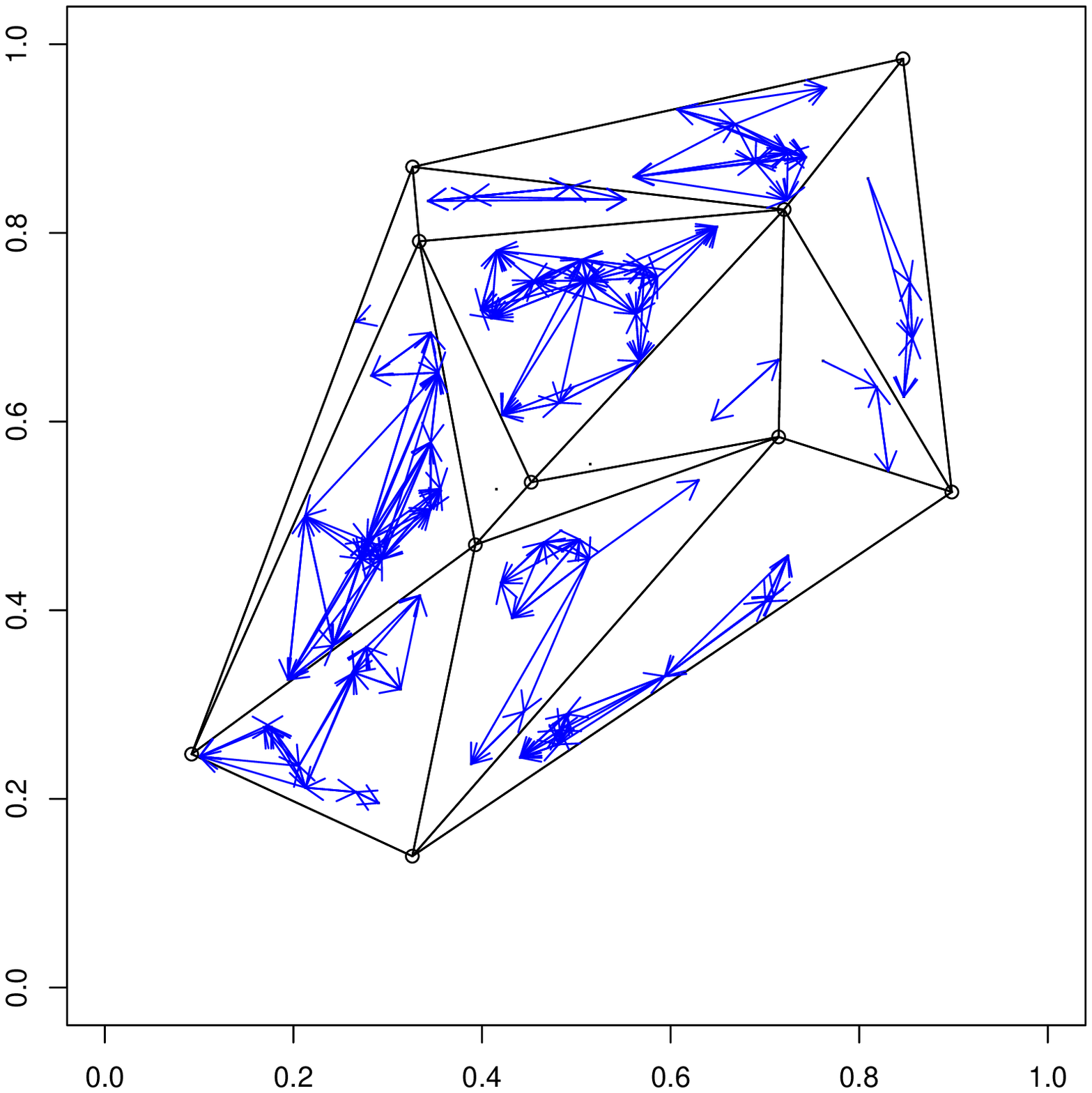} }}
\caption{
\label{fig:PE-arcs}
The arcs for $\NPE^{r=2}(x,M_C)$ the 7 $\X$ points in Figure \ref{fig:AS-arcs-1T},
and the arcs for $\NPE^{r=2}(x,M_C)$ for the 77 $\X$ points that lie in the $\C_H\left(\Y_{10}\right)$
in Figure \ref{fig:deltri}.
}
\end{figure}

Notice that $r \ge 1$ implies $x \in \NPE^r(x,M)$.
Furthermore,
$\lim_{r \rightarrow \infty} \NPE^r(x,M) = \TY$ for all $x \in \TY \setminus \Y_3$,
so $\NPE^{\infty}(x,M) := \TY$ for all such $x$.
For $x \in \Y_3$, $\NPE^r(x,M) := \{x\}$ for all $r \in [1,\infty]$.
See Figure \ref{fig:PE-arcs} for the arcs
based on $\NPE^{r=2}(x,M_C)$ in the one triangle
and the multi-triangle cases.

Notice that $X_i \stackrel{iid}{\sim} F$,
with the additional assumption
that the non-degenerate two-dimensional
pdf $f$ exists with support $\mS(F) \subseteq \TY$,
implies that the special case in the construction
of $\NPE^r$ ---
$X$ falls on the boundary of two vertex regions ---
occurs with probability zero.
Note that for such an $F$, $\NPE^r(X,M)$ is a triangle a.s.

The functional form of $\NPE^r(x,M)$ for $x=(x_0,y_0) \in T_b$ is given by
\begin{itemize}
\item[] $\NPE^r(x,M)=\left \{(x,y)\in T_b:
y \le \frac{r\,(y_0\,(1-c_1)+c_2\,x_0)-c_2\,x}{1-c_1} \right \}$ for $x \in R_M(\y_1)$,
\item[] $\NPE^r(x,M)=\left \{(x,y)\in T_b:
y \le \frac{r\,(y_0\,c_1-c_2\,(x_0+1))-c_2+c_2\,x}{c_1} \right\}$ for $x \in R_M(\y_2)$,
\item[] $\NPE^r(x,M)=\left \{(x,y)\in T_b:
y \ge 2\,y_0-c_2(r-1) \right \}$ for $x \in R_M(\y_3)$.
\end{itemize}
Of particular interest is $\NPE^r(x,M)$ with any $M$ and $r \in \left\{\sqrt{2},\,3/2,2\right\}$.
For $r=\sqrt{2}$, $\ell(v(x),x)$ divides $T_{\sqrt{2}}(x)$ into
two regions of equal area, hence $\NPE^{\sqrt{2}}(x,M)$ is also
referred to as \emph{double-area proximity region}. For $r=2$,
$\ell(v(x),x)$ divides the edges of $T_2(x)$ ---other than
$\ell_r(v(x),x)$ --- into two segments of equal length, hence
$\NPE^2(x,M)$ is also referred to as \emph{double-edge proximity region}.
For $r < 3/2$, $\RS(\NPE^r,M_C)=\emptyset$,
and for $r > 3/2$, $\RS(\NPE^r,M_C)$
has positive area; for $r=3/2$, $\RS(\NPE^r,M_C)=\{M_C\}$.
Therefore, $r=3/2$ is the threshold for $\NPE^r(\cdot,M_C)$
to satisfy \textbf{P6}.
Furthermore, $r=3/2$ is the value at which the asymptotic distribution of the
domination number of the PCD based on $\NPE^r(\cdot,M_C)$ will be nondegenerate
(see \cite{ceyhan:Phd-thesis} and \cite{ceyhan:dom-num-NPE-SPL}).

The properties \textbf{P1}, \textbf{P2}, \textbf{P4}, \textbf{P5},
and \textbf{P7} follow by definition for all $M$ and $r$.
Furthermore \textbf{P9} holds, since $\NPE^r$ is geometry invariant for uniform data.
Property \textbf{P5} holds
with similarity ratio of $\NPE^r(x,M)$ to $\TY$:
$\displaystyle \frac{\min(d(v(x),\,e(x)),r\,d(v(x),\, \ell(v(x),x)))}{d(v(x),\,e(x))}$;
that is, $\NPE^r(x,M)$ is similar to $\TY$ with the given ratio.
\textbf{P6} holds depending on the pair $M$ and $r$.
That is, there exists an $r_o:=r_o(M)$ so that $\NPE^{r_o}(x,M)$
satisfies \textbf{P6} for all $r \le r_o(M)$, and fails
to satisfy otherwise.
\textbf{P6} fails for all $M$ when $r=\infty$.
\textbf{P8} holds only when $M=M_C$.
With $CM$-vertex regions, for all $r \in [1,\infty]$,
the area $A\left(\NPE^r(x,M_C) \right)$ is a continuous
function of $d(\ell_r(v(x),x),v(x))$ which is a
continuous function of $d(\ell(v(x),x),v(x))$ which is a
continuous function of $x$.

Moreover, $\Lambda_0\left(\NPE^r,M\right)=\Y_3$ for all $r \in [1,\infty]$
and $M \in \R^2 \setminus \Y_3$, since the $\R^2$-Lebesgue measure
$\lambda(\NPE^r(x,M))=0$ iff $x\in \Y_3$.

\begin{figure}[ht]
\begin{center}
\psfrag{A}{\scriptsize{$\y_1$}}
\psfrag{B}{\scriptsize{$\y_2$}}
\psfrag{C}{\scriptsize{$\y_3$}}
\psfrag{CC}{\scriptsize{$M_{CC}$}}
\psfrag{x}{\scriptsize{$x$}}
\psfrag{M1}{\scriptsize{$M_3$}}
\psfrag{M2}{\scriptsize{$M_1$}}
\psfrag{M3}{\scriptsize{$M_2$}}
\epsfig{figure=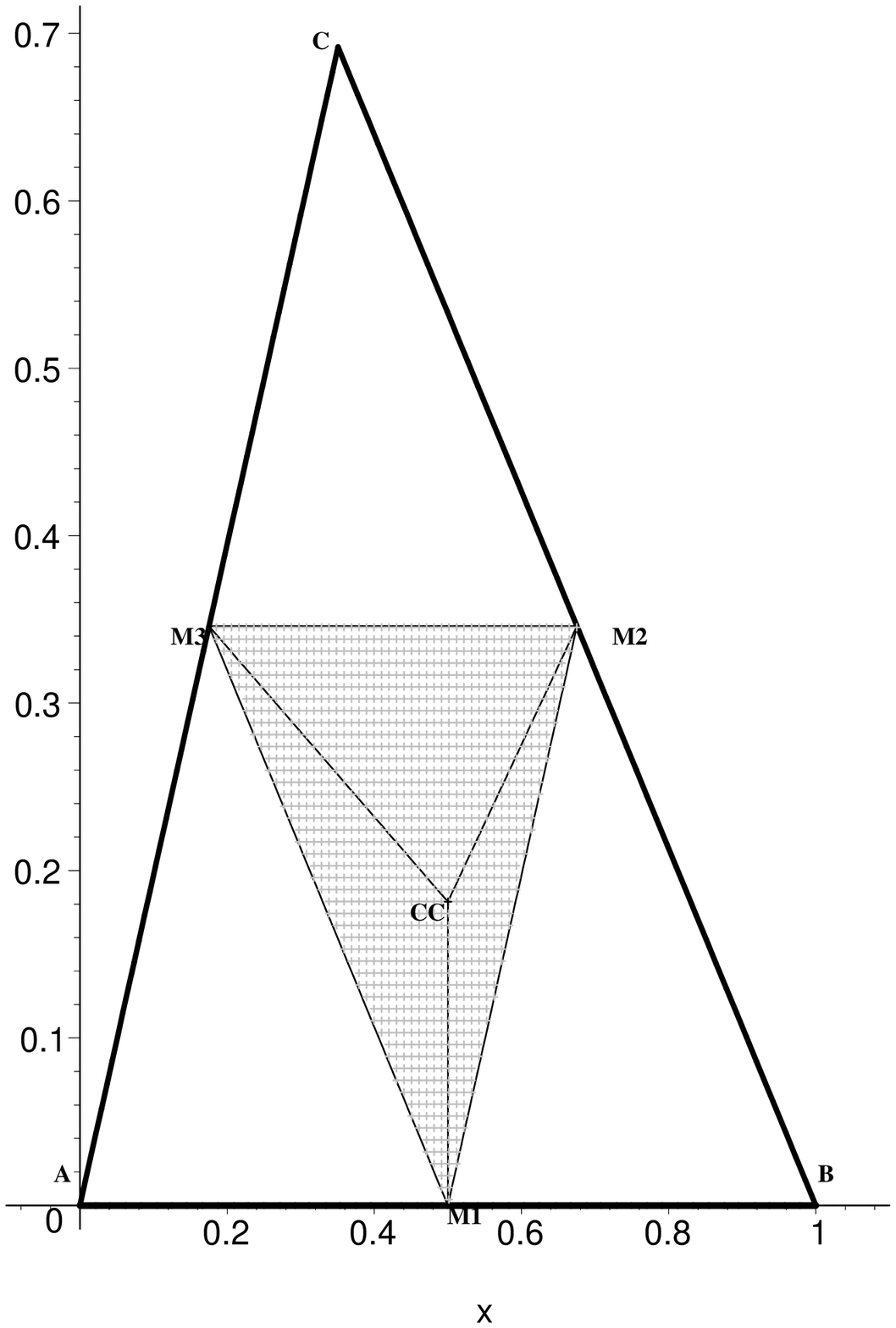,  height=140pt , width=200pt}
\psfrag{IC}{\scriptsize{$M_I$}}
\psfrag{P1}{\scriptsize{$P^{IC}_3$}}
\psfrag{P2}{\scriptsize{$P^{IC}_1$}}
\psfrag{P3}{\scriptsize{$P^{IC}_2$}}
\psfrag{Q2}{}
\psfrag{Q3}{}
\epsfig{figure=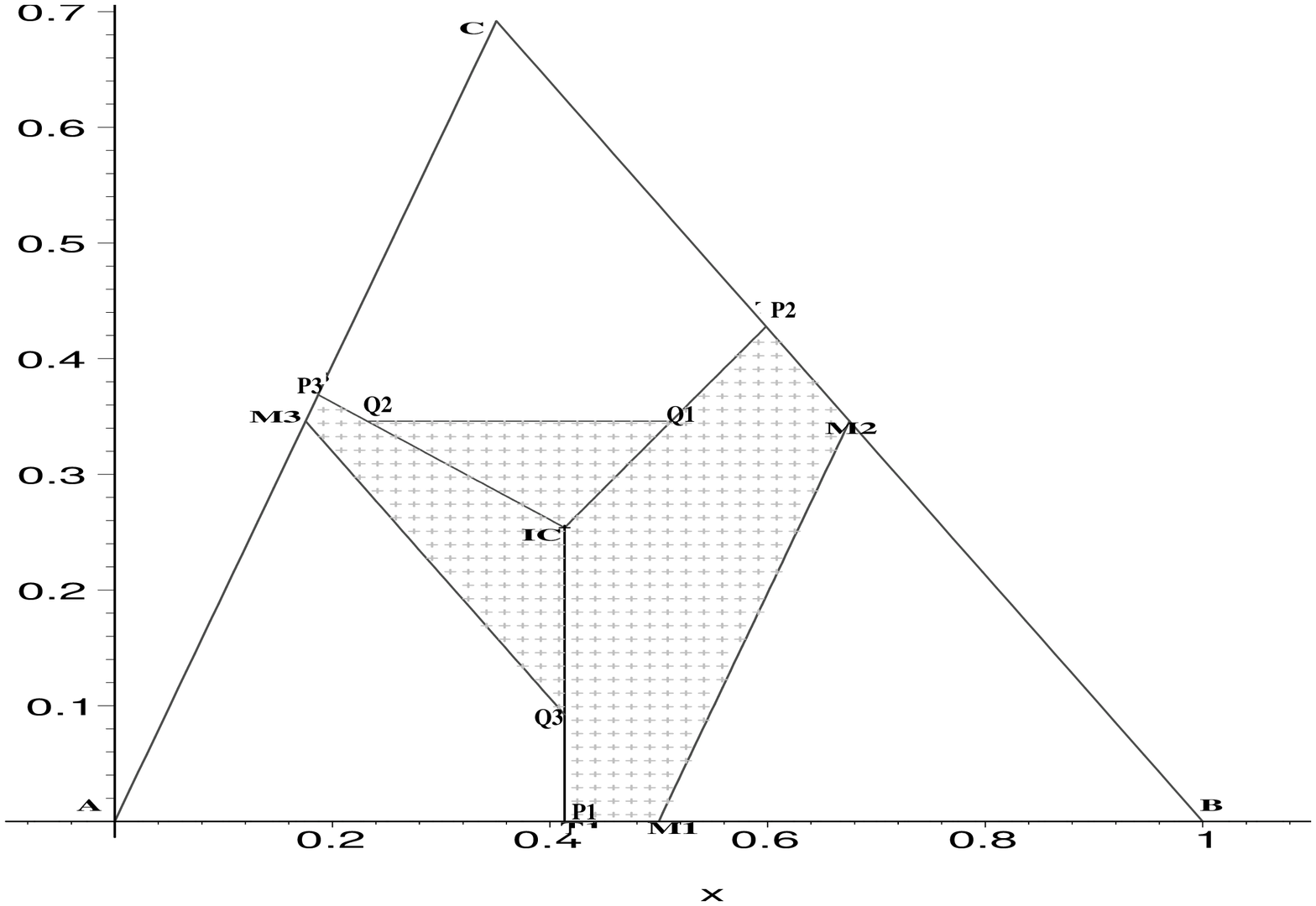, height=140pt , width=200pt}
\end{center}
\caption{Superset region $\RS\left(\NPE^2,M_{CC} \right)$ in an acute
triangle (left), superset region, $\RS^{\perp}(\NPE^2,M_I)$ (right)
}
\label{fig:superset-NDE}
\end{figure}

As for \textbf{P3}, for $T_2(x) \subseteq \TY$
one can loosen the concept of center by treating the line
$\ell(v(x),x)$ as the {\em edge-wise central line},
so \textbf{P3} is satisfied in this loose sense for $r=2$.
Notice that $x$ is not the unique center in this sense
but a point on a central line.
Let $M_i$, $i \in \{1,2,3\},$ be the midpoints of the edges of $\TY$,
and $T(M_1,M_2,M_3)$ be triangle whose vertices are these midpoints.
Then for any $x \in T(M_1,M_2,M_3)$, $\NPE^2(x,M) = \TY$, so
$T(M_1,M_2,M_3) \subseteq \RS(\NPE^2,M)$ where equality
holds for $M=M_C$ for all triangles and for $M=M_{CC}$ in
non-obtuse triangles (see Figure \ref{fig:superset-NDE} (left)).

\begin{figure}[ht]
\begin{center}
\psfrag{A}{\scriptsize{$\y_1$}}
\psfrag{B}{\scriptsize{$\y_2$}}
\psfrag{C}{\scriptsize{$\y_3$}}
\psfrag{CC}{\scriptsize{$M_{CC}$}}
\psfrag{x}{\scriptsize{$x$}}
\psfrag{E}{\scriptsize{}}
\psfrag{F}{\scriptsize{}}
\psfrag{G}{\scriptsize{}}
\psfrag{H}{\scriptsize{}}
\epsfig{figure=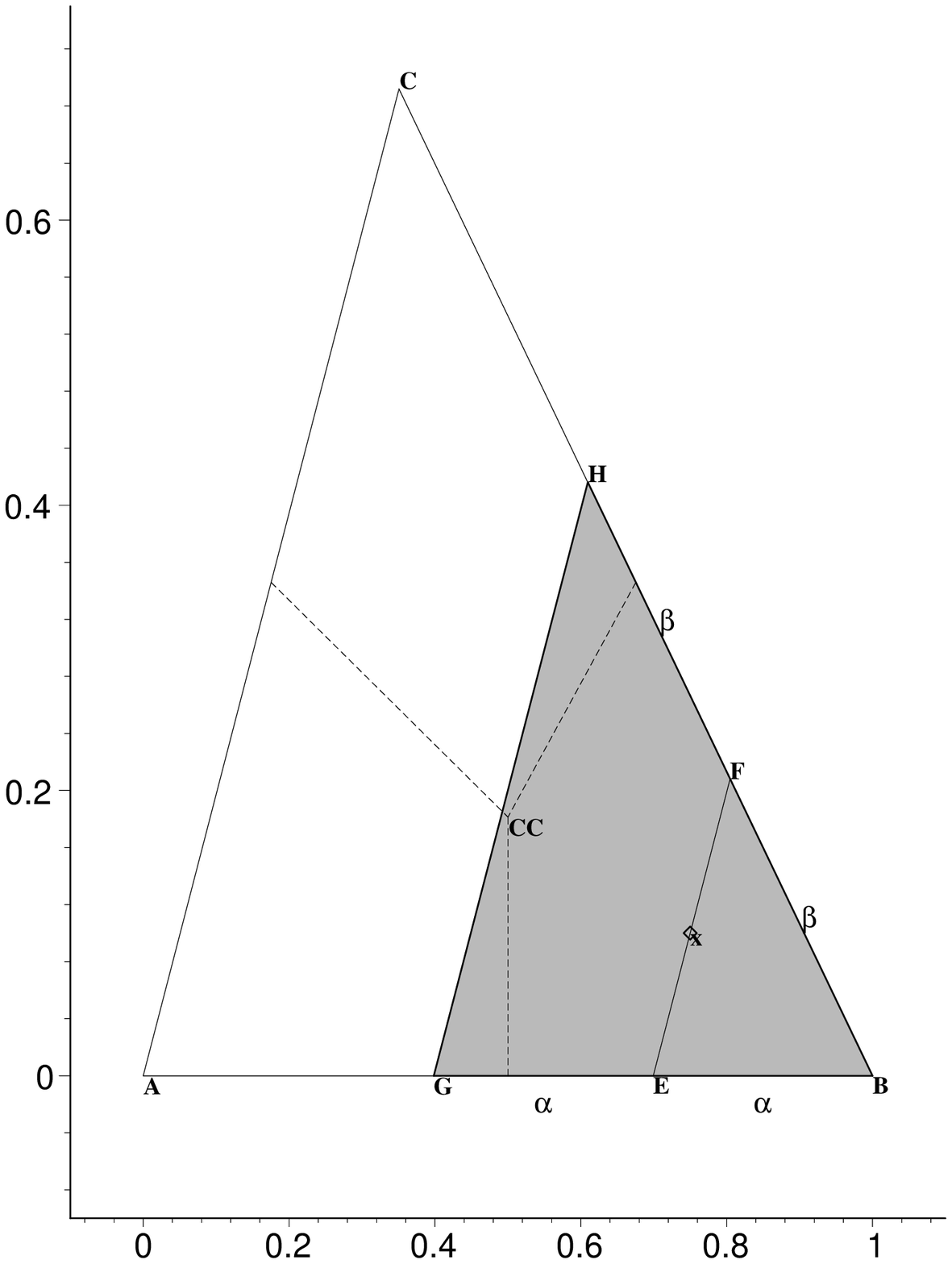,  height=140pt , width=200pt}
\psfrag{IC}{\scriptsize{$M_I$}}
\psfrag{S1}{}
\psfrag{S2}{}
\epsfig{figure=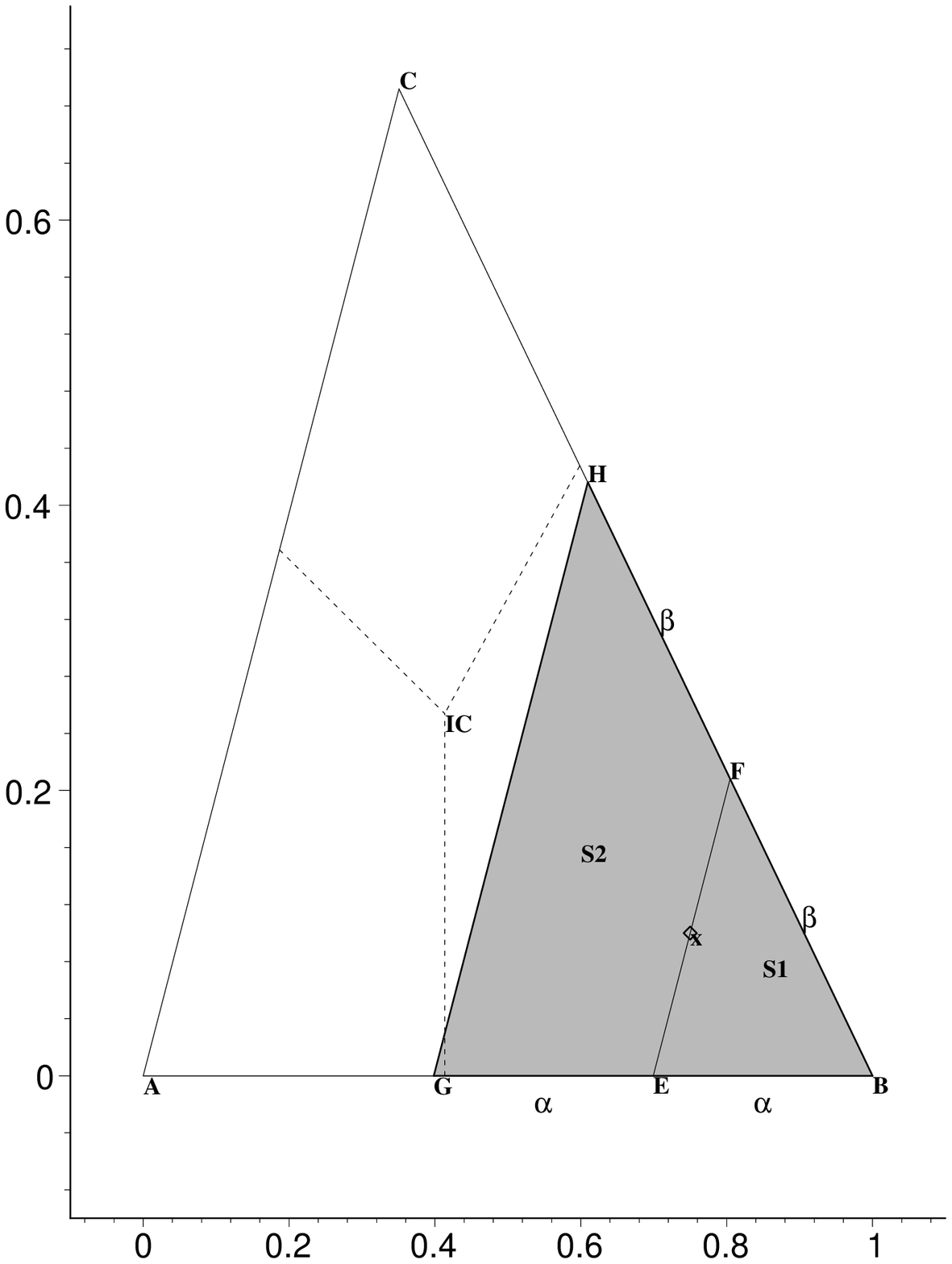, height=140pt , width=200pt}
\psfrag{S2}{}
\epsfig{figure=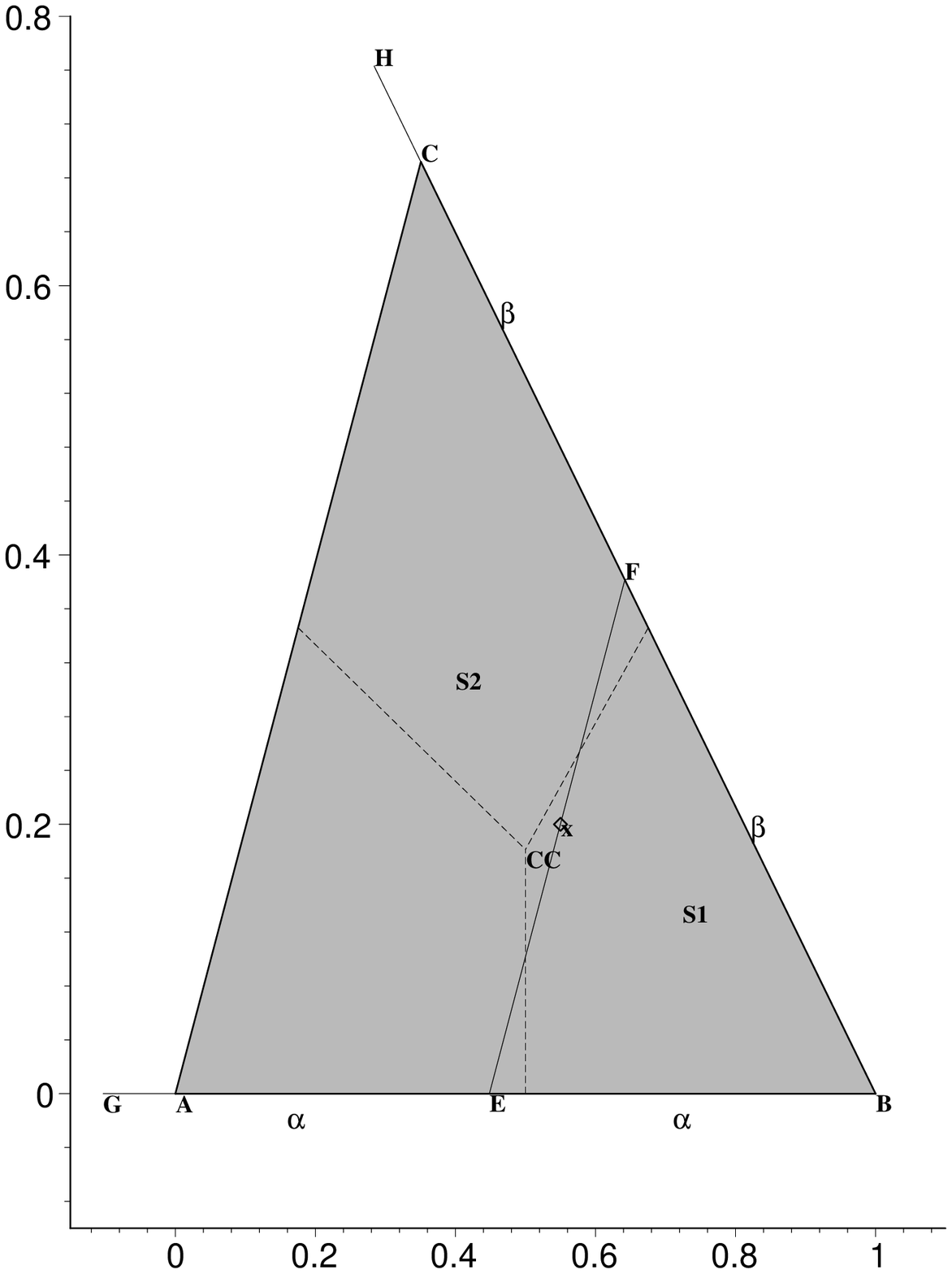, height=140pt, width=200pt}
\end{center}
\caption{
\label{fig:NDE-with-CC-IC}
Shaded regions are double-edge proximity regions $\NPE^2(x,M) \subsetneq \TY$
with $M=M_{CC}$ for an $x \in R^{\perp}_{CC}(\y_2)$ (top left),
with $M=M_I$ for an $x \in R^{\perp}_{IC}(\y_2)$ (top right).
Depicted in the bottom is an example of $\NPE^2(x,M)= \TY$
with $M=M_{CC}$ for an $x \in R^{\perp}_{CC}(\y_2)$.
}
\end{figure}

For an example of double-edge proximity regions $\NPE^2(x,M)$  with
$CC$-vertex regions with orthogonal projections, see Figure \ref{fig:NDE-with-CC-IC} (top left).
Notice that we use the vertex closest to $x$; i.e., $\argmin_{\y \in \Y_3}d(x,\y)$
for $\NPE^2(x,CC)$, i.e. vertex regions $R_{CC}(\cdot)$.
Furthermore, if $x$ is close enough to $M$,
it is possible to have $\NPE^2(x,M) = \TY$.
See Figure \ref{fig:NDE-with-CC-IC} (bottom) for an example
with $CC$-vertex regions with orthogonal projections.

In non-obtuse triangles, $\RS\left(\NPE^2,M_{CC}\right)=T(M_1,M_2,M_3)$.
But, in obtuse triangles, $\RS\left(\NPE^2,M_{CC}\right) \supsetneq T(M_1,M_2,M_3)$
and is a quadrilateral.
The functional form of the superset region, $\RS\left(\NPE^r,M\right)$,
in $T_b$ is given by
\begin{multline*}
\RS\left(\NPE^r,M\right)=\left \{(x,y) \in R_{M}(\y_1): y \ge  \frac{c_2\,(1-r\,x)}{r\,(1-c_1)}\right\}\bigcup \\
\left \{(x,y) \in R_{M}(\y_2): y \ge \frac{c_2\,(r\,(x-1)+1)}{r\,c_1} \right\}
\bigcup \left \{(x,y) \in R_{M}(\y_3): y \le c_2\frac{r-1}{r} \right\},
\end{multline*}
and
the functional form of $T(M_1,M_2,M_3)$ in $T_b$ is given by
$$T(M_1,M_2,M_3)=\left \{(x,y) \in T_b: y \le \frac{c_2}{2};\;
y \ge \frac{c_2\,(-1+2\,x)}{2\,c_1};\; y \ge  \frac{c_2\,(1-2\,x)}{2\,(1-c_1)} \right \}.$$

Let $\RS^{\perp}\left(\NPE^r,M\right)$ be the superset region for $\NPE^r$
based on $M$-vertex regions with orthogonal projections.
See Figure \ref{fig:superset-NDE}
for the superset region $\RS^{\perp}\left(\NPE^2,M_I \right)$.
Again $T(M_1,M_2,M_3) \subseteq \RS^{\perp}(\NPE^2,M_I)$ for all $\TY$
with equality holding when $\TY$ is an equilateral triangle.
For $\NPE^2(\cdot,M_C)$ constructed using the median lines
$\RS\left(\NPE^2,M_C \right)=T(M_1,M_2,M_3)$ and for
$\NPE^2(\cdot,M_C)$ constructed by the orthogonal projections,
$\RS^{\perp} \left(\NPE^2,M_C \right) \supseteq T(M_1,M_2,M_3)$
with equality holding when $\TY$ is an equilateral triangle.
An example of double-edge proximity regions is given in Figure
\ref{fig:NDE-with-CC-IC} (top right) where $IC$-vertex regions with orthogonal
projections to the edges is used.
We could also use $IC$-vertex
regions obtained by inner angle bisectors.
Note also that the superset region $\RS^{\perp}(\NPE^2,M_I)$ is as in
Figure \ref{fig:superset-NDE}.
Again $T(M_1,M_2,M_3) \subseteq \RS^{\perp}(\NPE^2,M_I)$ for all $\TY$ and
$T(M_1,M_2,M_3)=\RS^{\perp}(\NPE^2,M_I)$ iff $\TY$ is an
equilateral triangle.

\begin{figure}[ht]
\begin{center}
\psfrag{A}{\scriptsize{$\y_1$}}
\psfrag{B}{\scriptsize{$\y_2$}}
\psfrag{C}{\scriptsize{$\y_3$}}
\psfrag{CC}{\scriptsize{$M_{CC}$}}
\psfrag{x}{\scriptsize{$x$}}
\psfrag{E}{\scriptsize{}}
\psfrag{F}{\scriptsize{}}
\psfrag{G}{\scriptsize{}}
\psfrag{H}{\scriptsize{}}
\epsfig{figure=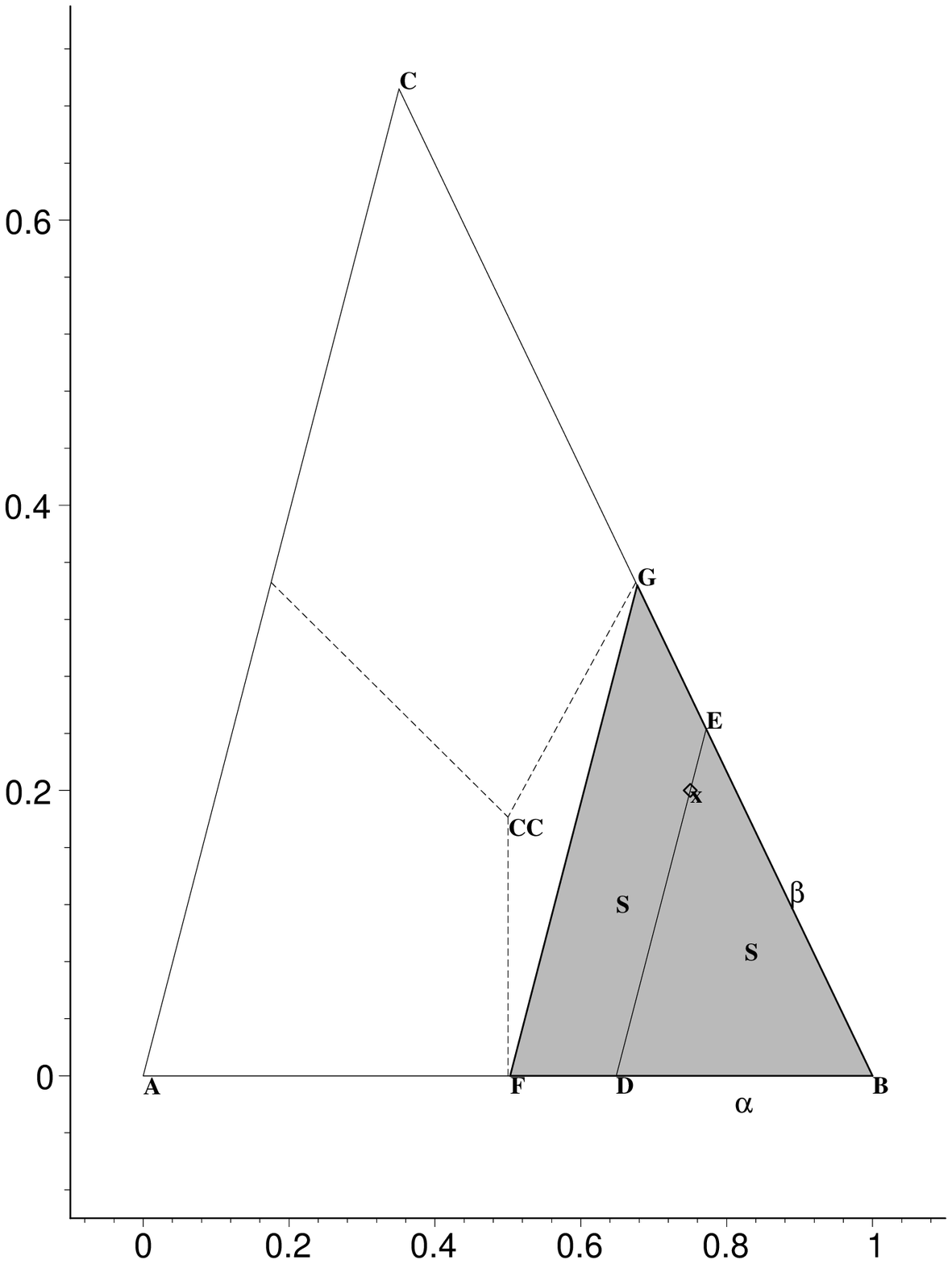,  height=140pt , width=200pt}
\psfrag{IC}{\scriptsize{$M_I$}}
\psfrag{S1}{}
\psfrag{S2}{}
\epsfig{figure=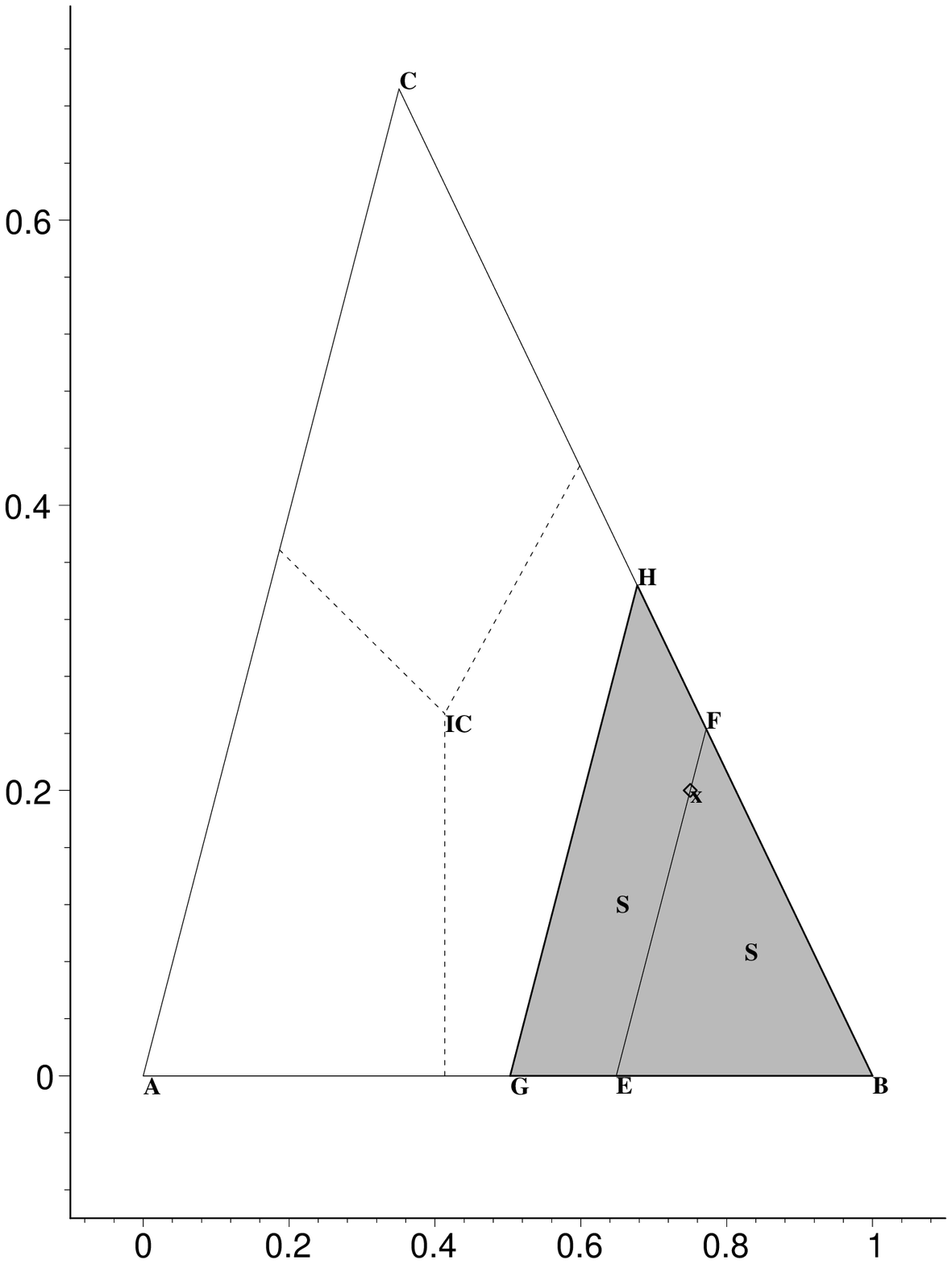, height=140pt , width=200pt}
\psfrag{G}{}
\psfrag{H}{}
\psfrag{S}{}
\epsfig{figure=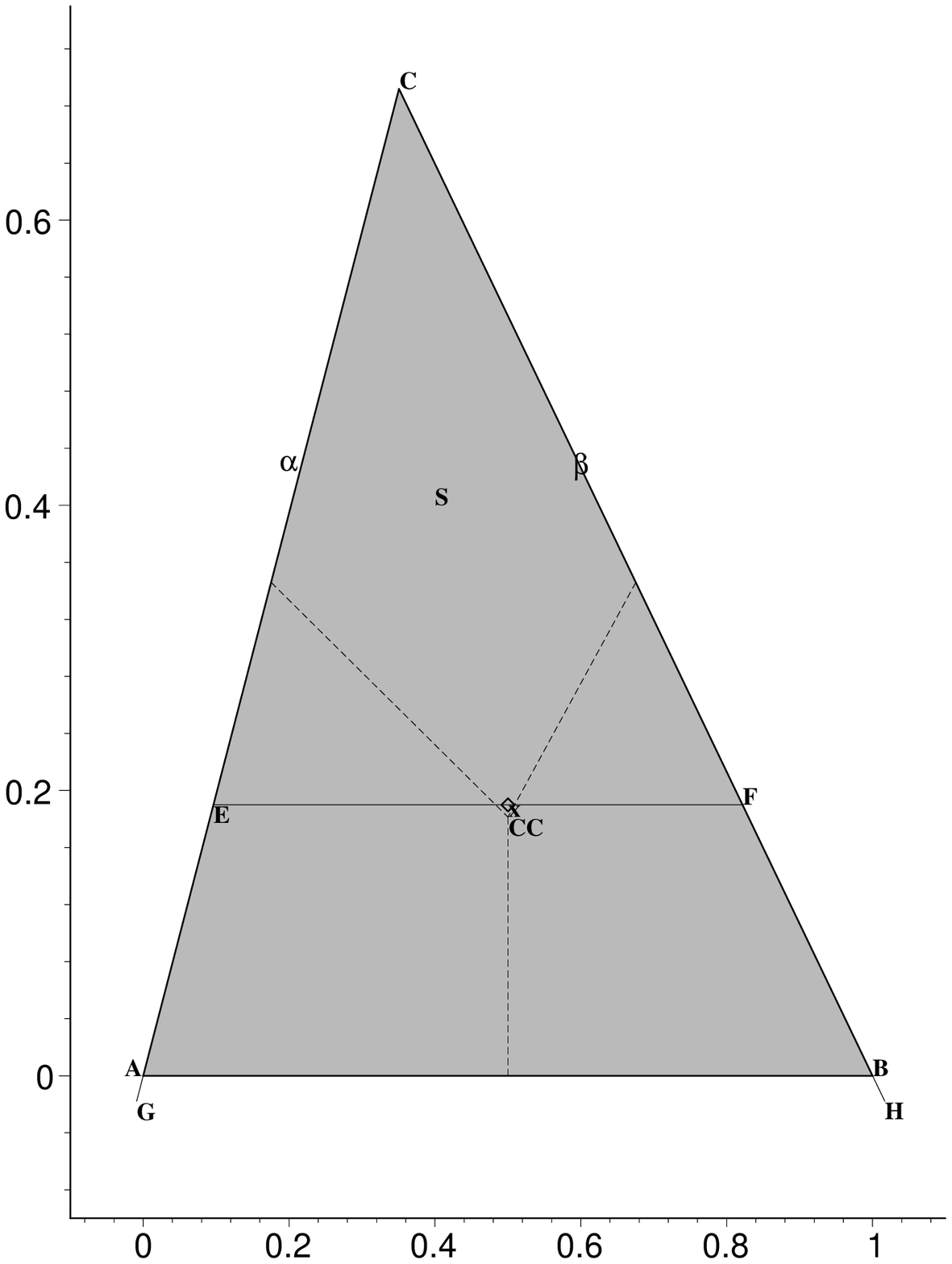, height=140pt , width=200pt}
\end{center}
\caption{
Shaded regions are double-area proximity regions $\NPE^{\sqrt{2}}(x,M) \subsetneq \TY$
with $M=M_{CC}$ for an $x \in R^{\perp}_{CC}(\y_2)$ (top left),
with $M=M_I$ for an $x \in R^{\perp}_{IC}(\y_2)$ (top right)
and $\NPE^{\sqrt{2}}(x,M)=\TY$ for an $x \in R^{\perp}_{CC}(\y_2)$ with $M=M_{CC}$ (bottom).
}
\label{fig:NDA-with-CC-IC}
\end{figure}

\begin{figure}[ht]
\begin{center}
\psfrag{A}{\scriptsize{$\y_1$}}
\psfrag{B}{\scriptsize{$\y_2$}}
\psfrag{C}{\scriptsize{$\y_3$}}
\psfrag{CC}{\scriptsize{$M_{CC}$}}
\psfrag{IC}{\scriptsize{$M_{I}$}}
\psfrag{x}{\scriptsize{$x$}}
\psfrag{D}{}
\psfrag{E}{}
\psfrag{F}{}
\psfrag{S}{}
\psfrag{Ab}{}
\psfrag{Ac}{}
\psfrag{Ba}{}
\psfrag{Bc}{}
\psfrag{Ca}{}
\psfrag{Cb}{}
\epsfig{figure=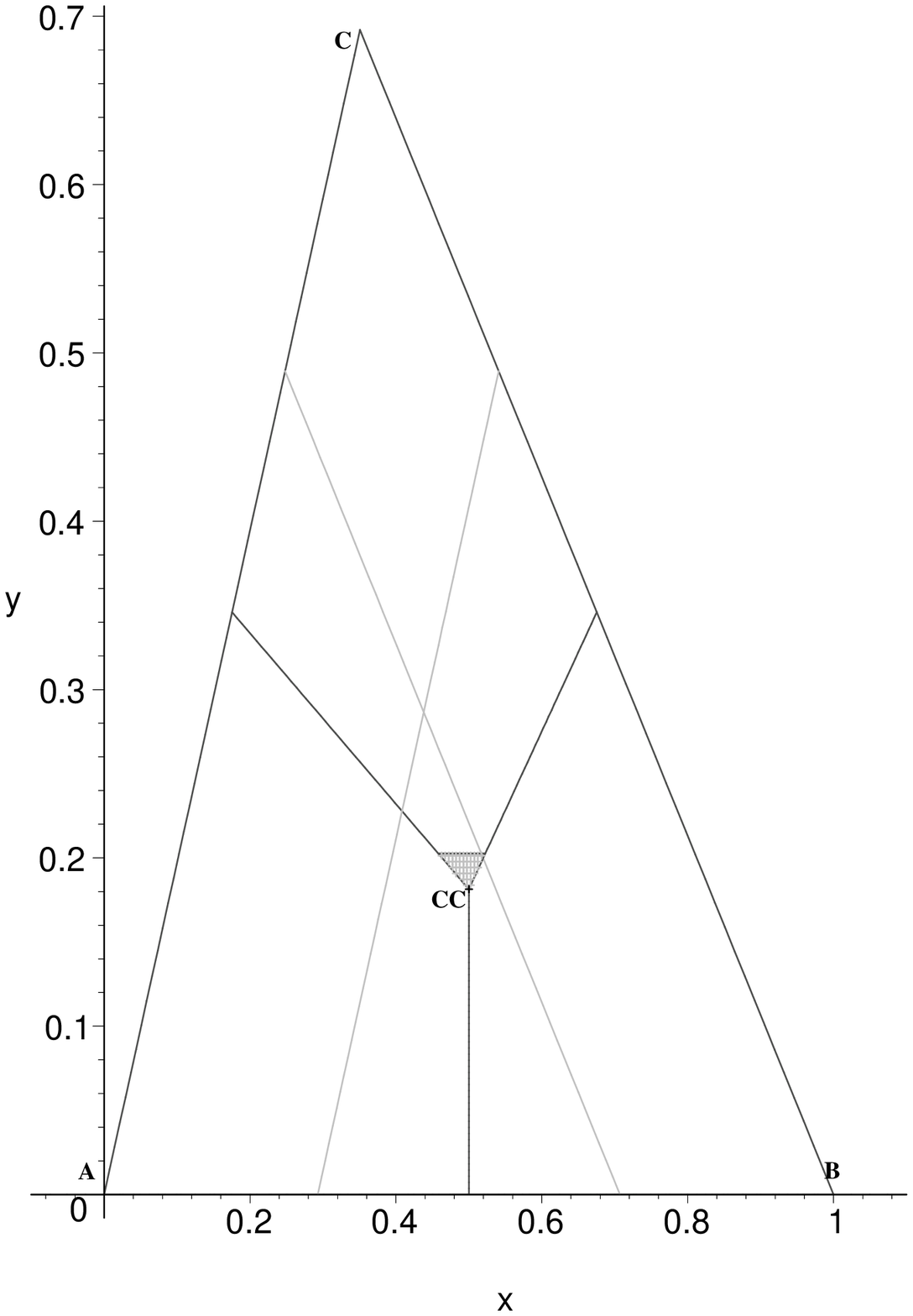, height=140pt , width=200pt}
\epsfig{figure=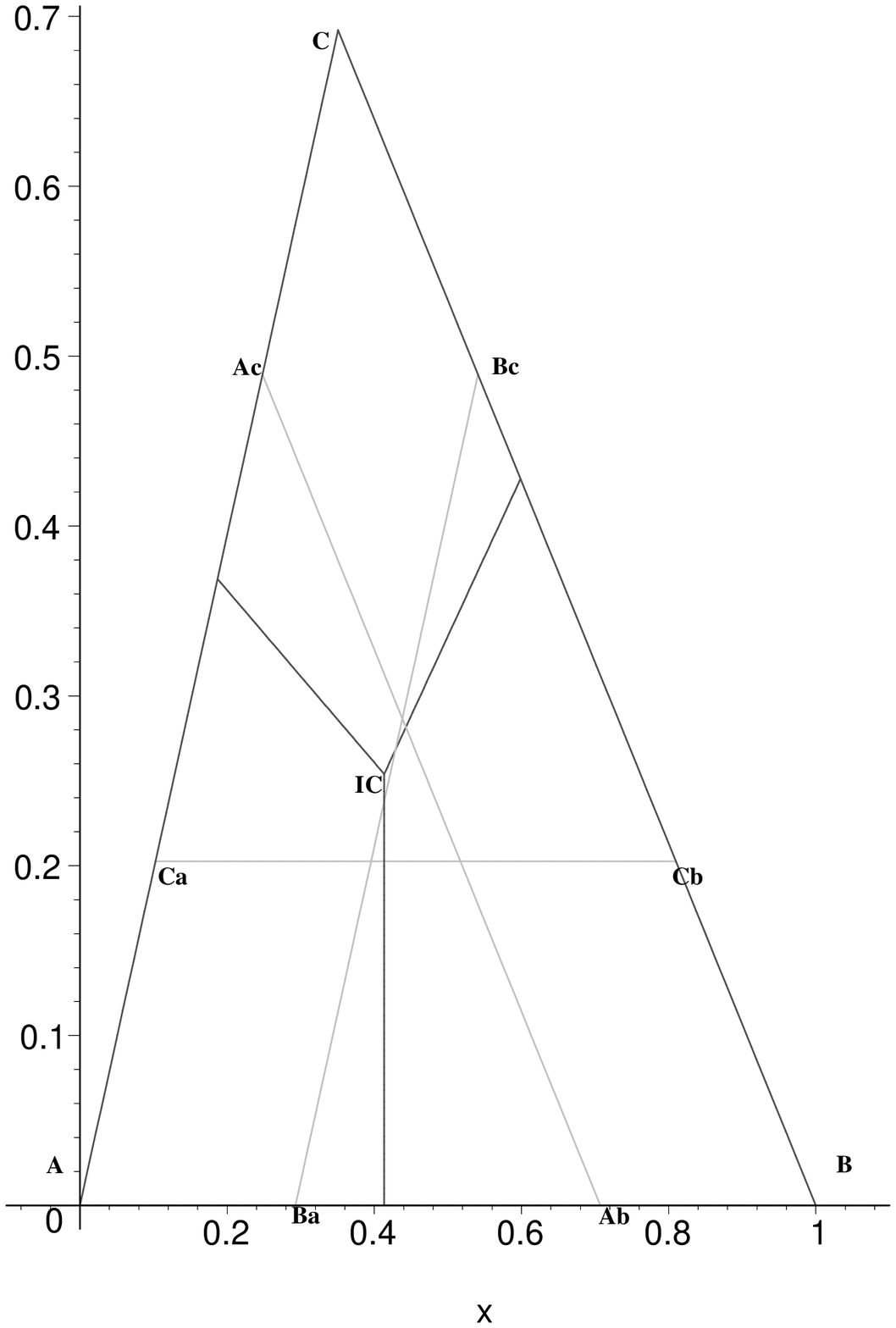, height=140pt , width=200pt}
\end{center}
\caption{
The superset region, $\RS\left( \NPE^{\sqrt{2}},M \right)$ with $M=M_{CC}$ (left)
and $M=M_I$ (right).
}
\label{fig:RS-NDA-CC-IC}
\end{figure}

For $r=\sqrt{2}$, one can loosen the concept of center by treating
the line $\ell(v(x),x)$ as the {\em area-wise central line} in
$\NPE^{\sqrt{2}}(x,M)$, so \textbf{P3} is satisfied in this loose sense.
For an example of $\NPE^{\sqrt{2}}(x,M)$ with $CC$-vertex regions with orthogonal projections,
see Figure \ref{fig:NDA-with-CC-IC} (top left).
$\RS\left(\NPE^{\sqrt{2}},M_{CC}\right)$ has positive area;
see Figure \ref{fig:RS-NDA-CC-IC}.
An example of double-area proximity region with $IC$-vertex regions
is given at Figure \ref{fig:NDA-with-CC-IC} (top right) with orthogonal projections to the edges.
Note that if $x$ is close enough to $M$,
it is possible to have $\NPE^{\sqrt{2}}(x,M) = \TY$.
See Figure \ref{fig:NDA-with-CC-IC} (bottom)
with $CC$-vertex regions with orthogonal projections.
We could also use $IC$-vertex regions obtained by inner
angle bisectors.

Note also that $\NPE^{\sqrt{2}}(x,M_I) = \TY$ might occur if $x$ is
close enough to $M_I$ when $M_I \not\in \mathscr T^r$.
Let $h_j$ be the altitude of $\TY$ at vertex $y_j$, for $j=1,2,3$.
If $r_{ic} < \left(\frac{\sqrt{2}-1}{\sqrt{2}}\right)h_j$ for some $j \in \{1,2,3\}$, then
$\RS\left(\NPE^{\sqrt{2}},M_I\right)$ has positive area.
See Figure \ref{fig:RS-NDA-CC-IC} where the superset region is barely noticeable.
If $r_{ic}  \ge \left( \frac{\sqrt{2}-1}{\sqrt{2}} \right) \max_{j \in \{1,2,3\}}\{h_j\}$,
then $\RS\left(\NPE^{\sqrt{2}},M_I\right)$ has zero area.
In $T_b$, $r_{ic}  \ge \left(\frac{\sqrt{2}-1}{\sqrt{2}}\right) \max\{h_1,h_3\}$ always hold,
but $r_{ic}  \ge \left(\frac{\sqrt{2}-1}{\sqrt{2}}\right)h_2$ holds iff
$c_1^2+c_2^2 < \left(\sqrt{2}-1\right)^2 \left(1+\sqrt{(1-c_1)^2+c_2^2}\right)^2$ iff
$|e_2|<\left(\sqrt{2}-1\right)(1+|e_1|)$.

\begin{figure}[ht]
\begin{center}
\psfrag{A}{\scriptsize{$\y_1$}}
\psfrag{B}{\scriptsize{$\y_2$}}
\psfrag{C}{\scriptsize{$\y_3$}}
\psfrag{CC}{\scriptsize{$M_{CC}$}}
\psfrag{IC}{\scriptsize{$M_{I}$}}
\psfrag{x}{\scriptsize{$x$}}
\psfrag{D}{}
\psfrag{E}{}
\psfrag{F}{}
\psfrag{S}{}
\psfrag{Ab}{}
\psfrag{Ac}{}
\psfrag{Ba}{}
\psfrag{Bc}{}
\psfrag{Ca}{}
\psfrag{Cb}{}
\epsfig{figure=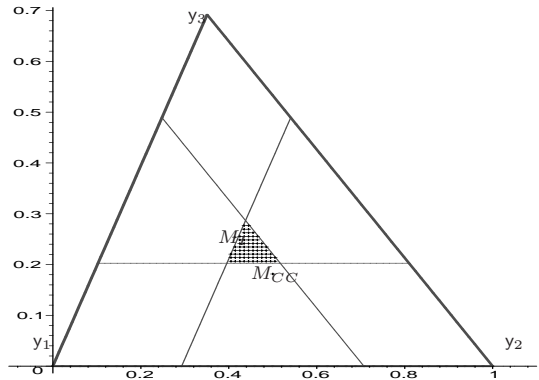, height=140pt , width=200pt}
\end{center}
\caption{
The triangle $\mathscr T^{r=\sqrt{2}}$.
%(left) and the superset region, $\RS\left(\NPE^{\sqrt{2}},M_{CC} \right)$ (right).
}
\label{fig:Tr-RS-NDA-CC}
\end{figure}

In $\TY$, drawing the lines $q_i(r,x)$ such that
$d(\y_i,e_i)=r\,d(q_i(r,x),\y_i)$ for $i \in \{1,2,3\}$  yields a triangle,
$\Tr$, for $r<3/2$ .
See Figure \ref{fig:Tr-RS-NDA-CC} for $\Tr$ with $r=\sqrt{2}$.
The functional form of $\Tr$ in $T_b$ is
{\small
\begin{align}
\label{eqn:T^r-def}
& \Tr = \left \{(x,y) \in T_b: y \ge \frac{c_2\,(r-1)}{r};\;
y \le \frac{c_2\,(1-r\,x)}{r\,(1-c_1)};\; y \le \frac{c_2\,(r\,(x-1)+1)}{r\,c_1} \right\}=
T\left(t_1(r),t_2(r),t_3(r)\right)\\
&=T\Biggl( \left(\frac{(r-1)\,(1+c_1)}{r},\frac{c_2\,(r-1)}{r} \right),
\left(\frac{2-r+c_1\,(r-1)}{r},\frac{c_2\,(r-1)}{r} \right),
\left(\frac{c_1\,(2-r)+r-1}{r},\frac{c_2\,(r-2)}{r} \right) \Biggr) \nonumber.
\end{align}
}
There is a crucial difference between $\Tr$ and $T(M_1,M_2,M_3)$:
$T(M_1,M_2,M_3) \subseteq \RS\left(\NPE^r,M\right)$ for all $M$ and $r \ge 2$,
but $(\Tr)^o$ and $\RS\left(\NPE^r,M\right)$ are disjoint regions
%(i.e., $(\Tr)^o \cap \RS\left(\NPE^r,M\right) =\emptyset$)
for all $M$ and $r$.
So if $M \in (\Tr)^o$, then $\RS\left(\NPE^r,M\right)=\emptyset$;
if $M \in \partial(\Tr)$, then $\RS\left(\NPE^r,M\right)=\{M\}$;
and if $M \not\in \Tr$, then $\RS\left(\NPE^r,M\right)$ has positive area.
Thus $\NPE^r(\cdot,M)$ fails to satisfy \textbf{P6} if $M \not\in \Tr$.
The triangle $\Tr$ defined above plays a crucial role in the analysis of the distribution
of the domination number of the proportional-edge PCD.
In fact, it has been shown that for $M \in \{t_1(r),t_2(r),t_3(r)\}$
there exists a specific value of $r$ for which the asymptotic distribution
of the domination number is non-degenerate (\cite{ceyhan:dom-num-NPE-MASA}).
The superset region $\RS\left(\NPE^r,M\right)$ will be important for both
the domination number and the relative density of the corresponding PCDs.

In non-acute triangles, the above condition holds for $M=M_{CC}$,
since in right and obtuse triangles,
$M_{CC} \notin \TY^o \Rightarrow M_{CC} \notin \mathscr T^r$
(since $\TY \supset \mathscr T^r$).
For an acute basic triangle,
if $\displaystyle y_{cc}<\frac{c_2\,\left( \sqrt{2}-2\,{x_{cc}}\right)}{2\,(1-c_1)}$ holds,
then $M_{CC} \notin \mathscr T^{r=\sqrt{2}}$ and the superset region for such
triangles is the triangle with vertices $M_{CC}$,
$\displaystyle \left(\frac{-c_2^2+c_1^2-1+c_2^2\sqrt{2}}{2\,(-1+c_1)},c_2-\frac{c_2}{\sqrt{2}} \right)$,
and
$\displaystyle \left(\frac{-c_2^2+c_1^2+c_2^2\sqrt{2}}{2\,c_1},c_2-\frac{c_2}{\sqrt{2}} \right)$.

\begin{remark}
In terms of the properties stated in Section \ref{sec:prox-maps},
$\NPE^{3/2}(\cdot,M_C)$ is the most appealing proximity map in the family
$\mathscr N_{PE}^r:=
\left \{\NPE^r(\cdot,M): r \in [1,\infty],\; M \in \R^2 \setminus \Y_3 \right \}$.
It is also noteworthy that the asymptotic distribution of the domination
number of the PCD based on $\NPE^{3/2}(\cdot,M_C)$ is nondegenerate.
$\square$
\end{remark}

\begin{remark}
\begin{itemize}
\item[]
\item
For $r_1 \le r_2$, $\NPE^{r_1}(x,M) \subseteq \NPE^{r_2}(x,M)$ for all $x \in \TY$.
For $r_1<r_2$, $\NPE^{r_1}(x,M) \subseteq \NPE^{r_2}(x,M)$ with equality holding
for only $ x\in \Y_3$ or $x \in \RS\left(\NPE^{\min(r_1,r_2)},M\right)$.
\item
For $3/2 < r_1 <r_2 $, $\RS\left(\NPE^{r_1},M\right) \subsetneq \RS\left(\NPE^{r_2},M\right)$ and
$\RS\left(\NPE^{r},M\right)=\emptyset$ for $r<3/2$.
\item
For $r_1 < r_2$, $A\left(\NPE^{r_1}(X,M)\right) \le ^{ST} A\left(\NPE^{r_2}(X,M)\right)$ for $X$
from a continuous distribution on $\TY$ where $\le ^{ST}$ stands
for ``stochastically smaller than". $\square$
\end{itemize}
\end{remark}

%\begin{figure}[ht]
%\centering
%%\rotatebox{-90}{ \resizebox{2.0 in}{!}{\includegraphics{DelaunayXY.ps} }}
%%\rotatebox{-90}{ \resizebox{2.0 in}{!}{\includegraphics{DelaunayCH.ps} }}
%\rotatebox{-90}{ \resizebox{2.0 in}{!}{\includegraphics{arcs_PE.ps} }}
%\caption{
%\label{fig:PE-arcs-multiT}
%%A realization of 10 $\Y$ points generated iid $\U(0,1)\times(0,1)$
%%and the convex hull $\C_H\left(\Y_{10}\right)$ of these $\Y$ points and 200 $\X$
%%points generated iid $\U(0,1)\times(0,1)$ (left),
%%77 $\X$ points that lie in the $\C_H\left(\Y_{10}\right)$ (middle),
%%and the corresponding arcs for $\NPE^{r=2}(x,M_C)$ (right).
%The arcs for $\NPE^{r=2}(x,M_C)$ for the 77 $\X$ points that lie in the $\C_H\left(\Y_{10}\right)$
%in Figure \ref{fig:AS-arcs-multiT}.
%}
%\end{figure}

\subsubsection{Extension of $\NPE^r$ to Higher Dimensions}
\label{sec:NYr-higher-D}
The extension to $\R^d$ for $d > 2$ is straightforward.
The extension with $M=M_C$ is given her,
but the extension for general $M$ is similar.
Let $\Y_{d+1} = \{\y_1,\y_2,\cdots,\y_{d+1}\}$ be $d+1$ points
that do not lie on the same $d-1$-dimensional hyperplane.
Denote the simplex formed by these $d+1$ points as $\mathfrak S (\Y_{d+1})$.
A simplex is the simplest polytope in $\R^d$
having $d+1$ vertices, $d\,(d+1)/2$ edges and $d+1$ faces of dimension $(d-1)$.
For $r \in [1,\infty]$, define the proximity map as follows.
Given a point $x$ in $\mathfrak S (\Y_{d+1})$,
let $v := \argmin_{\y \in \Y_{d+1}} \mbox{V}(Q_{\y}(x))$
where $Q_{\y}(x)$ is the polytope with vertices
being the $d\,(d+1)/2$ midpoints of the edges,
the vertex $v$ and $x$ and $\mbox{V}(\cdot)$ is the
$d$-dimensional volume functional.
That is, the vertex region for vertex $v$ is the polytope with vertices
given by $v$ and the midpoints of the edges.
Let $v(x)$ be the vertex in whose region $x$ falls.
If $x$ falls on the boundary of two vertex regions,
$v(x)$ is assigned arbitrarily.
Let $\varphi(x)$ be the face opposite to vertex $v(x)$,
and $\Upsilon(v(x),x)$ be the hyperplane parallel to $\varphi(x)$ which contains $x$.
Let $d(v(x),\Upsilon(v(x),x))$ be the (perpendicular)
Euclidean distance from $v(x)$ to $\Upsilon(v(x),x)$.
For $r \in [1,\infty)$, let $\Upsilon_r(v(x),x)$ be
the hyperplane parallel to $\varphi(x)$
such that
\begin{gather*}
d(v(x),\Upsilon_r(v(x),x))=r\,d(v(x),\Upsilon(v(x),x))\\
\text{ and }\\
d(\Upsilon(v(x),x),\Upsilon_r(v(x),x))< d(v(x),\Upsilon_r(v(x),x)).
\end{gather*}
Let $\mathfrak S_r(x)$ be the polytope similar to and with the same
orientation as $\mathfrak S$ having $v(x)$
as a vertex and $\Upsilon_r(v(x),x)$ as the opposite face.
Then the proximity region
$\NPE^r(x,M_C):=\mathfrak S_r(x) \cap \mathfrak S(\Y_{d+1})$.
Notice that $r \ge 1$ implies $x \in \NPE^r(x,M_C)$.

\subsection{Central Similarity Proximity Maps}
\label{sec:NCS-region}
The other type of triangular proximity map introduced is the
central similarity proximity map.
This will turn out to be the most appealing proximity map in terms of
the properties in Section \ref{sec:prox-maps}.
Furthermore, the relative density of the corresponding
PCD will have mathematical tractability (\cite{ceyhan:arc-density-CS}).
Alas, the distribution of the domination number of the associated PCD
is still an open problem (\cite{ceyhan:Phd-thesis}).

\begin{figure}[ht]
\centering
\scalebox{.45}{\input{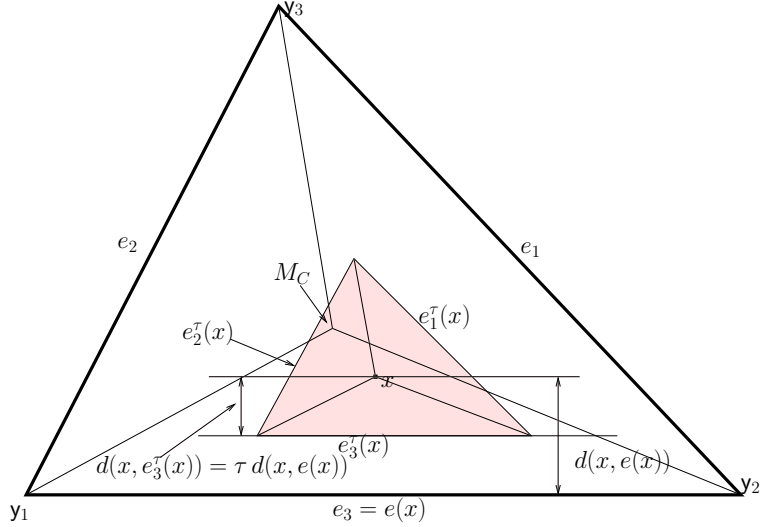}}
\caption{
\label{fig:ProxMapDefCS}
Construction of central similarity proximity region, $\NCS^{\tau=1/2}(x,M_C)$ (shaded region).}
\end{figure}

For $\tau \in [0,1]$,
define $\NCSt(\cdot,M):=N\left(\cdot,M;\tau,\Y_3\right)$ to be the
\emph{central similarity proximity map} with $M$-edge regions as follows;
see also Figure \ref{fig:ProxMapDefCS} with $M=M_C$.
For $x \in \TY\setminus \Y_3$, let $e(x)$ be the
edge in whose region $x$ falls; i.e., $x \in R_M(e(x))$.
If $x$ falls on the boundary of two edge regions,
$e(x)$ is assigned to $x$ arbitrarily.
For ${\tau} \in(0,1]$, the central similarity proximity region
$\NCSt(x,M)$ is defined to be the triangle
$T_{\tau}(x)$ with the following properties:
\begin{itemize}
\item[(i)] $T_{\tau}(x)$ has edges $e^{\tau}_i(x)$ parallel to $e_i$
for $i \in \{1,2,3\}$, and for $x \in R_M(e(x))$,
$d(x,e^{\tau}(x))=\tau\, d(x,e(x))$ and $d(e^{\tau}(x),e(x)) \le d(x,e(x))$
where $d(x,e(x))$ is the Euclidean (perpendicular) distance from $x$ to $e(x)$;
\item[(ii)] $T_{\tau}(x)$ has the same orientation as and is similar to $\TY$;
\item[(iii)] $x$ is the same type of center of $T_{\tau}(x)$ as $M$ is of $\TY$.
\end{itemize}
Note that (i) implies the parametrization of the PCD,
(ii) explains ``similarity", and
(iii) explains ``central" in the name,
{\em central similarity proximity map}.
For $\tau=0$, $\NCS^{\tau=0}(x,M):=\{x\}$ for all $x \in \TY$.
For $x \in \partial(\TY)$,
$\NCSt(x,M):=\{x\}$ for all $\tau \in [0,1]$.

\begin{figure}[ht]
\centering
\rotatebox{-90}{ \resizebox{2.0 in}{!}{\includegraphics{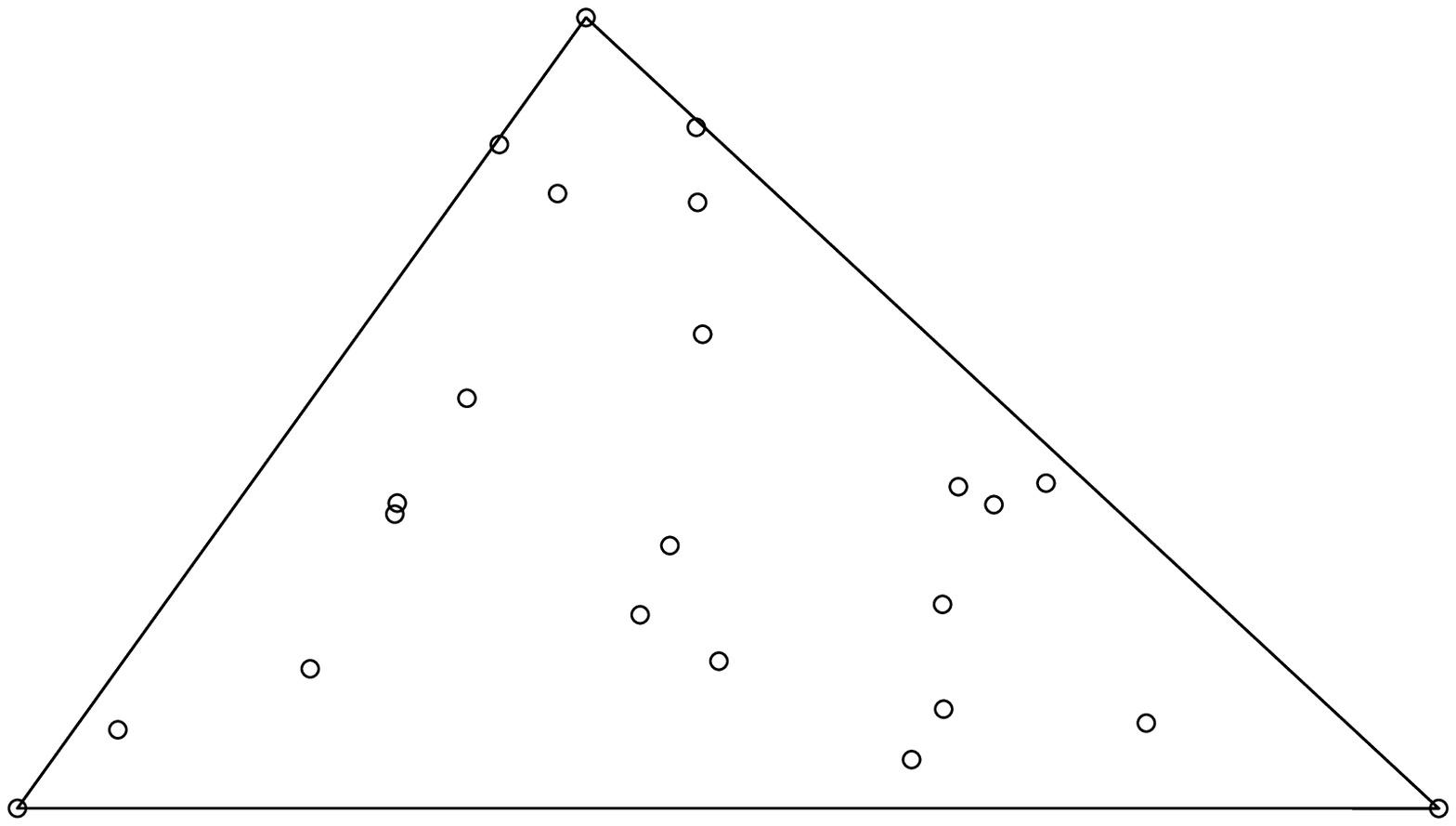} }}
\rotatebox{-90}{ \resizebox{2.0 in}{!}{\includegraphics{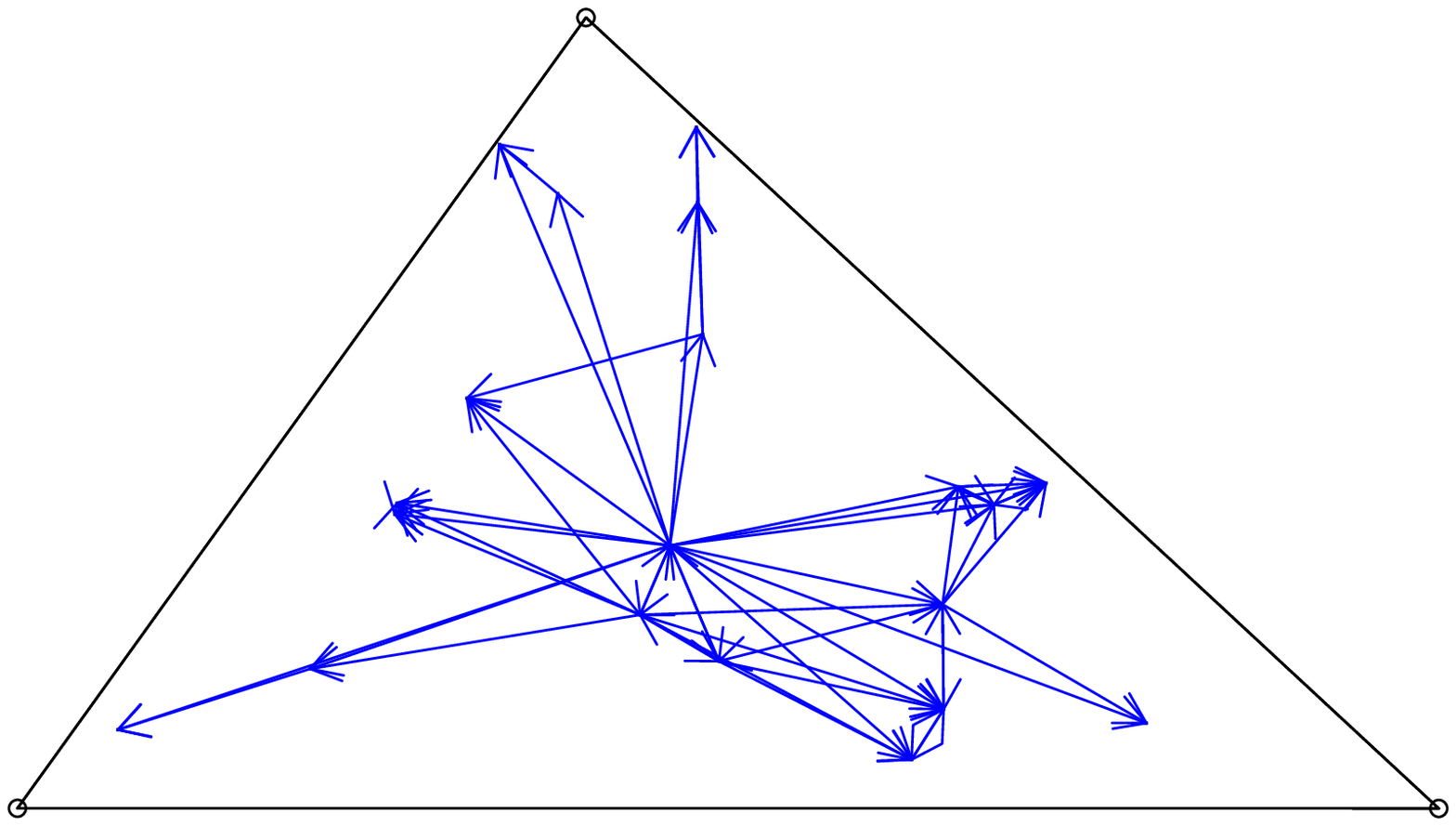} }}
\caption{
\label{fig:CS-arcs-1T}
A realization of 20 $\X$ points generated iid $\U(\TY)$ (left)
and the corresponding arcs for $\NCS^{\tau=1}(x,M_C)$ (right).
}
\end{figure}

Notice that by definition $x \in \NCSt(x,M)$ for all $x \in \TY$.
Furthermore, $\tau \le 1$ implies that $\NCSt(x,M)\subseteq \TY$
for all $x \in \TY$ and $M \in \TY^o$.
For all $ x\in \TY^o \cap R_M(e(x))$,
the edges $e^{\tau}(x)$ and $e(x)$ are coincident iff $\tau=1$.
See Figure \ref{fig:CS-arcs-1T} for
the arcs based on $\NCS^{\tau=1}(x,M_C)$ for 20 $\X$
points in the one triangle case.

Notice that $X_i \stackrel{iid}{\sim} F$,
with the additional assumption
that the non-degenerate two-dimensional
pdf $f$ exists with support $\mS(F)\subseteq \TY$,
implies that the special case in the construction
of $\NCSt(\cdot)$ ---
$X$ falls on the boundary of two edge regions ---
occurs with probability zero.
Note that for such an $F$, $\NCSt(X,M)$ is a triangle for $\tau > 0$ a.s.

Notice that central similarity proximity maps are defined
with $M$-edge regions for $M \in \TY^o$.
Among the four centers considered
in Section \ref{sec:triangle-centers},
 $M_C$ and $M_I$ are inside the triangle, so they can be used in
construction of the central similarity proximity map.

%In general, for central similarity proximity regions with
%$M$-edge regions,
%the similarity ratio of $\NCSt(x,M)$ to $\TY$ is
%$\displaystyle d(x,e^{\tau}(x))/d(M,e(x))$.
%See Figure \ref{fig:CS-IC-M} (right) for $\NCS^{\tau=1}(x,M)$ with $e(x)=e_3$.
%Notice that $\NCSt(\cdot,M)$ satisfies properties \textbf{P1}-\textbf{P8}.

With $M=M_C$, for $x \in R_{CM}(e)$, the similarity ratio of
$N^{\tau}_{CS}(x,M_C)$ to $\TY$ is $\displaystyle d(x,e_{\tau})/d(M_{C},e)$.

\begin{figure}[ht]
\begin{center}
\psfrag{A}{\scriptsize{$\y_1$}}
\psfrag{B}{\scriptsize{$\y_2$}}
\psfrag{C}{\scriptsize{$\y_3$}}
\psfrag{M}{\scriptsize{$M$}}
\psfrag{x}{\scriptsize{$x$}}
\psfrag{D}{}
\psfrag{E}{}
\psfrag{F}{}
\psfrag{d1}{}
\psfrag{d2}{}
\epsfig{figure=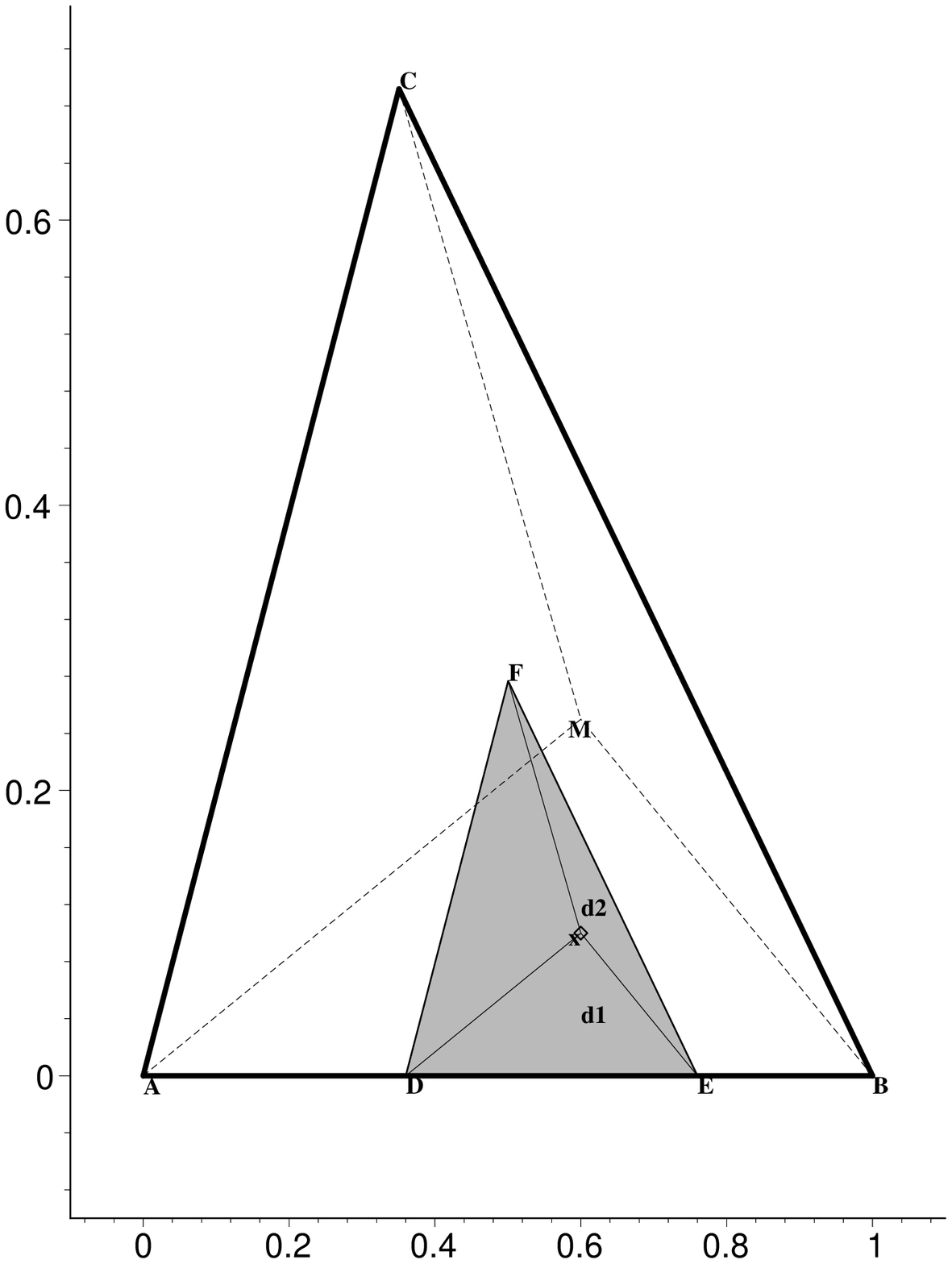, height=140pt , width=200pt}
\psfrag{CM}{\scriptsize{$M_C$}}
\epsfig{figure=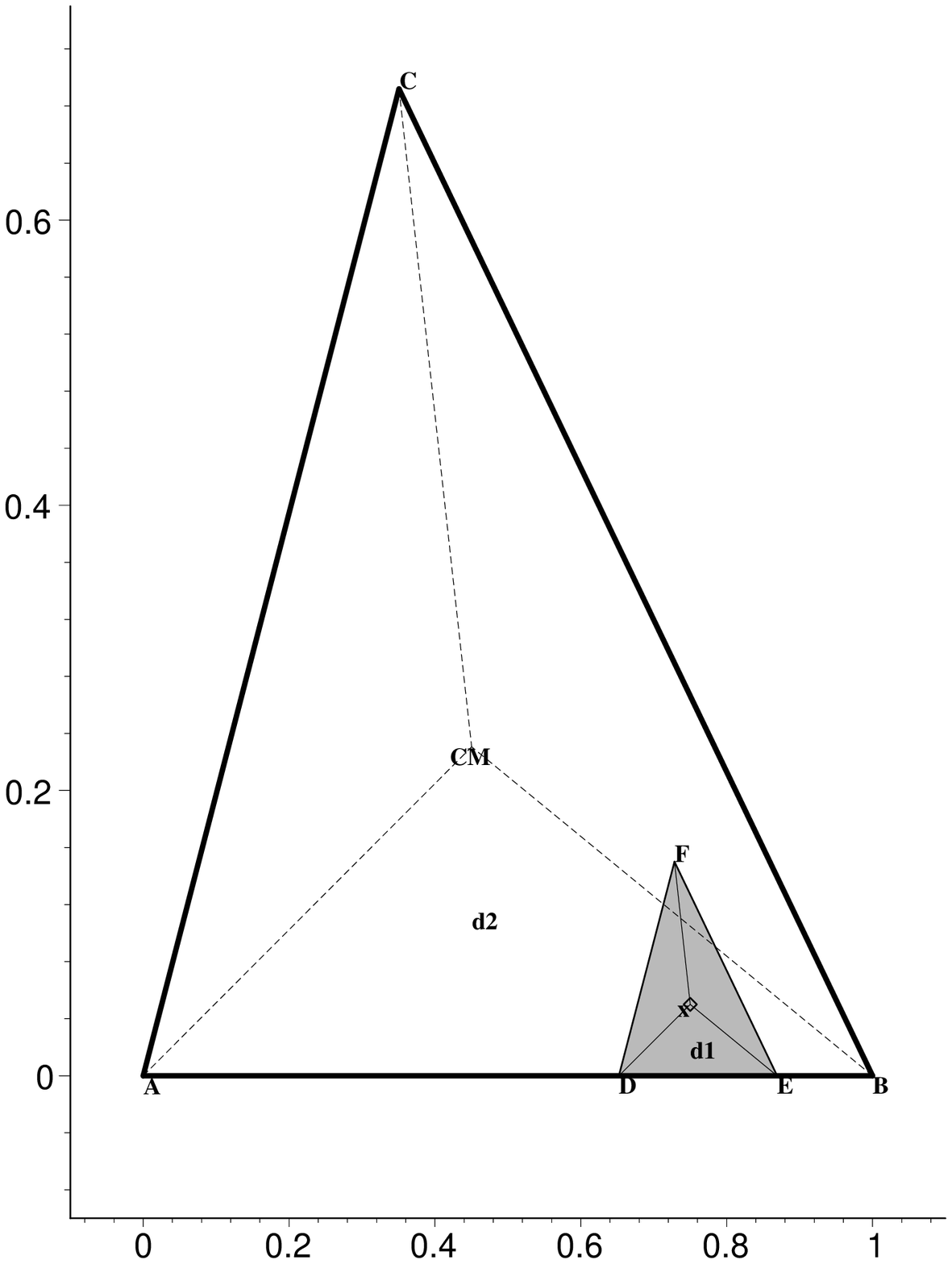, height=140pt , width=200pt}
\psfrag{IC}{\scriptsize{$M_{I}$}}
\epsfig{figure=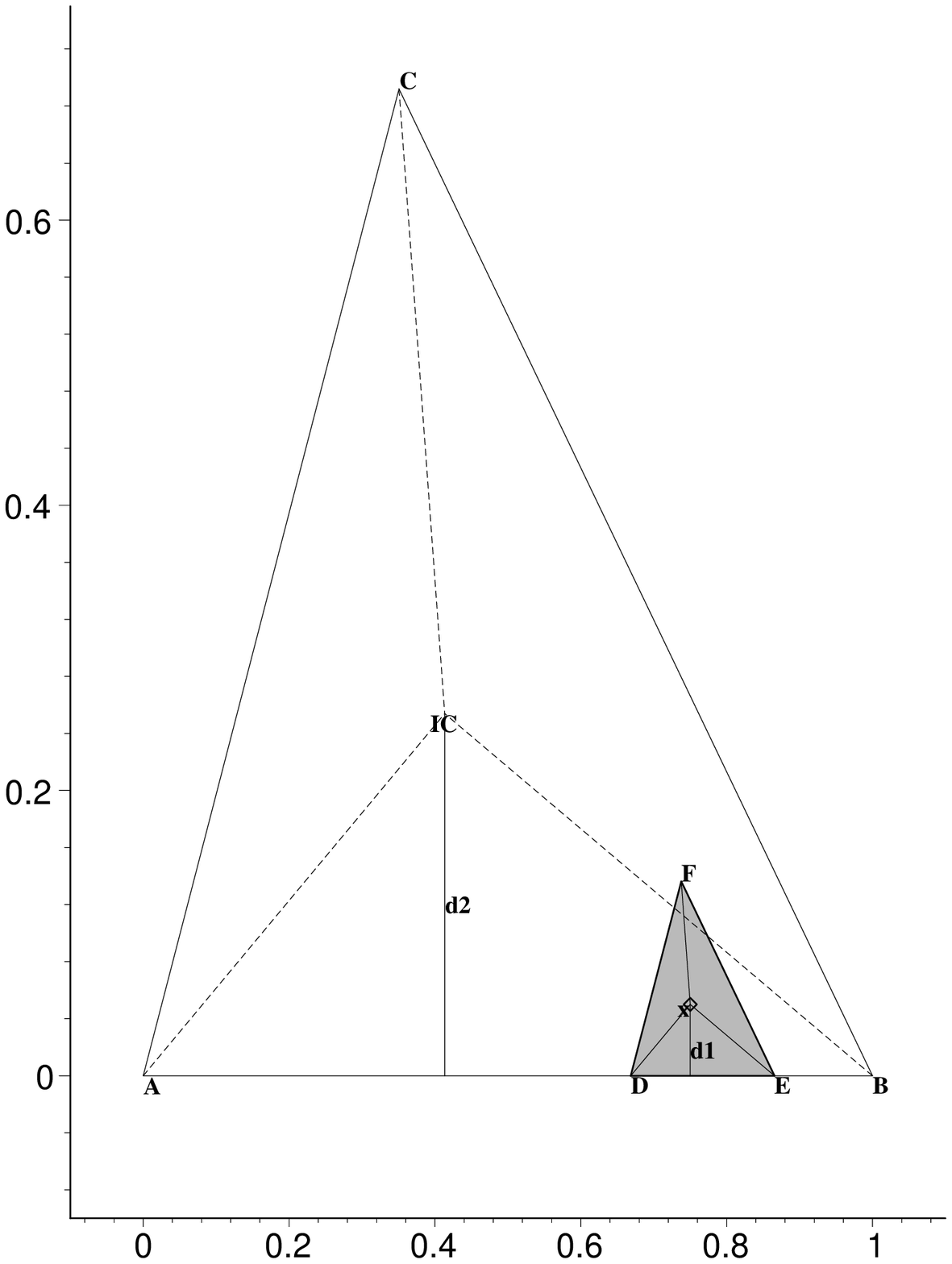, height=140pt , width=200pt}
\end{center}
\caption{$\NCS^{\tau=1}(x,M)$ with an $x \in R_{M}(e_3)$ (top left);
$\NCS^{\tau=1}(x,M_C)$ with an $x \in R_{M_C}(e_3)$ (top right);
and $\NCS^{\tau=1}(x,M_{IC})$ with an $x \in R_{M_{IC}}(e_3)$ (top right).}
\label{fig:CS-IC-M}
\end{figure}

See Figure \ref{fig:ProxMapDefCS} for $\NCS^{\tau=1/2}(x,M_C)$ with $e=e_3$
and Figure \ref{fig:CS-IC-M} for $\NCS^{\tau=1}(x,M_C)$ with $e=e_3$.
The functional form of $\NCSt(x,M_C)$ for an
$x=(x_0,y_0) \in R_{CM}(e)$ is as follows:
\begin{itemize}
\item[]
For $x \in R_{CM}(e_1)$,
\begin{multline*}
\NCSt(x,M_C)=\Biggl \{(x,y) \in T_b:
y \ge y_0+\tau\,(y_0\,(1-c_1)),\;y \le \frac{(1-\tau)\,(y_0\,(1-c_1)-x_0\,c_2)+c_2\,(\tau-x)}{1-c_1};\;\\
 y \le \frac{y_0\,(c_1\,(1+\tau)-\tau)-x_0\,c_2\,(1+\tau)+c_2\,(\tau+x)}{c_1} \Biggr \}.
\end{multline*}
\item[]
For $x \in R_{CM}(e_2)$,
\begin{multline*}
\NCSt(x,M_C)=\Biggl \{(x,y) \in T_b:
y\ge y_0+\tau\,(y_0\,c_1-c_2\,x_0),\; y \le \frac{(1-\tau)\,(y_0\,c_1-x_0\,c_2)+c_2\,x}{c_1};\\
\;y \le \frac{y_0\,(1-c_1\,(1+\tau))+x_0\,c_2\,(1+\tau)}{1-c_1} \Biggr\}.
\end{multline*}
\item[]
For $x \in R_{CM}(e_3)$,
\begin{multline*}
\NCSt(x,M_C)=\Biggl \{(x,y) \in T_b:
y \ge y_0\,(1-\tau);\; y \le \frac{(y_0\,(\tau+c_1)+c_2\,(x-x_0))}{c_1};\;\\
y \le \frac{y_0\,(1-c_1+\tau)-c_2\,(x-x_0)}{1-c_1} \Biggr\}.
\end{multline*}
\end{itemize}

%\clearpage

\subsubsection{$IC$-Central Similarity Proximity Regions}
With $M=M_I$, the similarity ratio is $d(x,e_{\tau})/d(M_I,e)$.
See Figure \ref{fig:CS-IC-M}  for
$N^{\tau=1}_{CS}(x,M_I)$ with $e=e_3$.
The functional form of
$N_{CS}^{\tau}(x,M_I)$ for an $x=(x_0,y_0) \in R_{IC}(e)$ is as
follows:
\begin{itemize}
\item[]
If $x \in R_{IC}(e_1)$, then
{\small
\begin{multline*}
N_{CS}^{\tau}(x,M_I)=\Bigg\{(x,y) \in T_b:
y \ge \frac{y_0\,\left(\sqrt{c_2^2+(1-c_1)^2}+\tau\,(1-c_1)\right)+c_2\,\tau\,(x_0-\tau)}{\sqrt{c_2^2+(1-c_1)^2}};\;\\
y \le \frac{(1-\tau)\,(y_0\,(1-c_1)+x_0\,c_2)+c_2\,(\tau-x)}{1-c_1};\;
y \le \frac{\left(\sqrt{c_2^2+(1-c_1)^2}+\tau\,\sqrt{c_1^2+c_2^2}\right)\,
(y_0\,c_1-x_0\,c_2)-y_0\,\tau\,\sqrt{c_1^2+c_2^2}}{c_1\,\sqrt{c_2^2+(1-c_1)^2}}\Bigg\}.
\end{multline*}
}
\item[]
If $x \in R_{IC}(e_2)$, then
\begin{multline*}
N_{CS}^{\tau}(x,M_I)=\Bigg\{(x,y) \in T_b: y \ge \frac{y_0\,\left(\sqrt{c_1^2+c_2^2}+
\tau\,c_1\right)-\tau\,c_2\,x_0}{\sqrt{c_1^2+c_2^2}};\;
y \le \frac{(1-\tau)\,(c_1\,y_0-c_2\,x_0)+c_2\,x}{c_1};\;\\
 y \le \frac{\left(\sqrt{c_1^2+c_2^2}+\tau\,\sqrt{c_2^2+(1-c_1)^2}\right)\,(x_0\,c_2-y_0\,c_1)+
\sqrt{c_1^2+c_2^2}\,(y_0-c_2\,x)}{(1-c_1)\,\sqrt{c_1^2+c_2^2}}\Bigg\}.
\end{multline*}
\item[]
If $x\in R_{IC}(e_3)$, then
\begin{multline*}
N_{CS}^{\tau}(x,M_I)=\Bigg\{(x,y) \in T_b: y \ge y_0\,(1-\tau);\;
y \le \frac{y_0\,\left(c_1+\tau\,\sqrt{c_1^2+c_2^2}\right)+c_2\,(x-x_0)}{c_1},\;\\
y \le \frac{y_0\,\left(1-c_1+\tau\,\sqrt{c_2^2+(1-c_1)^2}\right)+c_2\,(x_0-x)}{1-c_1}\Bigg\}.
\end{multline*}
\end{itemize}
$N_{CS}^{\tau}(\cdot,M_I)$ also satisfies all the properties
\textbf{P1}-\textbf{P9}.
\begin{remark}
For acute triangles we could use $CC$ or $OC$-edge regions in
central similarity proximity regions which will also
satisfy properties \textbf{P1}-\textbf{P9}.
But for obtuse triangles, \textbf{P2} is not satisfied and edge regions are not
defined in a natural manner.
\end{remark}

In general for $M$-central similarity proximity
regions, the similarity ratio is $d(x,e_{\tau}(x))/d(M,e(x))$.
See Figure \ref{fig:CS-IC-M}
for $N^{\tau=1}_{CS}(x,M)$ with $e=e_3$.
The functional form of
$N_{CS}^{\tau}(x,M)$ for an $x=(x_0,y_0) \in R_{M}(e)$ is as
follows:
\begin{itemize}
\item[]
If $x \in R_{M}(e_1)$, then\\
$\displaystyle N_{CS}^{\tau}(x,M)=\Bigg\{(x,y) \in T_b:
y \ge \frac{y_0\,(m_2\,(1-\tau)\,(c_1-1)+c_2\,(1-m_1))-\tau\,m_2\,c_2\,(1-x_0)}{m_2\,(1-c_1)+c_2\,(1-m_1)};\;\\
y \le \frac{(1-\tau)\,(y_0\,(1-c_1)+x_0\,c_2)+c_2\,(\tau-x)}{1-c_1};\;
y \le \Bigl[y_0\,(c_1\,m_2\,(1-\tau)(1-c_1)+m_1\,c_2\,(\tau+c_1)-c_1\,c_2\,(1+\tau\,m_1)) +
x_0\,(c_2^2\,(1-m_1\,(1-\tau))+c_2\,m_2\,(c_1\,(1-\tau)-c_2))+
\tau\,c_2\,(m_2\,c_1-m_1\,c_2)+(c_2\,m_2\,(1-c_1)-c_2^2\,(1-m_1))\,x\Bigr]\Big/
\Bigl[c_1\,(m_2\,(1-c_1)-c_1\,(1-m_1))\Bigr]\Bigg\}$.
\item[]
If $x \in R_{M}(e_2)$, then\\
$\displaystyle N_{CS}^{\tau}(x,M)=\Bigg\{(x,y) \in T_b:
y \ge \frac{y_0\,(c_2\,m_1-c_1\,m_2)-\tau\,m_2\,(c_2\,x_0-c_1\,y_0)}{m_1\,c_2-c_1\,m_2};\;
y \le \frac{(1-\tau)\,(y_0\,c_1-x_0\,c_2)+c_2\,x}{c_1};\;\\
y \le \Bigl[y_0\,(c_1\,m_2-c_2\,m_1)+\tau\,c_1\,(c_2\,(1-m_1)-m_2\,(1-c_1))+
x_0\,(c_2^2\,(-\tau-m_1\,(1-\tau))+m_2\,c_2\,(\tau\,(1-c_1)+c_1))+
c_2\,(m_1\,c_2-m_2\,c_1)\,x\Bigr]\Big/
\Bigl[(1-c_1)\,(c_1\,m_2-c_2\,m_1)\Bigr] \Bigg\}$.
\item[]
If $x \in R_{M}(e_3)$, then\\
$\displaystyle N_{CS}^{\tau}(x,M)=\Bigg\{(x,y) \in T_b:
y \ge y_0\,(1-\tau);\;y \le \frac{y_0\,(c_1\,m_2\,(1-\tau)+
c_2\,\tau\,m_1)+c_2\,m_2\,(x-x_0)}{c_1\,m_2};\;\\
y\le\frac{y_0\,(m_2\,(1-c_1)\,(1-\tau)+c_2\,\tau)+c_2\,m_2\,(x_0-c_2)}{(1-c_1)\,m_2} \Bigg\}$.
\end{itemize}

Notice that $N_{CS}^{\tau}(\cdot,M)$ also satisfies properties
\textbf{P1}-\textbf{P9}.
$M$-central similarity proximity regions with $m_1 < c_1$ can be defined in a similar fashion.
Furthermore, $\Lambda_0\left(\NCSt(\cdot,M) \right)=\partial(\TY)$
for all $\tau \in (0,1]$ and $\Lambda_0\left(N^{\tau=0}_{CS}(\cdot,M) \right)=\TY$,
since $\lambda\left(\NCSt(x) \right)=0$ iff $x \in e_i$ for $i \in \{1,2,3\}$ or $\tau=0$.

\begin{remark}
Among the family
$\mathscr N^{\tau}_{CS}:=\bigl\{\NCSt(\cdot,M):\;\tau \in [0,1],\;M \in \TY^o \bigr\}$,
every $\NCSt(\cdot,M)$ with $\tau \in (0,1]$ satisfies all the properties in
Section \ref{sec:prox-maps}.
$\square$
\end{remark}

\begin{remark}
\begin{itemize}
\item[]
\item
For $\tau_1 \le \tau_2$, $\NCS^{\tau_1}(x,M) \subseteq \NCS^{\tau_2}(x,M)$ for all $x \in \TY$.
For $\tau_1 < \tau_2$, $\NCS^{\tau_1}(x,M) \subseteq \NCS^{\tau_2}(x,M)$
with equality holding only for $ x \in \partial(\TY)$.
\item
The superset region $\RS\left(\NCSt,M\right)=\emptyset$ for
$\tau \in [0,1)$ and $\RS\left(\NCS^{\tau=1},M\right)=\{M\}$.
\item
For $\tau_1 < \tau_2$, $A\left(\NCS^{\tau_1}(X,M) \right) \le ^{ST} A\left(\NCS^{\tau_2}(X,M) \right)$
for $X$ from a continuous distribution on $\TY$. $\square$
\end{itemize}
\end{remark}

\subsubsection{Extension of $N^{\tau}_{CS}$ to Higher Dimensions}
\label{sec:NCS-higher-D}
The extension of $N^{\tau}_{CS}$ to $\R^d$ for $d > 2$ is straightforward.
the extension for $M=M_C$ is described, the extension for general $M$ is similar.
Let $\Y_{d+1} = \{\y_1,\y_2,\cdots,\y_{d+1}\}$ be $d+1$ points
that do not lie on the same $(d-1)$-dimensional hyperplane.
Denote the simplex formed by these $d+1$ points as $\mathfrak S(\Y_{d+1})$.
%A simplex is the simplest polytope in $\R^d$
%having $d+1$ vertices, $d\,(d+1)/2$ edges and $d+1$ faces of dimension $(d-1)$.
For $\tau \in (0,1]$,
define the central similarity proximity map as follows.
Let $\varphi_i$ be the face opposite vertex $\y_i$
for $i \in \{1,2,\ldots,(d+1)\}$, and ``face regions''
$R_{CM}(\varphi_1),\ldots,R_{CM}(\varphi_{d+1})$
partition $\mathfrak S(\Y_{d+1})$ into $d+1$ regions, namely the $d+1$
polytopes with vertices being the center of mass together
with $d$ vertices chosen from $d+1$ vertices.
For $x \in \mathfrak S(\Y_{d+1}) \setminus \Y_{d+1}$, let $\varphi(x)$ be the
face in whose region $x$ falls; $x \in R(\varphi(x))$.
If $x$ falls on the boundary of two face regions,
$\varphi(x)$ is assigned arbitrarily.
For ${\tau} \in (0,1]$,
the central similarity proximity region
$\NCSt(x,M_C)=\mathfrak S_{\tau}(x)$ is defined to be the simplex
$\mathfrak S_{\tau}(x)$ with the following properties:
\begin{itemize}
\item[(i)] $\mathfrak S_{\tau}(x)$ has faces $\varphi^{\tau}_i(x)$ parallel to
$\varphi_i(x)$ for $i \in \{1,2,\ldots,(d+1)\}$, and for $x \in R_{CM}(\varphi(x))$,
$\tau\, d(x,\varphi(x))=d(\varphi^{\tau}(x),x)$
 where $d(x,\varphi(x))$ is the Euclidean (perpendicular)
distance from $x$ to $\varphi(x)$;
\item[(ii)] $\mathfrak S_{\tau}(x)$ has the same orientation
as and similar to $\mathfrak S(\Y_{d+1})$;
\item[(iii)] $x$ is the center of mass of
$\mathfrak S_{\tau}(x)$, as $M_C$ is of $\mathfrak S (\Y_{d+1})$.
Note that $\tau>1$ implies that $x \in \NCSt(x)$.
\end{itemize}

\subsection{The Behavior of Proximity Regions}
In this section, we provide the conditions for $x$, which,
if satisfied, will imply some sort of increase in the size
of the proximity regions we have defined.
Let $N(\cdot)$ be any proximity map defined on the
measurable space $\Omega$ with measure $\mu$,
and let $\{x_n\}_{n=1}^{\infty}$ be a sequence of points in $\Omega$.
We say $N(x_n)$ \emph{gets larger} if $N(x_n)\subseteq N(x_m)$ for $m \ge n$,
and $N(x_n)$ \emph{gets strictly larger} if $N(x_n)\subsetneq N(x_m)$ for $m>n$.

In the following theorems we will assume $\Omega=\R^2$ with
$\mu$ being the $\R^2$-Lebesgue measure $\lambda$ and
$M$-vertex regions are defined with points $M \in \R^2 \setminus \Y_3$.

\begin{theorem}
For arc-slice proximity regions with $M$-vertex regions for an
$M \in \R^2 \setminus \Y_3$, as $d(x,\y)$ (strictly) increases
for $x$ lying on a ray from $\y$ in $R_M(\y) \setminus \RS(\NAS,M)$,
$\NAS(x,M)$ gets (strictly) larger.
\end{theorem}
\textbf{Proof:}
For $x,y$ lying on a ray from $\y$
in $R_M(\y) \setminus \RS(\NAS,M)$,
if $d(x,\y) \le d(y,\y)$, then $B(x,r(x)) \subseteq B(y,r(y))$,
which implies $\NAS(x,M) \subseteq \NAS(y,M)$,
hence $\NAS(x,M)$ gets larger as $d(x,\y)$ increases
for $x$ lying on a ray from $\y$ in $R_M(\y) \setminus \RS(\NAS,M)$.
The strict version follows similarly.
If $x,\, y \in R_M(\y)\cap \RS(\NAS,M)$,
then $\NAS(x,M)=\NAS(y,M)=\TY$. $\blacksquare$

Let $\ell(\y,x)$ be the line at $x$ parallel to $e(x)$ for $x \in R_M(\y)$
where $e(x)$ is the edge opposite vertex $\y$.

\begin{theorem}
For the proportional-edge proximity maps with
$M$-vertex regions for an $M \in \R^2 \setminus \Y_3$,
as $d(\ell(\y,x),\y)$ (strictly) increases for
$x \in R_M(\y) \setminus \RS\left(\NPE^r,M\right)$, $\NPE^r(x,M)$ gets (strictly) larger for $r < \infty$.
\end{theorem}
\textbf{Proof:} For $x,\, y \in R_M(\y) \setminus \RS\left(\NPE^r,M\right)$,
if $d(\ell(\y,x),\y) \le d(\ell(\y,y),\y)$,
then by definition $\NPE^r(x,M) \subseteq \NPE^r(y,M)$,
hence the result follows.
The strict version follows similarly.
If $x,\, y \in R_M(\y) \cap \RS\left(\NPE^r,M\right)$,
then $\NPE^r(x,M)=\NPE^r(y,M)=\TY$,
and if $r=\infty$ and $x,y \in \TY \setminus \Y_3$,
$\NPE^r(x,M)=\NPE^r(y,M)=\TY$.$\blacksquare$

Note that as $d(\ell(\y,x),\y)$ increases for $x\in R_M(\y)$,
$d(\ell(\y,x),M)$ decreases, provided that $M \in \TY^o$
and $M$-vertex regions are convex.

We define the $M$-edge regions, $R_M(e)$, with points $M \in \TY^o$.
\begin{theorem}
\label{thm:CS-depart}
For central similarity proximity regions with
$M$-edge regions for an $M \in \TY^o$,
as $d(x,e)$ (strictly) increases for $x \in R_M(e)$,
the area $A\left(\NCSt(x,M) \right)$ (strictly) increases for $\tau \in (0,1]$.
\end{theorem}
\textbf{Proof:}
For $x,\, y \in R_M(e)$ and $\tau \in (0,1]$,
if $d(x,e) \le d(y,e)$ then the similarity ratio of
$N^{\tau}_{CS}(y,M)$ to $\TY$ is larger than or equal to that of
$N^{\tau}_{CS}(x,M)$, which in turn implies that
$A\left(\NCSt(x,M) \right) \le A\left(\NCSt(y,M) \right)$.
The strict version follows similarly. $\blacksquare$

Observe that the statement of Theorem \ref{thm:CS-depart}
is about the area $A\left(\NCSt(x,M) \right)$.
We need further conditions for $\NCSt(x,M)$ to get larger.

\begin{figure}[ht]
\centering
\scalebox{.3}{\input{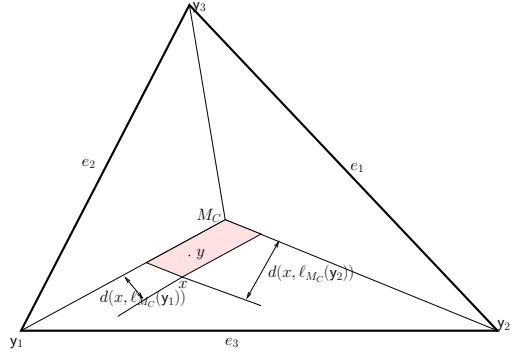}}
\caption{The figure for $x,\,y \in R_M(e_3)$ described in Theorem \ref{thm:quadrilateral-depart}.}
\label{fig:quadrilateral-departure}
\end{figure}

\begin{theorem}
\label{thm:quadrilateral-depart}
Let $\ell_M(\y)$ be the line joining $M$ and vertex $\y \in \Y_3$.
As $d(x,\ell_M(\y_j))$ and $d(x,\ell_M(\y_k))$ both (strictly)
decrease for $x \in R_M(e_l)$ where $j,k,l$ are distinct,
$\NCSt(x,M)$ (strictly) increases for $\tau \in (0,1]$.
\end{theorem}
\textbf{Proof:}
Suppose, without loss of generality, that $x,\, y \in R_M(e_3)$.
Consider the set
$$S(e_3,x):=\left\{y \in R_M(e_3):d(y,\ell_M(\y_1)) \le d(x,\ell_M(\y_1)) \text{ and }
d(y,\ell_M(\y_2)) \le d(x,\ell_M(\y_2))\right\},$$
which a parallelogram.
See Figure \ref{fig:quadrilateral-departure} for an example
of $S(e_3,x)$ with $M=M_C$ and $e=e_3$.
Given $x$, for $y\in S(e_3,x)$, by construction,
$\NCSt(x,M) \subseteq \NCSt(y,M)$.
Then the desired result follows for $\tau \in (0,1]$.
Observe that if $x_{n+1}$ is in $S(e_3,x_n)$, then
$d(x_n,\ell_M(\y_1))$ and $d(x_n,\ell_M(\y_2))$ both decrease.
The strict version follows similarly. $\blacksquare$

\begin{remark}
For $\RS(\NY)$ with positive area,  by definition, as $x \rightarrow y \in \RS(\NY)$,
$\NY(x) \rightarrow \TY$ and hence $\argsup_{x \in \TY}A(\NY(x)) \in \RS(\NY)$
with $\sup_{x \in \TY}A(\NY(x))=A(\TY).$
\begin{itemize}
\item As $x \rightarrow M_{CC}$ in a non-obtuse triangle $\TY$,
$\NAS(x,M_{CC}) \rightarrow \TY$ and
$$\argsup_{x \in \TY}A\left(\NAS(x,M_{CC}) \right)=M_{CC}
\text{ with }
\sup_{x \in \TY}A\left(\NAS(x,M_{CC}) \right)=A(\TY).$$

\item For $r>3/2$, as $x \rightarrow y \in (\RS\left( \NPE^r,M \right))^o$,
$\NPE^r(x,M) \rightarrow \TY$ hence
$$\argsup_{x \in \TY} A(\NPE^r(x,M)) \in \RS\left( \NPE^r,M \right)
\text{ with }
\sup_{x \in \TY}A(\NPE^r(x,M))=A(\TY).$$

\item For $r<3/2$, if $M \notin \Tr$, then as $x \rightarrow M$,
$\NPE^r(x,M) \rightarrow \TY$ and
$$\argsup_{ x \in \TY} A(\NPE^r(x,M))=M
\text{ with }
\sup_{x \in \TY}A(\NPE^r(x,M))=A(\TY).$$
If $M \in (\Tr)^o$, then as $x \rightarrow M$,
$\NPE^r(x,M) \rightarrow \NPE^r(M,M)\subsetneq \TY$,
but still
$$\argsup_{ x \in \TY} A(\NPE^r(x,M))=M
\text{ with }
\sup_{x \in \TY}A(\NPE^r(x,M))=A(\NPE^r(M,M)).$$
If $M \in \partial(\Tr)$, then as $x \rightarrow M$,
$\NPE^r(x,M) \rightarrow \NPE^r(M,M)\subseteq \TY$,
and
$$\argsup_{ x \in \TY} A(\NPE^r(x,M))=M
\text{ with }
\sup_{x \in \TY}A(\NPE^r(x,M))=A(\NPE^r(M,M))$$
which might be $\TY$ or a proper subset of $\TY$.

\item As $x \rightarrow M_C$, $\NPE^{3/2}(x,M_C) \rightarrow \TY$ and
$$\argsup_{ x \in \TY}A\left(\NPE^{3/2}(x,M_C) \right)=M_C
\text{ with }
\sup_{x \in \TY}A\left(\NPE^{3/2}(x,M_C) \right)=A(\TY).$$

\item As $x \rightarrow M$, $\NCS^{\tau=1}(x,M) \rightarrow \TY$ and
$$\argsup_{ x \in \TY}A(\NCS^{\tau=1}(x,M))=M
\text{ with }
\sup_{x \in \TY}A(\NCS^{\tau=1}(x,M))=A(\TY).$$
For $\tau<1$, as $x \rightarrow M$, $\NCSt(x,M) \rightarrow \NCSt(M,M)$ and
$$\argsup_{x \in \TY}A\left(\NCSt(x,M) \right)=M
\text{ with }
\sup_{x \in \TY}A\left(\NCSt(x,M) \right)=A(\NCSt(M,M)).~~\square $$
\end{itemize}
\end{remark}

Although the comments in the above remark follow by definition,
they will be indicative of whether the asymptotic distribution of
the domination number of the associated PCD is degenerate or not.

\section{Relative Arc Density and Domination Number of PCDs}
\label{sec:rel-dens-dom-num}
\subsection{Relative Arc Density}
The \emph{relative arc density} of a digraph $D=(\V,\A)$
of order $|\V| = n$,
denoted as $\rho(D)$,
is defined as
$$
\rho(D) = \frac{|\A|}{n(n-1)}
$$
where $|\cdot|$ denotes the set cardinality functional (\cite{janson:2000}).
Thus $\rho(D)$ represents the ratio of the number of arcs
in the digraph $D$ to the number of arcs in the complete symmetric
digraph of order $n$, which is $n(n-1)$.
For brevity of notation we use \emph{relative density}
rather than relative arc density henceforth.

If $X_1,\cdots,X_n \stackrel{iid}{\sim} F$
the relative density
of the associated data-random PCD $D$,
denoted as $\rho(\X_n;h,\NY)$, is a $U$-statistic,
\begin{eqnarray}
\rho(\X_n;h,\NY) =
  \frac{1}{n(n-1)}
    \sum\hspace*{-0.1 in}\sum_{i < j \hspace*{0.25 in}}   \hspace*{-0.1 in}
      \hspace*{-0.1 in}\;\;h(X_i,X_j;\NY)
\end{eqnarray}
where
\begin{eqnarray}
h(X_i,X_j;\NY)&=& \I\{(X_i,X_j) \in \A\}+ \I\{(X_j,X_i) \in \A\} \nonumber \\
       &=& \I\{X_j \in \NY(X_i)\}+ \I\{X_i \in \NY(X_j)\},
\end{eqnarray}
where $\I(\cdot)$ is the indicator function.
We denote $h(X_i,X_j;\NY)$ as $h_{ij}$ for brevity of notation.
Since the digraph is asymmetric, $h_{ij}$ is defined as
the number of arcs in $D$ between vertices $X_i$ and $X_j$,
in order to produce a symmetric kernel with finite variance (\cite{lehmann:1988}).

The random variable $\rho_n := \rho(\X_n;h,\NY)$ depends on $n$ and $\NY$ explicitly
and on $F$ implicitly.
The expectation $\E[\rho_n]$, however, is independent of $n$
and depends on only $F$ and $\NY$:
\begin{eqnarray}
0 \leq \E[\rho_n] = \frac{1}{2}\E[h_{12}] \leq 1 \text{ for all $n\ge 2$}.
\end{eqnarray}
The variance $\Var[\rho_n]$ simplifies to
\begin{eqnarray}
\label{eqn:var-rho}
0 \leq
  \Var[\rho_n] =
     \frac{1}{2n(n-1)} \Var[h_{12}] +
     \frac{n-2}{n(n-1)} \Cov[h_{12},h_{13}]
  \leq 1/4.
\end{eqnarray}
A central limit theorem for $U$-statistics
(\cite{lehmann:1988})
yields
\begin{eqnarray}
\sqrt{n}\bigl(\rho_n-\E[\rho_n]\bigr) \stackrel{\mathcal{L}}{\longrightarrow} \N\bigl(0,\Cov[h_{12},h_{13}]\bigr)
\end{eqnarray}
provided $\Cov[h_{12},h_{13}] > 0$.
The asymptotic variance of $\rho_n$, $\Cov[h_{12},h_{13}]$,
depends on only $F$ and $\NY$.
Thus, we need determine only
$\E[h_{12}]$
and
$\Cov[h_{12},h_{13}]$
in order to obtain the normal approximation
\begin{eqnarray}
\rho_n \stackrel{\text{approx}}{\sim}
\N\bigl(\E[\rho_n],\Var[\rho_n]\bigr) =
\N\left(\frac{\E[h_{12}]}{2},\frac{\Cov[h_{12},h_{13}]}{n}\right) \text{ for large $n$}.
\end{eqnarray}

\subsection{Domination Number}
%\begin{definition}
In a digraph $D=(\V,\A)$, a vertex $v \in \V$ \emph{dominates}
itself and all vertices of the form $\{u: vu \in \A\}$.
A \emph{dominating set} $S_D$ for the digraph $D$ is a subset of
$\V$ such that each vertex $v \in \V$ is dominated by a vertex in
$S_D$. A \emph{minimum dominating set} $S^*_{D}$ is a dominating
set of minimum cardinality and the \emph{domination number}
$\g(D)$ is defined as $\g(D):=|S^*_{D}|$ (see, \cite{lee:1998}) where
$|\cdot|$ denotes the set cardinality functional.
See \cite{chartrand:1996} and \cite{west:2001} for more on graphs and digraphs.
If a minimum dominating set is of size one,
we call it a \emph{dominating point}.
%\end{definition}

Note that for $|\V|=n>0$, $1 \le \g(D) \le n$, since $\V$ itself is always
a dominating set.

\subsection{Asymptotic Distribution of Relative Arc Density of PCDs}
\label{asy-dist-relative-density}
By detailed geometric probability calculations,
provided in \cite{ceyhan:arc-density-PE} and \cite{ceyhan:arc-density-CS}
the mean and the asymptotic variance of the relative density of the
proportional-edge and central similarity PCDs can explicitly be computed.
The central limit theorem for $U$-statistics then establishes
the asymptotic normality under the uniform null hypothesis.
These results are summarized in the following theorems.

\begin{theorem}
For $r \in [1,\infty)$,
\begin{eqnarray}
 \frac{\sqrt{n}\,\bigl(\rho_n(r)-\mu(r)\bigr)}{\sqrt{\nu(r)}}
 \stackrel{\mathcal{L}}{\longrightarrow}
 \N(0,1)
\end{eqnarray}
where
\begin{eqnarray}
\label{eqn:Asymean}
\mu(r) =
 \begin{cases}
  \frac{37}{216}r^2                                 &\text{for} \quad r \in [1,3/2), \\
  -\frac{1}{8}r^2 + 4 - 8r^{-1} + \frac{9}{2}r^{-2}  &\text{for} \quad r \in [3/2,2), \\
  1 - \frac{3}{2}r^{-2}                             &\text{for} \quad r \in [2,\infty), \\
 \end{cases}
\end{eqnarray}
and
\begin{equation}
\label{eqn:Asyvar}
\nu(r) =\nu_1(r) \,\I(r \in [1,4/3)) + \nu_2(r) \,\I(r \in [4/3,3/2))+
\nu_3(r) \,\I(r \in [3/2,2)) + \nu_4(r) \,\I( r \in [2,\infty])
\end{equation}
with
{\small
\begin{align*}
  \nu_1(r) &=\frac{3007\,r^{10}-13824\,r^9+898\,r^8+77760\,r^7-117953\,r^6+48888\,r^5
-24246\,r^4+60480\,r^3-38880\,r^2+3888}{58320\,r^4},\\
  \nu_2(r) &=\frac{5467\,r^{10}-37800\,r^9+61912\,r^8+46588\,r^6-191520\,r^5+13608\,r^4
+241920\,r^3-155520\,r^2+15552}{233280\,r^4},  \\
  \nu_3(r) &=-[7\,r^{12}-72\,r^{11}+312\,r^{10}-5332\,r^8+15072\,r^7+13704\,r^6-
139264\,r^5+273600\,r^4-242176\,r^3\\
& +103232\,r^2-27648\,r+8640]/[960\,r^6],\\
  \nu_4(r) &=\frac{15\,r^4-11\,r^2-48\,r+25}{15\,r^6}.
\end{align*}
}
For $r=\infty$, $\rho_n(r)$ is degenerate.
\end{theorem}
%See Appendix 1 for the proof.

%Similarly, by detailed geometric  probability calculations provided in (\cite{ceyhan:arc-density-CS}),
%the mean and the asymptotic variance of the relative
%density of the PCD can be
%calculated explicitly.
%The central limit theorem for $U$-statistics
%then establishes the asymptotic normality under the uniform null hypothesis.
%These results are summarized in the following theorem.

\begin{theorem}
For $\tau \in (0,1]$, the relative density of the
central similarity proximity digraph converges in law
to the normal distribution;
i.e., as $n \rightarrow \infty$,
\begin{eqnarray}
 \frac{\sqrt{n}(\rho_n(\tau)-\mu(\tau))}{\sqrt{\nu(\tau)} }
 \stackrel{\mathcal{L}}{\longrightarrow}
 \N(0,1)
\end{eqnarray}
where
\begin{eqnarray}
\label{eqn:CSAsymean} \mu(\tau) = \tau^2/6
\end{eqnarray}
and
\begin{eqnarray}
\label{eqn:CSAsyvar} \nu(\tau) =\frac
{\tau^4(6\,\tau^5-3\,\tau^4-25\,\tau^3+\tau^2+49\,\tau+14)}{45\,(\tau+1)(2\,\tau+1)(\tau+2)}.
\end{eqnarray}
For $\tau=0$, $\rho_n(\tau)$ is degenerate for all $n>1$.
\end{theorem}
%See the Appendix for the derivation.

\subsection{Asymptotic Distribution of Domination Number of PCDs}
\label{asy-dist-domination-number}
Recall the triangle $\Tr$ defined in Equation \eqref{eqn:T^r-def}
(see also Figure \ref{fig:Tr-RS-NDA-CC} for $\Tr$ with $r=\sqrt{2}$).
Let $\g_n(r,M):=\g\left(\X_n, \NPE^r,M \right)$ be the domination number
of the PCD based on $\NPE^r$ with $\X_n$, a set of iid random
variables from $\UT$, with $M$-vertex regions.

\begin{figure}[ht]
\begin{center}
\psfrag{r}{\scriptsize{$r$}}
\psfrag{pr}{\scriptsize{$p_r$}}
\epsfig{figure=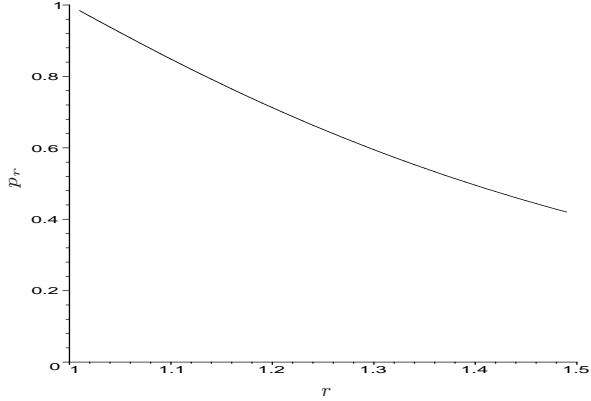, height=150pt, width=220pt}
\end{center}
\caption{ Plotted is the probability $p_r=\lim_{n\rightarrow
\infty}P\left( \g_n(r,M)=2 \right)$ given in Equation \eqref{eqn:p_r-form} as a function of $r$ for $r \in [1,3/2)$ and $M
\in \{t_1(r),t_2(r),t_3(r)\}$.} \label{fig:pglt2ofr}
\end{figure}

The domination number $\g_n(r,M)$ of the PCD has the following
asymptotic distribution (\cite{ceyhan:dom-num-NPE-MASA}).
As $n \rightarrow \infty$,
\begin{equation}
\label{eqn:asymptotic-NYr}
\g_n(r,M) \stackrel{\mathcal L}{\longrightarrow}
\left\lbrace \begin{array}{ll}
       2+\BER(1-p_r)& \text{for $r \in [1,3/2)$ and $M \in \{t_1(r),t_2(r),t_3(r)\}$,}\\
       1            & \text{for $r>3/2$ and $M \in \TY^o$,}\\
       3            & \text{for $r \in [1,3/2)$ and $M \in \Tr\setminus \{t_1(r),t_2(r),t_3(r)\}$,}\\
\end{array} \right.
\end{equation}
where $\stackrel{\mathcal L}{\longrightarrow}$ stands for ``convergence in law"
and $\BER(p)$ stands for Bernoulli distribution with probability of
success $p$, $\Tr$ and $t_i(r)$ are defined in Equation \eqref{eqn:T^r-def}, and for $r \in [1,3/2)$ and $M \in
\{t_1(r),t_2(r),t_3(r)\}$,
\begin{equation}
\label{eqn:p_r-form}
p_r=\int_0^{\infty}\int_0^{\infty}\frac {64\,r^2}{9\,(r-1)^2}\,w_1\,w_3\,\exp\left(\frac{4\,r}
{3\,(r-1)}\,(w_1^2+w_3^2+2\,r\,(r-1)\,w_1\,w_3)\right)\,dw_3w_1,
\end{equation}
%is given in Theorem \ref{thm:g_nNYr=2-for-r-vertex-of-Tr}
and for $r=3/2$ and $M=M_C=\left\{(1/2,\sqrt{3}/6)\right\}$, $p_r \approx 0.7413$,
which is not computed as in Equation \eqref{eqn:p_r-form};
for its computation, see \cite{ceyhan:dom-num-NPE-SPL}.
For example, for $r=5/4$ and
$M \in \left\{t_1(r)=\left(3/10,\sqrt{3}/10\right),t_2(r)=\\
\left(7/10,\sqrt{3}/10\right),t_3(r)=\left(1/2,3\sqrt{3}/5\right)\right\}$,
$p_r \approx 0.6514$.
See Figure \ref{fig:pglt2ofr} for the plot of the numerically computed values
(i.e., the values computed by numerical integration of the expression in
Equation \eqref{eqn:p_r-form}) of $p_r$ as a function of $r$.
Notice that in the nondegenerate case in \eqref{eqn:asymptotic-NYr},
$\E[\g_n(r,M)]=3-p_r$ and $\Var[\g_n(r,M)]=p_r(1-p_r)$.

The results in Theorem \ref{thm:geo-inv-arc-prob} and
Corollaries \ref{cor:geo-inv-arc-prob-1} and
\ref{cor:geo-inv-arc-prob-2} also hold for relative arc density
and the domination number of PCDs based $\NY$.
That is, we have the following corollary.
\begin{corollary}
Given any triangle $T_o$ and $\X_n$ a random sample from $\U(T_o)$.
Suppose the PCD, $D_o$ is defined in such a way that
the ratio of the area of $N(x)$ to the area of the triangle $T_o$
is preserved under the uniformity preserving transformation,
then the distributions of the relative arc density and
the domination number of $D_o$ are geometry invariant.
\end{corollary}
%{\bf Proof:}
%A composition of translation, rotation, reflections, and scaling
%will take any given triangle $T_o = T\left(\y_1,\y_2,\y_3\right)$
%to the ``basic'' triangle $T_b = T\left((0,0),(1,0),(c_1,c_2)\right)$
%with $0 < c_1 \le 1/2$, $c_2 > 0$ and $(1-c_1)^2+c_2^2 \le 1$,
%preserving uniformity.
%The transformation $\phi_e: \R^2 \rightarrow \R^2$ given by
%$\displaystyle \phi_e(u,v) = \left( u+\frac{1-2\,c_1}{\sqrt{3}}\,v,\frac{\sqrt{3}}{2\,c_2}\,v \right)$
%takes $T_b$ to the equilateral triangle
%$T_e = T\left((0,0),(1,0),\left(1/2,\sqrt{3}/2\right)\right)$.
%Investigation of the Jacobian shows that $\phi_e$
%also preserves uniformity.
%Furthermore, the composition of $\phi_e$ with the rigid motion transformations
%maps the boundary of the original triangle, $T_o$,
%to the boundary of the equilateral triangle, $T_e$,
%the straight lines in $T_o$
%to the straight lines in $T_e$,
%and lines parallel to the edges of $T_o$
%to lines parallel to the edges of $T_e$.
%Since the distribution of $\g(\X_n;\U(\TY),\NY^r)$
%involves only probability content of unions and intersections
%of regions bounded by precisely such lines,
%and the probability content of such regions is preserved since uniformity is preserved,
%the desired result follows.
%$\blacksquare$

\section{Two New Proximity Maps}
\label{sec:new-prox-map}
In this section, we introduce two new proximity maps
and investigate their properties.

\subsection{Directional-Doubling Proximity Maps}
\label{sec:N-DD}
\begin{figure}[ht]
\begin{center}
\psfrag{A}{\scriptsize{$\y_1$}}
\psfrag{B}{\scriptsize{$\y_2$}}
\psfrag{C}{\scriptsize{$\y_3$}}
\psfrag{IC}{\scriptsize{$M_{C}$}}
\psfrag{x}{\scriptsize{$z$}}
\psfrag{M}{\scriptsize{$M$}}
\psfrag{D}{}
\psfrag{E}{}
\epsfig{figure=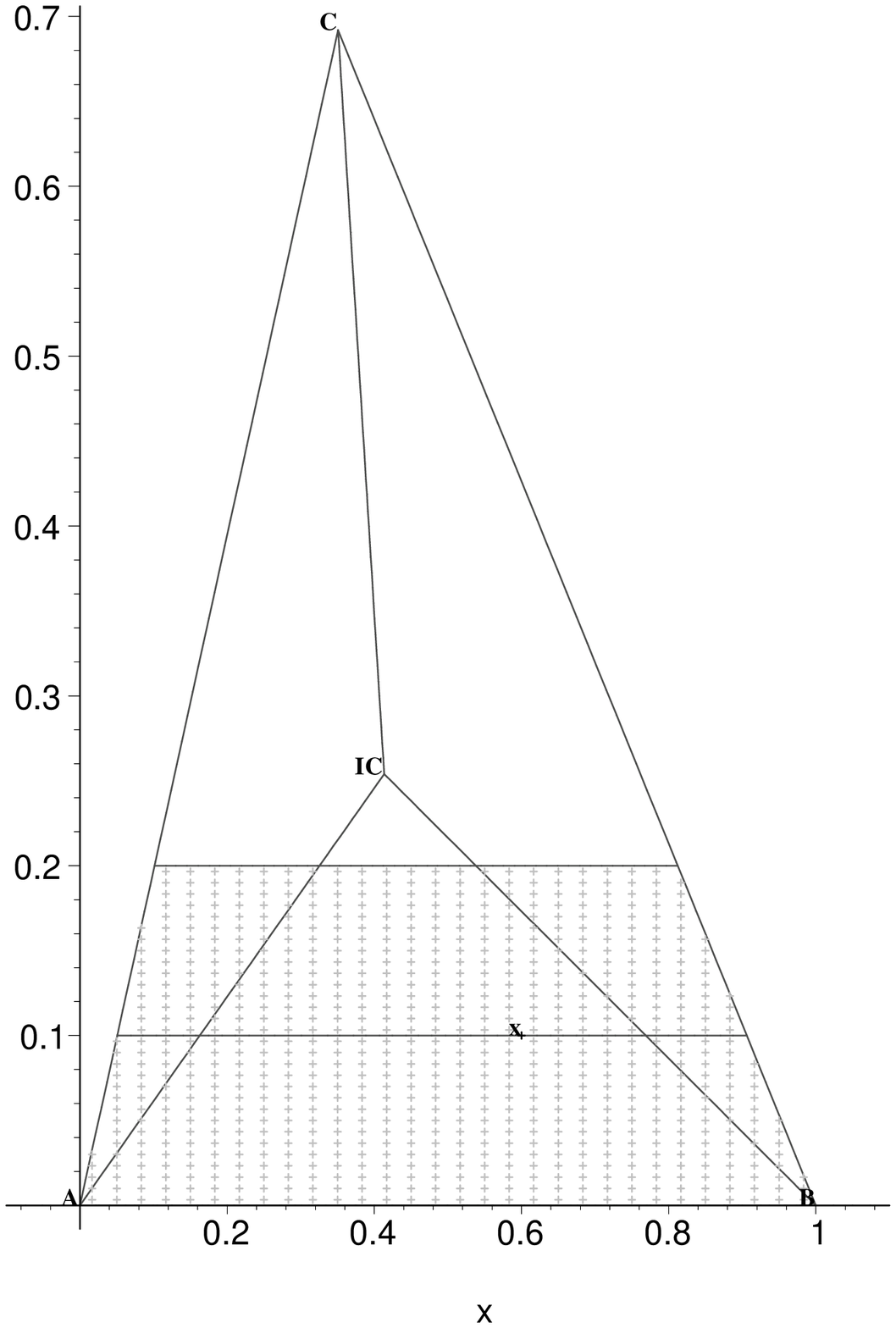, height=140pt , width=200pt}
\epsfig{figure=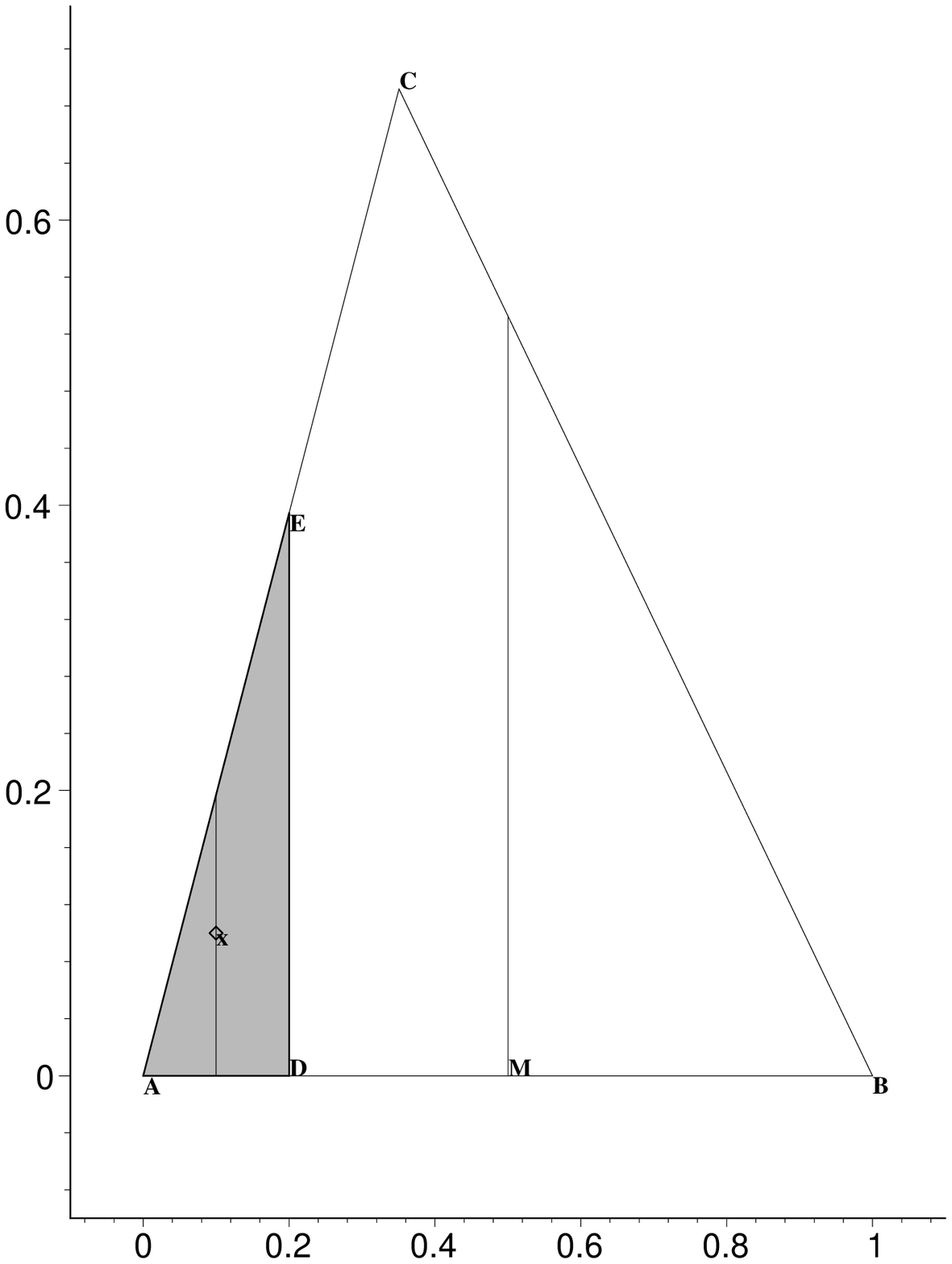, height=140pt , width=200pt}
\end{center}
\caption{
An example of directional-doubling proximity region (left)
and double-$X$ proximity region (right).}
\label{fig:new-prox-regions}
\end{figure}

Without loss of generality,
we can assume that $\TY=T_b$.
Partition the triangle $T_b$ by $M$-edge regions
to obtain $R_M(e_i)$ for $i=1,2,3$.
For $z \in R_M(e_i)$,
directional-doubling proximity map is defined as
$$N_{DD}(z,M):=\left\{\mathbf{x} \in T_b: d(x,e_i) \le 2\,d(z,e_i)\right\}.$$
See Figure \ref{fig:new-prox-regions} (left) for $M=M_C$.
If $z \in e_i$ then $N_{DD}(z,M):=e_i.$
Notice that if $z \not\in e_i$,
then $N_{DD}(z,M)$ is a quadrilateral.
Among the properties,
\textbf{P1} and \textbf{P2} follows trivially.
The line at $z \in R_M(e_i)$ parallel to $e_i$ divides the region into two
pieces (half-way in the perpendicular direction to $e_i$)
so \textbf{P3} holds in this special sense.
\textbf{P4} and \textbf{P5} both fail, since $N_{DD}(z,M)$ is a quadrilateral.
\textbf{P6} holds if $M \in T(M_1,M_2,M_3)$,
otherwise it fails since $\RS(N_{DD},M)$ will have positive area.
\textbf{P7} also follows by definition.
However, \textbf{P8} holds only when $M=M_C$.

Property \textbf{P9} follows for $N_{DD}$,
since $N_{DD}(z,M)$ is constructed with the boundary of $\TY$
and parallel lines to the edges,
by Corollary \ref{cor:geo-inv-arc-prob-1},
geometry invariance for uniform data follows.
That is, the distributions of relative arc density and
the domination number of the corresponding PCD do not depend on
the geometry of the triangle $\TY$.
Hence, it suffices to compute them for the standard equilateral triangle only.
Furthermore,
$\Lambda_0(N_{DD})=\partial(\TY)$ since $N_{DD}(x,M)$ has zero area
iff $x \in \partial(\TY)$.

\subsection{Double-$X$ Proximity Maps}
\label{sec:N-DX}
Without loss of generality,
we can assume that $\TY=T_b$.
Partition the triangle $T_b$
using the vertical line at $M_3=(1/2,0)$ as in Figure \ref{fig:new-prox-regions} (right).
Let $R_1:=\{(x,y) \in T_b: x < 1/2\}$ and
$R_2:=\{(x,y) \in T_b: x > 1/2\}$.
If $(x,y) \in T_b$ with $x=1/2$, assign $(x,y)$
arbitrarily to one of $R_1$ or $R_2$.
We define the double-$X$ proximity map as follows.
For $z=(x_o,y_o) \in T_b \setminus \{\y_1,\y_2\}$
\begin{equation*}
\label{eqn:}
N_{DX}(z):=
\begin{cases}
\{(x,y) \in T_b: x \le 2\,x_o\}        &\text{if} \quad z \in R_1, \\
\{(x,y) \in T_b: 1-x \le 2\,(1-x_o)\}  &\text{if} \quad z \in R_2.
\end{cases}
\end{equation*}
If $z=(x_o,y_o) \in \{\y_1,\y_2\}$, then $N_{DX}(z):=\{z\}$.
See also Figure \ref{fig:new-prox-regions} (right).
Notice that if $z \not \in \{\y_1,\y_2\}$,
then $N_{DD}(z,M)$ is a right triangle or a quadrilateral.
Among the properties,
\textbf{P1} and \textbf{P2} follows trivially.
The vertical line at $z$ divides the region into two
pieces (half-way along the $x$-coordinate)
so \textbf{P3} holds in this special sense.
\textbf{P4} and \textbf{P5} fails to hold,
since $N_{DD}(z,M)$ may be a quadrilateral for some $z \in T_b$.
\textbf{P6} holds if $M=(1/2,y)$, otherwise $\RS(N_{DX},M)$ has positive area.
\textbf{P7} also follows by definition.
However, \textbf{P8} holds only when the regions
$R_1$ and $R_2$ are constructed at a point where
the vertical line divides the area into two equal pieces.

Property \textbf{P10} fails,
since $N_{DX}(z)$ is constructed with the boundary of $\TY$
and a line with a specific angle (perpendicular to the largest edge),
by Corollary \ref{cor:geo-inv-arc-prob-2}, geometry invariance for uniform data does not hold.
That is, the distributions of relative arc density and
the domination number of the corresponding PCD depend on
the geometry of the triangle $\TY$.
Hence, it does not suffice to compute them for the standard equilateral triangle only,
but instead one should compute them for each pair of $(c_1,c_2)$.
Moreover,
$\Lambda_0(N_{DX})=\{\y_1,\y_2\}$ since $N_{DX}(x,M)$ has zero area
iff $x \in \{\y_1,\y_2\}$.

\section{Discussion and Conclusions}
\label{sec:disc-conc}
In this article,
we discuss the construction of proximity catch digraphs (PCDs) in multiple dimensions.
PCDs are a special type of proximity graphs which have applications in various fields.
The class cover catch digraph (CCCD) is the first type of PCD family in literature (\cite{priebe:2001})
which is based on spherical proximity maps and has ``nice properties" for uniform data in $\R$,
in the sense that, the exact and asymptotic distribution of the domination number for CCCDs
is available for one-dimensional uniform data.
We determine some of the properties of the spherical proximity maps in $\R$
(called \emph{appealing properties}),
and use them as guidelines for extending PCDs to higher dimensions.
We also characterize the geometry invariance for PCDs based on uniform data.
Geometry invariance is important since it facilitates the computation of
quantities (such as relative arc density or domination number) related to PCDs.

We discuss four PCD families in literature and introduce two new PCD families in this article.
We investigate these PCD families in terms of the appealing properties and in particular
geometry invariance for uniform data.
We provide the asymptotic distribution of relative arc density and domination number
for some of the PCD families.
These tools have applications in spatial point pattern analysis and statistical pattern classification.
We have demonstrated that the more the properties are satisfied,
the better the asymptotic distribution of relative density.
Furthermore, the availability of the asymptotic distribution of domination number also is
highly correlated with the number of properties satisfied.

The spherical proximity regions were defined with (open) balls only,
whereas the new proximity maps are not based on a particular geometric shape
or a functional form; that is,
the new proximity maps admit any type of region, e.g., circle
(ball), arc slice, triangle, a convex or nonconvex polygon, etc.
In this sense, the PCDs are defined in a more general setting compared to CCCD.
On the other hand, the types of PCDs we introduce in this article
are well-defined for points restricted to the convex hull of $\Y_m$, $\C_H(\Y_m)$.
Moreover, the new families of proximity maps we introduce
will yield closed regions.
Furthermore, the CCCDs based on balls use proximity regions
which are defined by the obvious metric, while the PCDs do not suggest an obvious metric.

The mechanism to define the proximity maps provided in this article
can be also used for defining new (perhaps with better properties) proximity map families.

\section*{Acknowledgments}
%I would like to thank an anonymous associate editor and two referees,
%whose constructive comments and suggestions greatly improved the presentation
%and flow of the paper.
%Most of the Monte Carlo simulations presented in this article
%were executed on the Hattusas cluster of
%Ko\c{c} University High Performance Computing Laboratory.
This research was supported by the research agency TUBITAK via the Kariyer Project \# 107T647.

%{\small
%\bibliography{References}
%\bibliographystyle{apalike}
%\bibliographystyle{plain}
%}

\end{document}